%% file: main_arxiv.tex
\documentclass[onefignum]{siamart190516}

\input{ex_shared}

\ifpdf
\hypersetup{
  pdftitle={Normal Map-based Proximal SGD},
  pdfauthor={J. Qiu, L. Jiang, and A. Milzarek}
}
\fi

\begin{document}

\maketitle
\begin{abstract}
The proximal stochastic gradient method ($\PSGD$) is one of the state-of-the-art approaches for stochastic composite-type problems. In contrast to its deterministic counterpart, $\PSGD$ has been found to have difficulties with the correct identification of underlying substructures (such as supports, low rank patterns, or active constraints) and it does not possess a finite-time manifold identification property. Existing solutions rely on convexity assumptions or on the additional usage of variance reduction techniques. In this paper, we address these limitations and present a simple variant of $\PSGD$ based on Robinson's normal map. The proposed normal map-based proximal stochastic gradient method ($\NSGD$) is shown to converge globally, i.e., accumulation points of the generated iterates correspond to stationary points almost surely. In addition, we establish complexity bounds for $\NSGD$ that match the known results for $\PSGD$ and we prove that $\NSGD$ can almost surely identify active manifolds in finite-time in a general nonconvex setting. Our derivations are built on almost sure iterate convergence guarantees and utilize  analysis techniques based on the Kurdyka-{\L}ojasiewicz inequality. 
\end{abstract}



\section{Introduction}
In this work, we propose and investigate a novel normal map-based proximal stochastic gradient method for the composite problem
\begin{equation}
	\label{eq:SO} \min_{x\in \Rn}~\psi(x) := f(x) + \vp(x),
\end{equation}
where $\vp: \Rn \to \Rex$ is a convex, lower semicontinuous (lsc.), and proper mapping and $f: \Rn \to \R$ is continuously differentiable on an open set containing the effective domain $\dom{\vp} := \{x \in \Rn: \vp(x) < \infty\}$. 

The model \cref{eq:SO} has gained increasing attention during the recent decades and is frequently used in large-scale applications and stochastic optimization tasks, including, e.g.,  machine learning \cite{Bis06,BotCurNoc18,LeCun2015}, statistical learning and sparse regression \cite{friedman2001elements,SST2011}, and stochastic programming \cite{nemjudlansha08,ShaDenRus14}. 
In these applications, the (potentially nonsmooth) function $\vp$ is typically chosen to promote specific structural properties, such as sparsity, group sparsity, low rank, etc., while the smooth (not necessarily convex) part $f$ corresponds to a data-driven learning model or loss function.  
Since a full evaluation of the function and gradient values $f(x)$ and $\nabla f(x)$ can be prohibitively expensive, sampling schemes or stochastic approximation techniques are often employed in practice to generate more tractable information. 
This integral mechanism, pioneered by Robbins and Monro \cite{robbins1951stochastic}, is the basis of the stochastic gradient descent method ($\SGD$) and many other successful stochastic algorithms.  

Proximal stochastic gradient descent methods, \cite{duchi2009efficient,xiao2014proximal,nitanda2014stochastic,GhaLanZha16,AtcForMou17,MajMiaMou18,rosasco2020convergence}, extend the basic stochastic approximation principles utilized in $\SGD$ to the composite problem \cref{eq:SO}. 
The update rule of the basic proximal stochastic gradient method ($\PSGD$) for \cref{eq:SO} is given by 
\be \label{eq:update-PSGD}
 x^{k+1} = \prox{\alpha_k \vp}(x^k - \alpha_k g^k),
\ee
where the sequence $\revise{\{\alpha_k\}_{k\in\N}}\subseteq \R_{++}$ denotes the selected step sizes, $g^k$ is a (stochastic) approximation of $\nabla f(x^k)$, and the well-known proximity operator of $\vp$ is defined as $\prox{\alpha \vp}(x) := {\argmin}_{y\in\Rn} \, \vp(y) + \frac{1}{2\alpha}\|x-y\|^2$, $\alpha >0$. \vspace{1mm}   

Throughout this work, we consider the following base set of conditions for \cref{eq:SO}:
\vspace{1.5mm}
    \begin{enumerate}[label=\textup{\textrm{(H.\arabic*)}},leftmargin=3em,topsep=0ex,itemsep=1ex,partopsep=0ex]
        \item \label{A1} The gradient $\nabla f : \Rn \to \Rn$ is Lipschitz continuous on $\dom{\vp}$ with modulus $\sL>0$ and $\vp$ is convex, lsc., and proper.
        \item \label{A2} The objective function $\psi$ is bounded from below, i.e., $\bar \psi := \inf_{x \in \Rn}\,\psi(x) > -\infty$. 
    \end{enumerate}
    \vspace{1.5mm}

In addition, we assume that there is a filtered probability space $(\Omega,\mathcal F,\revise{\{\mathcal F_k\}_{k\in\N}},\Prob)$ that can describe the stochastic components of the update \cref{eq:update-PSGD} and of other stochastic algorithms studied in this work in a unified way. In particular, each approximation $g^k$ is understood as a realization of an $\mathcal F_{k+1}$-measurable random vector $\sg^k : \Omega \to \Rn$. Thus, the stochastic process $\revise{\{\sx^k\}_{k\in\N}}$, associated with \cref{eq:update-PSGD} (or other stochastic methods), is adapted to the filtration $\revise{\{\mathcal F_k\}_{k\in\N}}$. We further define the stochastic error terms $\se^k : \Omega \to \Rn$ via $\se^k := \sg^k - \nabla f(\sx^k)$ and set $\Exp_k[\cdot] := \Exp[\cdot | \cF_k]$. We mainly work with the following conditions on the errors $\revise{\{\se^k\}_{k\in\N}}$ and step sizes $\revise{\{\alpha_k\}_{k\in\N}}$:
\vspace{1.5mm}
\begin{enumerate}[label=\textup{\textrm{(S.\arabic*)}},leftmargin=3em,topsep=0ex,itemsep=1ex,partopsep=0ex]
        \item \label{B1} The stochastic error $\se^k$ has mean zero, i.e., $\Exp_k[\se^k]=0$ almost surely (a.s.).
        \item \label{B2} There exists $\revise{\{\sigma_k\}_{k\in\N}}\subseteq \R_+$ such that $\Exp_k[\|\se^k\|^2]\leq \sigma_k^2$ a.s. for all $k$. 
        \item \label{B3} The step sizes $\revise{\{\alpha_k\}_{k\in\N}}$ satisfy $\alpha_k\to0$, $\sum_{k=0}^\infty \alpha_k = \infty$, and $\sum_{k=0}^{\infty}\alpha_k^2\sigma_k^2<\infty$.
    \end{enumerate}
    \vspace{1.5mm}
The assumptions \ref{A1}--\ref{A2} and \ref{B1}--\ref{B3} are standard in the analysis of stochastic optimization methods, see, e.g., \cite{GhaLanZha16,duchiruan2018,davdru19,davis2020stochastic}.

\begin{figure}
	\setlength{\abovecaptionskip}{-3pt plus 3pt minus 0pt}
	\setlength{\belowcaptionskip}{-10pt plus 3pt minus 0pt}
    \centering
\begin{tikzpicture}[scale=1]
    \node[right] at (2.97,2.45) {\footnotesize $\NSGD$};
    \node[right] at (6.935,2.45) {\footnotesize $\PSGD$};
    \node[right] at (10.9,2.45) {\footnotesize $k \mapsto \frac3k$};
    \draw [draw=MyGray] (0.89,2.2) rectangle (12.78,2.7);
    \draw [draw=DeepBlue,fill=DeepBlue,line width = 1pt] (1.87,2.45) -- (2.77,2.45);
    \draw [draw=MyGrayTwo,fill=MyGrayTwo,line width = 1pt] (5.835,2.45) -- (6.735,2.45);
    \draw [draw=MyOrange,fill=MyOrange,line width = 1pt] (9.8,2.45) -- (10.1,2.45);
    \draw [draw=MyOrange,fill=MyOrange,line width = 1pt] (10.2,2.45) -- (10.3,2.45);
    \draw [draw=MyOrange,fill=MyOrange,line width = 1pt] (10.4,2.45) -- (10.7,2.45);
	\node[right] at (0,0) {\includegraphics[width=0.48\linewidth]{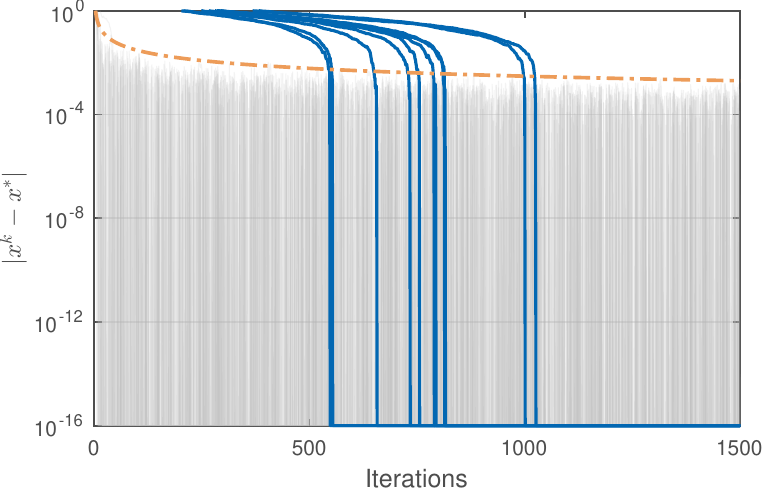}};     
    \node[right] at (6.5,0) {\includegraphics[width=0.48\linewidth]{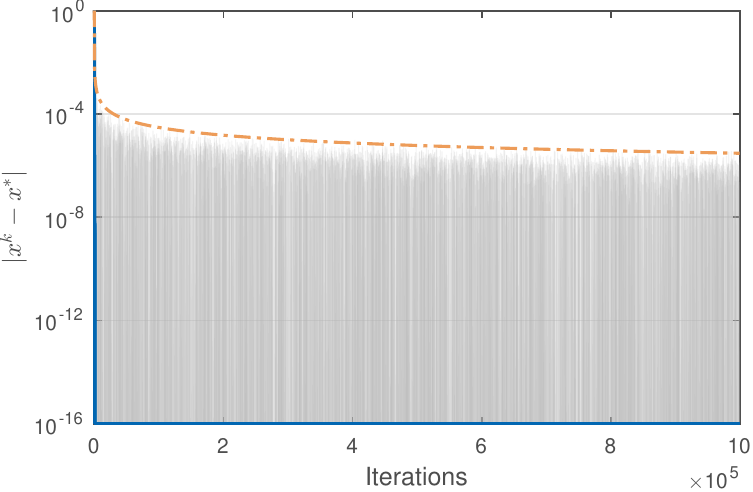}}; 
\end{tikzpicture}    
\caption{Behavior of $\PSGD$ and $\NSGD$ on ${\min}_{x \in \R}~f(x)+ \vp(x)$ where $f(x) := x$ and $\vp(x) := \mathds{1}_{[-1,1]}(x)$. This example originates from \cite{duchi2021asymptotic}. $\PSGD$ and $\NSGD$ use the stochastic gradients $g^k = f^\prime(x^k) + e^k = 1+e^k$ with iid Gaussian noise $e^k\sim \mathcal{N}(0,1)$ and $\alpha_k := \frac{1}{k}$, $\lambda = 1$. In this scenario, both $\PSGD$ and $\NSGD$ converge to the global solution $x^* = -1$ a$.$s$.$. 
The plots depict the distance $|x^k-x^*|$ (in logarithmic scale) for $10$ independent runs with $x^0 = z^0 = 100$. \revise{The plot on the left is a zoomed-in view for the first $1500$ iterations.} The dash-dotted line, $k\mapsto\frac{3}{k}$, approximates the convergence trend of $\PSGD$. Notably, $\PSGD$ repeatedly escapes from $x^*$ (i.e., from the active manifold $\mathcal M_{x^*} := \{x^*\}$), whereas $\NSGD$ remains at $x^*$ after identifying the active constraint. Indeed, in \cite[p$.$ 30]{duchi2021asymptotic}, it is shown that the stochastic process $\revise{\{\sx^k\}_{k\in\N}}$ generated by $\PSGD$ satisfies $\Prob(\sx^{k} \geq -1 + \alpha_k) \geq \varepsilon > 0$. Thus, $\revise{\{\sx^k\}_{k\in\N}}$ will not stay on $\mathcal M_{x^*}$ a$.$s$.$ even if the optimal solution is correctly identified.}
\label{fig:toy1}
\end{figure}

Despite its popularity and convincing theoretical properties  \cite{GhaLanZha16,duchiruan2018,davdru19,rosasco2020convergence,LiMil22,LiaXu24}, Xiao \cite{xiao10a} observed limitations of $\PSGD$ when solving sparse problems with $\vp(x) = \nu \|x\|_1$, $\nu>0$. Specifically, the solutions produced by $\PSGD$ exhibit lower levels of sparsity compared to those obtained from the deterministic proximal gradient descent method ($\PGD$). This behavior is well recognized and not just limited to $\PSGD$. In \cite{duchi2021asymptotic}, Duchi and Ruan provide an example showing that the projected $\SGD$ method fails to identify the underlying active constraints (cf$.$ \Cref{fig:toy1}). In \cite{pooliasch18}, Poon et al$.$ study this effect from the perspective of manifold identification \cite{lewis2002active,lewis2013partial,lewis2016proximal}. As in \cite{xiao10a,duchi2021asymptotic}, their analysis reveals that $\PSGD$ generally does not have an identification property, which can explain this undesirable behavior. Several variants of $\PSGD$ have been shown to possess a manifold identification property, albeit under specific conditions such as (strong) convexity and requiring variance reduction techniques. 

In \cite{lee2012manifold}, Lee and Wright prove that the regularized dual averaging method ($\RDA$) \cite{xiao10a} can identify optimal manifolds with high probability. The authors consider the model \cref{eq:SO} with $f$ having the form $f(x) = \Exp[F(x,\sxi)]$, where $F : \Rn \times \Xi \to \R$ is a loss function and $\sxi : \Omega \to \Xi$ is a random variable. The mappings $x \mapsto F(x,\xi)$ and $x \mapsto \nabla F(x,\xi)$ are assumed to be convex and $\sL$-Lipschitz continuous for all $\xi$, respectively. In addition, the gradient $\nabla F(x,\xi)$ is supposed to be uniformly bounded for all $x$ and $\xi$. To ensure identification, Lee and Wright further assume the existence of a locally strong (hence unique) minimizer $x^*$ of \cref{eq:SO} such that the strict complementarity condition, $0 \in \mathrm{ri}(\partial\psi(x^*))$, holds and $\psi$ is partly smooth at $x^*$ relative to the optimal manifold $\mathcal M_{x^*}$ in the sense of Lewis \cite{lewis2002active} (cf$.$ \Cref{sec:manifold}). 
In \cite{duchi2021asymptotic}, Duchi and Ruan establish additional results for $\RDA$ in a similar setting,  when $\vp = \mathds{1}_{\mathcal C}$ is the indicator function of a convex set $\mathcal C := \{x: f_i(x) \leq 0, i= 1,\dots,m\}$. Under a bounded variance-type assumption and a quadratic growth condition, the authors prove a$.$s$.$ convergence of the iterates to the unique solution $x^*$. Moreover, under the LICQ, a second-order sufficient condition, and the strict complementarity condition, it is shown that $\RDA$ a$.$s$.$ identifies the active constraints in finite-time. 
Huang and Lee propose a variant of $\RDA$ with momentum \cite{huang2022training}. The identification results in \cite{huang2022training} are applicable to general nonconvex problems. However, the iterates $\revise{\{\sx^k\}_{k\in\N}}$ are assumed to converge a$.$s$.$ and the difference of consecutive iterates is required to converge to zero a$.$s$.$ at a sufficiently fast rate. 
In \cite{pooliasch18}, Poon et al$.$ provide a framework for the stochastic method \cref{eq:update-PSGD} to ensure identification. According to \cite{pooliasch18}, the scheme \cref{eq:update-PSGD} has a manifold identification property if: (i) $\liminf_{k\to\infty} \alpha_k > 0$; (ii) $\se^k \to 0$ a$.$s$.$; (iii) there is $x^* \in \argmin_{x \in \Rn}\psi(x)$ such that $\sx^k \to x^*$ a$.$s$.$; (iv) $\psi$ is partly smooth at $x^*$ (relative to the optimal manifold) and $0 \in \mathrm{ri}(\partial\psi(x^*))$. The authors apply this framework to $\PSVRG$ and $\PSAGA$ and verify the conditions (i)--(iii) when $f$ has the finite-sum form $f(x) = \frac1N\sum_{i=1}^N f_i(x)$ and each $f_i$ is convex. Sun et al$.$ \cite{sun2019manifold} establish identification results for $\PSVRG$, $\PSAGA$, and $\RDA$ that hold in expectation. It is assumed that $f$ is a convex finite-sum term, $\nabla f$ is Lipschitz, $\vp$ is separable, problem \cref{eq:SO} has a unique solution $x^*$, and a strict complementarity-type condition has to hold.  Depending on $\Exp[\|\sx^k-x^*\|]$, $\Exp[\|\se^k\|]$, and other terms, bounds are derived for the iteration index $k$ after which the optimal manifold is identified.
Dai et al$.$ \cite{dai2023} study identification properties of a proximal variant of $\STORM$, \cite{cutkosky2019storm}. High probability finite-time identification results are shown if $f$ is strongly convex, $\vp$ is fixed to a non-overlapping group regularization, the generated momentum directions are uniformly bounded a$.$s$.$, and the strict complementarity condition is satisfied.
In the recent preprint \cite{qin2025partial}, identifications properties are studied in a relaxed setting and without requiring the strict complementarity condition $0 \in \mathrm{ri}(\partial\psi(x^*))$. 
The authors establish identification results for a mini batch-variant of $\PSGD$ applied to LASSO when the mini-batch size is sufficiently large and $\|\sx^{k+1}-\sx^k\|/{\alpha_k} \to 0$.

In this work, our aim is to develop a simple stochastic ($\PSGD$-type) algorithm for the nonconvex composite problem \cref{eq:SO} that has comprehensive convergence properties, including a$.$s$.$ finite-time manifold identification guarantees, without relying on common variance reduction techniques. An overview of the proposed normal map-based method is shown below (see \Cref{subsec:design} for further details). \vspace{2mm}

\begin{mdframed}[
  linecolor=white!100,
  linewidth=.2mm,          
  backgroundcolor=gray!10, 
  roundcorner=0pt,         
  innerleftmargin=10pt,    
  innerrightmargin=15pt,   
  innertopmargin=5pt,     
  innerbottommargin=5pt   
]
\textbf{A normal map-based proximal stochastic gradient method ($\NSGD$)} \\[2mm]
Let $\lambda>0$, $\revise{\{\alpha_k\}_{k\in\N}}\subseteq \R_{++}$ be given. Choose $z^0\in\Rn$ and set $x^0=\proxl(z^0)$. \\[1mm]
\textbf{For} $k = 0,1,\dots$ \textbf{do:} 
\vspace{-1mm}
\[ z^{k+1} = z^k - \alpha_k(g^k + \lambda^{-1}(z^k-x^k)) \quad \text{and} \quad
    x^{k+1} = \proxl(z^{k+1}).\]
\end{mdframed}
\vspace{2mm}
In each iteration, $\NSGD$ generates one stochastic gradient $g^k \approx \nabla f(x^k)$ and performs one proximal step. Thus, the computational costs of $\NSGD$ and $\PSGD$ are essentially identical. We now summarize some of the fundamental properties of $\NSGD$ and our core contributions: 

\begin{itemize}[leftmargin = 15pt]
    \item \emph{(Stationarity)} The design of $\NSGD$ is based on Robinson's \emph{normal map} \cite{robinson1992normal}: 
    \[\Fnor(z):=\nabla f(x) + \lambda^{-1}(z-x)\quad  \text{where} \quad x=\proxl(z), \quad \lambda > 0.\]
    The normal map is connected to the standard stationarity measure $\dist(0,\partial\psi(x))$ in the sense that the condition $\|\Fnor(z)\| \leq \varepsilon$ implies $\dist(0,\partial\psi(\proxl(z))) \leq \varepsilon$ (cf$.$ \Cref{sec:stationary}). This fact plays a key role in the identification results and allows us to establish stronger stationarity guarantees for $\NSGD$. 
    \item \emph{(Complexity)} We derive complexity bounds for $\NSGD$ in expectation and in terms of the normal map that align with the existing results for $\PSGD$ \cite{davdru19}. 
    \item \emph{(Global convergence)} Under \ref{A1}--\ref{A2} and \ref{B1}--\ref{B3}, we show that accumulation points of the stochastic process $\revise{\{\sx^k\}_{k\in\N}}$ generated by $\NSGD$ are stationary points of $\psi$ a$.$s$.$. This allows us to recover existing global convergence results for $\PSGD$ \cite{duchiruan2018,MajMiaMou18,davis2020stochastic,LiMil22} under weaker assumptions.
    \item \emph{(Iterate convergence and identification)} We prove a$.$s$.$ iterate convergence, $\sx^k \to \sx^*$, when $\psi$ is additionally assumed to be definable \revise{\cite{coste2000introduction,van1998tame}} and if the event $\{\omega \in \Omega: \liminf_{k\to\infty}\|\sx^k(\omega)\| < \infty\}$ occurs with probability $1$. To our knowledge, such general convergence guarantees seem to be new for proximal-type stochastic algorithms. Based on these results, we further verify that $\NSGD$ can identify the underlying active manifold structure, \cite{lewis2002active,lewis2013partial,liang2017activity}, a$.$s$.$ in finite time.
\end{itemize}

\subsection{\texorpdfstring{Designing $\NSGD$}{Designing Norm-SGD}}\label{subsec:design}
To motivate $\NSGD$, we first revisit the deterministic proximal gradient method and express it in a form that separates the gradient and the proximal step. Specifically, the basic update of $\PGD$ can be written as   
\begin{equation}\label{eq:equiv-form}
    \left[ 
    \begin{aligned}
        z^{k+1} &= x^{k} - \lambda \nabla f(x^{k})\\
        x^{k+1} & = \proxl(z^{k+1})
    \end{aligned} 
    \right. \quad \Longleftrightarrow \quad 
    \left[\begin{aligned}
        z^{k+1} &= z^{k} - \textcolor{red}{\alpha} [\nabla f(x^{k}) + \textcolor{blue}{\lambda}^{-1}(z^{k} - x^k)]\\
        x^{k+1} & = \prox{\textcolor{blue}{\lambda}\vp}(z^{k+1}), \quad \textcolor{red}{\alpha} = \textcolor{blue}{\lambda}.
    \end{aligned}\right. 
\end{equation}
This equivalent formulation naturally introduces the normal map $\Fnor$, which has been widely used in the field of variational inequalities and generalized equations, \cite{facchinei2007finite}. 

Based on \cref{eq:equiv-form}, we can decouple the proximal parameter \textcolor{blue}{$\lambda$} and the step size $\textcolor{red}{\alpha} = \textcolor{blue}{\lambda}$. In particular, we propose to replace the constant step size $\textcolor{red}{\alpha} = \textcolor{blue}{\lambda}$ by a sequence of varying step sizes $\revise{\{\alpha_k\}_{k\in\N}}$ without adapting the proximal parameter \textcolor{blue}{$\lambda$}. This leads to the following proximal stochastic update that defines $\NSGD$: 
\be\label{eq:update-short}
    \sz^{k+1} = \sz^k - \alpha_k(\sg^k + \lambda^{-1}(\sz^k-\sx^k)) \quad \text{and} \quad \sx^{k+1} = \proxl(\sz^{k+1}),
\ee
where $\sg^k$ is a (unbiased) stochastic approximation of $\nabla f(\sx^k)$.

Our convergence analysis of $\NSGD$ is based on and focuses on the auxiliary iterates $\revise{\{\sz^k\}_{k\in\N}}$. \revise{The update rule of $\NSGD$ can be interpreted as a fixed-point iteration in the sense of the stochastic Krasnoselskii--Mann method ($\SKM$), \cite{ComPes15,bravo2024stochastic}. Specifically, introducing the operator $T(z):= \proxl(z) - \lambda \nabla f(\proxl(z))$, the iteration \cref{eq:update-short} takes the form
\begin{equation} \label{eq:fixed-point}
\sz^{k+1} = ( 1 - {\alpha_k}{\lambda^{-1}}) \sz^k + {\alpha_k}{\lambda^{-1}} \cdot (T(\sz^k) - \lambda \se^k), \end{equation}
where $\se^k = \sg^k - \nabla f(\sx^k)$, $\sx^k = \proxl(\sz^k)$, is exactly the stochastic error. $\NSGD$ differs from classical proximal stochastic methods \cite{GhaLanZha16,davdru19}, which typically focus on the main iterates $\{\sx^k\}_{k\in\N}$ and cannot be directly or easily recast as fixed-point iterations involving a fixed operator $T$ (which is independent of the iteration $k$) and mean-zero error terms. In $\PSGD$, this difficulty is primarily caused by the proximity operator $\mathrm{prox}_{\alpha_k\vp}$ and its dependence on the diminishing step sizes ${\{\alpha_k\}_{k\in\N}}$. In particular and in addition, we  generally have $\Exp_k[\prox{\alpha_k \vp}(\sx^k - \alpha_k \sg^k)] \neq \prox{\alpha_k \vp}(\sx^k - \alpha_k \nabla f(\sx^k))$. By contrast, in $\mathsf{Norm}$-$\mathsf{SGD}$}, we can leverage the  unbiasedness of the stochastic gradients $\sg^k$, i.e., taking conditional expectation in \cref{eq:update-short}, it follows
\[\Exp_k[\sz^{k+1}] =  \Exp_k[\sz^k - \alpha_k (\sg^k + \lambda^{-1}(\sz^k-\sx^k)) ] = \sz^k - \alpha_k \Fnor(\sz^k). \]
This property will play a distinctive and important role in our derivations.

\revise{As noted earlier, the regularized dual averaging method ($\RDA$) also has comprehensive identification guarantees in the convex setting. We provide a more detailed and formal comparison of $\NSGD$ and $\RDA$ in \Cref{app:comparison}.}

\subsection{Stationarity measures}\label{sec:stationary}
We now briefly discuss common stationarity measures for proximal gradient-type methods and for the problem \cref{eq:SO}. The first-order optimality conditions associated with \cref{eq:SO} are given by $0 \in \partial \psi(x) = \nabla f(x) + \partial \vp(x)$. Here, $\partial \psi$ denotes the limiting subdifferential of $\psi$ and $\partial \vp$ is the standard subdifferential for convex functions, cf. \cite{rocwet98,baucom11}. A point $\bar x \in \dom{\vp}$ is called a stationary point of \cref{eq:SO} if $0 \in \partial \psi(\bar x)$ and we use $\crit(\psi)$ to denote the set of all stationary points of $\psi$. We further introduce the notation $\|\partial \psi(x)\|_{-} := \dist(0,\partial\psi(x)) := {\min}_{v\in\partial \psi(x)} \;\|v\|$. A point $\bar x \in \dom{\vp}$ is an $\varepsilon$-stationary point of $\psi$, $\varepsilon \geq 0$, if we have $\|\partial\psi(\bar x)\|_- \leq \varepsilon$. 

\emph{Natural residual.} The first-order optimality conditions for problem \cref{eq:SO} can be equivalently expressed using the so-called \emph{natural residual}; it holds that:
\[ x \in \crit(\psi) \quad \iff \quad \Fnat(x) := \lambda^{-1}(x-\proxl(x-\lambda \nabla f(x))) = 0, \quad \lambda > 0. \]
The natural residual $\Fnat$ serves as an important stationarity measure and is a basic tool in the design of algorithms for composite problems, \cite{baucom11,Bec17}. Moreover, based on the characterization $p = \proxl(y) \iff  p \in y - \lambda\partial\vp(p)$, $y \in \Rn$, it follows that
\[
    \Fnat(x) =\lambda^{-1}(x-x^+)\in \nabla f(x) + \partial \vp(x^+), \quad \text{where} \quad x^+:=\proxl(x-\lambda \nabla f(x)). \]
 
\begin{figure}[t]
 \setlength{\belowcaptionskip}{-10pt plus 3pt minus 0pt}
       \centering
\includegraphics[height=4.4cm]{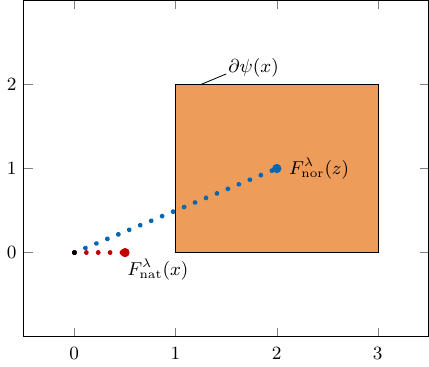} \hspace{3.5ex}
\includegraphics[height=4.4cm]{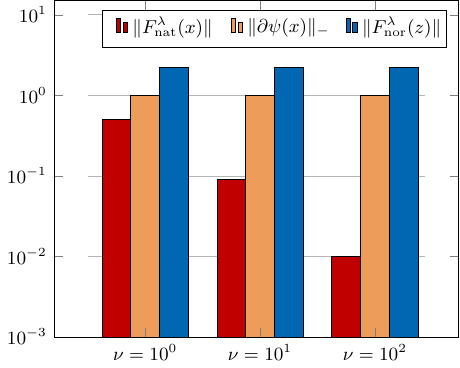}

\caption{Relation between different stationary measures. We consider ${\min}_{x \in \R^2}~f(x)+ \vp(x)$, where $f(x) := \frac12\|x+[2,1]^\top\|^2$ and $\vp(x) := \nu \|x\|^2 + \|x\|_1$, $\nu>0$. We further set $\lambda = \frac12$, $z=0$, and $x=\proxl(z) = 0$.  Left: Simple calculations yield $\Fnat(x) = [\frac{1}{1+\nu},0]^\top$, $\nu = 1$, $\Fnor(z) = [2,1]^\top$, and $\partial \psi(x)= [1,3]\times[0,2]$. 
Clearly, it holds that $\Fnor(z) \in \partial \psi(x)$ but $\Fnat(x) \notin \partial \psi(x)$. Right: Comparison of $\|\Fnat(x)\|$, $\|\Fnor(z)\|$, and $\|\partial\psi(x)\|_{-}$ for different regularization parameters $\nu$. 
       }
       \label{fig2:example}
\end{figure}

\emph{Normal map.} Robinson's normal map is an alternative stationarity measure for \cref{eq:SO}. It can be interpreted as a special subgradient of $\psi$. In particular, we have  
\begin{equation} \label{eq:norm-subdiff} \Fnor(z) = \nabla f(\proxl(z)) + \lambda^{-1}(z-\proxl(z)) \in \partial\psi(\proxl(z)) \quad \forall~z\in\Rn. \end{equation}
Thus, if $\bar z$ is a zero of $\Fnor$, then, by \cref{eq:norm-subdiff}, $\bar x = \proxl(\bar z)$ is a stationary point of \cref{eq:SO}. Conversely, if $\bar x$ satisfies $\Fnat(\bar x) = 0$, then $\bar z = \bar x-\lambda\nabla f(\bar x)$ is a zero of the normal map, see \cite{ouyang2021trust}. Applying \cite[Lemma 4.1.6]{milzarek2016numerical} or \cite[Theorem 3.5]{drusvyatskiy2018error}, it holds that
\begin{equation} \label{eq:stat-am} \|\Fnat(x)\| \leq \|\partial\psi(x)\|_- \quad \text{and} \quad \|\partial\psi(\prox{\lambda\vp}(z))\|_- \leq \|\Fnor(z)\| \end{equation}
for all $x \in \dom{\partial\vp}$ and $z \in \Rn$. Hence, the condition $\|\Fnor(z)\| \leq \varepsilon$ implies that $x = \proxl(z)$ is an $\varepsilon$-stationary point. Clearly, if $x$ is a point satisfying $\|\Fnat(x)\| \leq \varepsilon$, then it does not necessarily need to be an $\varepsilon$-stationary point of $\psi$. We illustrate this potential mismatch between $\|\Fnat(x)\|$ and $\|\partial\psi(x)\|_-$ in \Cref{fig2:example} on a toy example.

\emph{Moreau envelope}. The gradient of the Moreau envelope, $\nabla \env{\theta \psi}(x):=\theta^{-1}(x-\hat x)$, $\hat x:=\prox{\theta \psi}(x)$, has received increasing attention as a near-stationarity measure for nonsmooth problems, \cite{drusvyatskiy2019efficiency,davdru19,bot2023alternating}. \revise{Here, the Moreau envelope $\env{\theta \psi}$ is defined via $\env{\theta \psi}(x) := \min_{y \in \Rn} \psi(y) + \frac1{2\theta}\|x-y\|^2$}. 
In the composite setting \cref{eq:SO}, it is well known that $\|\nabla \env{\theta \psi}(x)\|$ is proportional to $\|\Fnat(x)\|$ for appropriate $\theta, \lambda$, \cite{drusvyatskiy2018error,davdru19}. Moreover, this measure reflects the near-stationarity of the proximal point $\hat x$, i.e., the condition $\|\nabla \env{\theta \psi}(x)\| \leq \varepsilon$ implies $\|\partial \psi(\hat x)\|_{-}\leq\varepsilon$ rather than  $\|\partial \psi(x)\|_{-}\leq\varepsilon$.
\subsection{\texorpdfstring{Further background on $\PSGD$}{Additional literature on Prox-SGD}}\label{subsec:rel-work}
While the convergence of $\PSGD$, \cref{eq:update-PSGD}, is well-understood in the strongly convex case, general results for stochastic and nonconvex composite  problems seem to be more limited. In \cite{davdru19}, under \ref{A1}--\ref{A2} and \ref{B1}--\ref{B2}, Davis and Drusvyatskiy present one of the first complexity results for $\PSGD$ in the nonconvex setting using the Moreau envelope as a merit function:
\begin{equation} \label{eq:comp-proxsgd} {\min}_{i=0,\dots,k-1} \ \Exp[\|\nabla\env{\theta\psi}(\sx^k)\|^2] = \mathcal O(\textstyle{\sum}_{i=0}^{k-1}\alpha_i^2 \,/\, \textstyle{\sum}_{i=0}^{k-1}\alpha_i). \end{equation}
Earlier studies of $\PSGD$ for nonconvex $f$ appeared in \cite{GhaLanZha16}, where convergence of $\PSGD$ is shown if the variance $\Exp[\|\se^k\|^2] \to 0$ vanishes as $k \to \infty$ which can be ensured via suitable variance reduction techniques, \cite{xiao2014proximal,j2016proximal}. 
It is not fully known whether $\PSGD$ converges globally under the base conditions \ref{A1}--\ref{A2} and \ref{B1}--\ref{B3}. Following the classical ODE-based analysis of $\SGD$, Majewski et al$.$ use a differential inclusion approach to study convergence of $\PSGD$, \cite{MajMiaMou18}. The authors introduce additional compact constraints to ensure applicability of the differential inclusion techniques. The analyses in \cite{duchiruan2018,davis2020stochastic} establish convergence for a broader class of proximal algorithms including $\PSGD$ as a special case. Both works apply differential inclusion-based mechanisms and require (almost) sure boundedness of $\revise{\{\sx^k\}_{k\in\N}}$ and a density-type condition to show that accumulation points of $\revise{\{\sx^k\}_{k\in\N}}$ correspond to stationary points. More recently, Li and Milzarek \cite{LiMil22} derive asymptotic convergence results for $\PSGD$ of the form $\Exp[\|\Fnat(\sx^k)\|] \to 0$ and $\|\Fnat(\sx^k)\| \to 0$ a$.$s$.$ under the additional assumption that $\varphi$ is globally Lipschitz continuous. It is not known whether this Lipschitz assumption can be relaxed or dropped. The easier convex and strongly convex cases have been investigated, e.g., in \cite{GhaLanZha16,AtcForMou17,rosasco2020convergence,patrascu2021stochastic,LiaXu24}. If $f$ is convex and $\psi$ is strongly convex, then convergence in expectation, $\Exp[\|\sx^k - \sx^*\|^2] = \mathcal O(k^{-1})$, can be ensured, see \cite{rosasco2020convergence,patrascu2021stochastic}. Moreover, when both $f$ and $\vp$ are convex, Rosasco et al$.$ \cite{rosasco2020convergence} show that there is a subsequence of $\revise{\{\sx^k\}_{k\in\N}}$ converging to $\sx^*$ a$.$s$.$. 

\section{Preliminaries and basic properties} 
In this section, we derive an approximate descent property for $\NSGD$, which will be one of our main tools in the convergence analysis. In \Cref{subsec:computation and lemma}, we state preliminary estimates for the distance between iterates and for certain accumulated stochastic error terms. In \Cref{subsec:normal map bound and descent}, we provide two basic bounds that involve the normal map and the function values $\revise{\{\psi(\sx^k)\}_{k\in\N}}$. Based on these results and using a ``normal map-friendly'' merit function for \cref{eq:SO}, we then establish approximate descent of $\NSGD$ in \Cref{subsec:merit function}. 

\subsection{Preparatory lemmas}\label{subsec:computation and lemma}
We first present several basic estimates and properties that facilitate the analysis of $\NSGD$. Since $\vp$ is assumed to be convex, lsc., and proper (cf$.$ \ref{A1}), the proximity operator $\proxl$ is firmly nonexpansive, \cite{mor65,baucom11}. Hence, recalling $\sx^k = \proxl(\sz^k)$, it holds that
\begin{equation}\label{prox-prop}
    \|\sx^n-\sx^m\|^2 \leq \iprod{\sz^n-\sz^m}{\sx^n-\sx^m}\leq \norm{\sz^n-\sz^m}^2 \quad \forall~\revise{m<n} \quad \text{(a$.$s$.$)}.
\end{equation}
In addition, by the update rule \cref{eq:update-short}, we can decompose $\sz^n$ as follows:
\begin{align}
    \nonumber \sz^n&=\sz^{n-1} -\alpha_{n-1}[\sg^{n-1}+ \tfrac{1}{\lambda}(\sz^{n-1} - \sx^{n-1})] =\sz^{n-1} - \alpha_{n-1}[\Fnor(\sz^{n-1}) +\se^{n-1}]	\\
     \nonumber &=\sz^m - {\sum}_{i=m}^{n-1}\alpha_i\Fnor(\sz^i) - {\sum}_{i=m}^{n-1} \alpha_i\se^i \\
    \nonumber &=\sz^m - {{\sum}_{i=m}^{n-1}\alpha_i}\Fnor(\sz^m) - {\sum}_{i=m}^{n-1} \alpha_i[\Fnor(\sz^i)-\Fnor(\sz^m)] - {\sum}_{i=m}^{n-1} \alpha_i\se^i \\ & = \sz^m - \tau_{m,n} \Fnor(\sz^m) + \scp^{m,n} \revise{\quad \forall~m<n,}\label{eq:pretty-important}
\end{align}
where $\tau_{m,n}:= {\sum}_{i=m}^{n-1}\alpha_i$ and $\scp^{m,n}:=-{\sum}_{i=m}^{n-1} \alpha_i[\Fnor(\sz^i)-\Fnor(\sz^m)] - {\sum}_{i=m}^{n-1} \alpha_i\se^i$.
This representation will play a key role in our derivations. Let us further introduce the stochastic aggregated error term
\[
\scs_{m,n} :=\max_{m<j\leq n}\Big\|{\sum}_{i=m}^{j-1} \alpha_i \se^i\Big\|.
\]
In the following, we state upper bounds for the iterates and $\scp^{m,n}$. The verification of \Cref{lem:est-err} is deferred to \cref{proof:lem:est-err}.
\begin{lemma}[Iterate bounds]
    \label{lem:est-err}
     Let \ref{A1}--\ref{A2} and \ref{B1}--\ref{B2} hold and let $\revise{\{\sx^k\}_{k\in\N}}$ and $\revise{\{\sz^k\}_{k\in\N}}$ be generated by $\NSGD$. Let us set $\Llambda := \sL+2\lambda^{-1}$. For all $0\leq m<n$ with $\tau_{m,n} \leq 1/(2\Llambda)$, we have
    \begin{align*}
        \max_{m<i\leq n}\|\sx^i-\sx^m\| &\leq \max_{m<i\leq n}\|\sz^i-\sz^m\| \leq 2\tau_{m,n} \norm{\Fnor(\sz^m)} + 2\scs_{m,n}\quad \text{and}\\
        \|{\scp^{m,n}}\| &\leq 2\Llambda\tau_{m,n}^2\norm{\Fnor(\sz^m)} + 2\scs_{m,n}. 
    \end{align*}
\end{lemma}

Next, we provide a standard martingale-type estimate for the aggregated noise terms $\scs_{m,n}$. The proof of \cref{lem:err_estimate} is based on a routine application of the Burk- holder-Davis-Gundy inequality and is presented in \cref{app:proof-martingale}. We refer to \cite[Proposition 4.2]{benaim1999dynamics} and \cite[Lemma 6.1]{tadic2015convergence} for related results. 

\begin{lemma}
    \label{lem:err_estimate}
    Let \ref{B1}--\ref{B2} hold and let $\revise{\{\alpha_k\}_{k\in\N}}$ be given. Then for any (arbitrary) subsequence $\revise{\{t_k\}_{k\in\N}} \subseteq \N$, we have $\Exp[{{\scs_{t_k,t_{k+1}}^2}}] \leq 4 \sum_{i=t_k}^{t_{k+1}-1}\alpha_i^2\sigma_i^2$. Moreover, if there is a non-decreasing sequence $\revise{\{\beta_k\}_{k\in\N}}\subset\R_{++}$ such that ${\sum}_{k=0}^\infty \alpha_k^2\beta_k^2\sigma_k^2 < \infty$, it holds that
    \[ {\sum}_{k=0}^\infty \beta_{t_k}^2 \Exp[{{\scs_{t_k,t_{k+1}}^2}}] < \infty \quad \text{and}\quad {\sum}_{k=0}^\infty \beta_{t_k}^2 \scs_{t_k,t_{k+1}}^2 < \infty \quad \text{a$.$s$.$.}\] 
\end{lemma}

\subsection{Bounding the normal map and function values}\label{subsec:normal map bound and descent}
We now present basic bounds for the objective function values $\psi(\sx^n)$ and  for $\|\Fnor(\sz^n)\|$. Before stating the results, we introduce the random variable $\su^{m,n}$:
\begin{subequations}
\begin{align}
         \su^{m,n}  & :=  \lambda^{-1}(\sz^n-\sz^m) + \Fnor(\sz^m) \label{eq:u-def1}\\  &\;=  (\lambda^{-1}-\tau_{m,n}^{-1})(\sz^n-\sz^m) + \tau_{m,n}^{-1} \scp^{m,n} \label{eq:u-def2}\\
          &\;=  (1-\lambda^{-1}\tau_{m,n}) \Fnor(\sz^m) + \lambda^{-1}\scp^{m,n}, \quad \revise{m < n}.\label{eq:u-def3}
     \end{align}
\end{subequations}
Equations \cref{eq:u-def2} and \cref{eq:u-def3} directly follow from \cref{eq:pretty-important}. 

\begin{lemma}
	\label{lem:basic bounds}
	Let the assumptions \ref{A1}--\ref{A2} and \ref{B1}--\ref{B2} be satisfied and let the sequences $\revise{\{\sx^k\}_{k\in\N}}$, $\revise{\{\sz^k\}_{k\in\N}}$ be generated by $\NSGD$. Let $0 \leq m<n$ be given with $\tau_{m,n} < \lambda$. Then, the following two bounds are valid:
    \begin{enumerate}[label=\textup{\textrm{(\alph*)}},leftmargin=3em,topsep=2pt,itemsep=.5ex,partopsep=0ex]
    \item It holds that $\norm{\Fnor(\sz^n)}^2 \leq  ({1-\frac{\tau_{m,n}}{2\lambda}})\norm{\Fnor(\sz^m)}^2 - \frac{2}{\lambda}\iprod{\su^{m,n}}{\sx^n-\sx^m} + [({\sL+\frac{1}{\lambda}})^2 + \frac{2\lambda \sL^2}{\tau_{m,n}}] \norm{\sx^n-\sx^m}^2 + \frac{3}{2\lambda\tau_{m,n}} \norm{\scp^{m,n}}^2$.
    \item We have $\psi(\sx^n) \leq \psi(\sx^m) + ({\frac{\sL}{2}-\frac{1}{\lambda}})\norm{\sx^n-\sx^m}^2 + \iprod{\su^{m,n}}{\sx^n-\sx^m}.$	
    \end{enumerate}
\end{lemma}

The proof of \cref{lem:basic bounds} is relegated to \cref{app:lem:basic bounds}. Based on \cref{lem:basic bounds}, we see that proper control of the inner product $\langle \su^{m,n}, \sx^n-\sx^m \rangle$ is crucial to achieve descent (of $\psi$). The following estimate illustrates that while this term can generate a desirable decrease, it also introduces an additional nonnegative error term: 
    \begin{align}
    \nonumber \iprod{\su^{m,n}}{\sx^n-\sx^m} & = (\lambda^{-1}-\tau_{m,n}^{-1}) \iprod{\sz^{n}-\sz^m}{\sx^n-\sx^m} + \tau_{m,n}^{-1}\iprod{\scp^{m,n}}{\sx^n-\sx^m}\\
    \nonumber & \leq (\lambda^{-1}-\tau_{m,n}^{-1}) \norm{\sx^n-\sx^m}^2 + \tau_{m,n}^{-1} ({\norm{\scp^{m,n}}^2 + 0.25\norm{\sx^n-\sx^m}^2})\\
    \label{u-descent} & = (\lambda^{-1}-0.75\tau_{m,n}^{-1}) \norm{\sx^n-\sx^m}^2 + \tau_{m,n}^{-1}\norm{\scp^{m,n}}^2.
    \end{align}
Here, we applied \cref{eq:u-def2}, \cref{prox-prop}, and Young's inequality. Thankfully and as we will verify in \Cref{subsec:merit function}, the extra term $\tau_{m,n}^{-1}\norm{\scp^{m,n}}^2$ in \cref{u-descent} can be effectively compensated for by leveraging the normal map bound and choosing $\tau_{m,n}$ sufficiently small. 

\subsection{Merit function and approximate descent}\label{subsec:merit function}
We now establish a descent-type property for the stochastic process $\revise{\{\sz^k\}_{k\in\N}}$ and for $\NSGD$. Our derivations rely on the following merit function which was first introduced in \cite{ouyang2021trust}.
\begin{definition}[Merit function]\label{def:mer-fun} Given $\lambda > 0$, we define the merit function: 
 \[\mer(z):=\psi(\prox{\lambda\vp}(z)) + \frac{\xi\lambda}{2}\|\Fnor(z)\|^2,\quad \text{where} \quad \xi:={(2+2\lambda^2 \sL^2)^{-1}}.\]
\end{definition}
Next, we state the promised approximate descent property of $\NSGD$. 

\begin{proposition}[Approximate descent property]
\label{lem:merit-descent-1}
Let \ref{A1}--\ref{A2} and \ref{B1}--\ref{B2} hold and let $\revise{\{\sz^k\}_{k\in\N}}$ be generated by $\NSGD$. Suppose, \revise{for $0\leq m <n$, the} step sizes $\revise{\{\alpha_k\}_{k\in\N}}$ satisfy $\tau_{m,n} = {\sum}_{i=m}^{n-1}\alpha_i\leq \frac{\sqrt{\xi}}{8\Llambda}$ where $\Llambda := \sL+\frac{2}{\lambda}$. We then have{}
\begin{equation}
    H(\sz^n) - H(\sz^m) \leq  -\frac{\xi\tau_{m,n}}{8} \norm{\Fnor(\sz^m)}^2 + \frac{8}{\tau_{m,n}} \scs_{m,n}^2.
\end{equation}
\end{proposition}

\begin{proof}
By \cref{lem:basic bounds} and setting $\tau = \tau_{m,n}$, $\su = \su^{m,n}$, $\scp = \scp^{m,n}$, and $\scs = \scs_{m,n}$, it holds that
\begin{align}
 \nonumber   H(\sz^n) - H(\sz^m) & \leq  -\frac{\xi\tau}{4} \norm{\Fnor(\sz^m)}^2 + (1-\xi) \iprod{\su}{\sx^n-\sx^m} + \frac{3\xi}{4\tau} \norm{\scp}^2\\
 \label{H-descent1}   &\hspace{4ex} + \Big[{\Big({\frac{\sL}{2}-\frac{1}{\lambda}}\Big) + \frac{\xi\lambda}{2} \Big({\sL+\frac{1}{\lambda}}\Big)^2 + \frac{\xi \lambda^2\sL^2}{\tau}}\Big] \norm{\sx^n-\sx^m}^2.
\end{align}
Invoking the upper bound for $\iprod{\su}{\sx^n-\sx^m}$ in \cref{u-descent}, this yields 
\begin{align}
\nonumber    H(\sz^n) - H(\sz^m) & \leq -\frac{\xi\tau}{4} \norm{\Fnor(\sz^m)}^2 + \frac{4-\xi}{4\tau}\norm{\scp}^2\\
\label{H-descent2}    &\hspace{-6ex} + \underbracket{\Big[\frac{\sL}{2}-\frac{\xi}{\lambda} + \frac{\xi\lambda}{2} \Big({\sL+\frac{1}{\lambda}}\Big)^2 + \frac{4\xi \lambda^2\sL^2 - 3(1-\xi)}{4\tau}\Big]}_{=:T_1} \norm{\sx^n-\sx^m}^2.
\end{align}
Next, we bound $T_1$. Using $\xi= {(2+2\lambda^2\sL^2)^{-1}}$ and $(a+b)^2 \leq 2(a^2+b^2)$, $a,b \in \R$, we have $\frac12\xi\lambda(\sL+\lambda^{-1})^2 \leq \xi \lambda (\sL^2+\lambda^{-2}) = \frac{1}{2\lambda}$ and $4\xi \lambda^2\sL^2 - 3(1-\xi) \leq 4\xi(\lambda^2\sL^2+1)-3 = -1$. Applying the bound $\tau \leq \frac{\sqrt{\xi}}{8\Llambda} \leq \frac{1}{2(\sL+\lambda^{-1})}$, we then obtain $T_1 \leq \frac{\sL+\lambda^{-1}}{2} - \frac{\xi}{\lambda} - \frac{1}{4\tau} \leq - \frac{\xi}{\lambda} < 0$. Therefore, due to $\frac{4-\xi}{4\tau}\leq \frac{1}{\tau}$, the estimate \cref{H-descent2} simplifies to
\[
H(\sz^n) - H(\sz^m) \leq  -{(\xi\tau/4)} \norm{\Fnor(\sz^m)}^2 + {\tau^{-1}} \norm{\scp}^2.
\]
By \cref{lem:est-err}, we further have $\|\scp\|^2 \leq 8\Llambda^2 \tau^4 \|\Fnor(\sz^m)\|^2+8\scs^2$. This implies
\begin{align*}
    H(\sz^n) - H(\sz^m) \leq  - \Big[\frac{\xi}{4} - 8\Llambda^2 \tau^2 \Big]\tau\|{\Fnor(\sz^m)}\|^2 + \frac{8\scs^2}{\tau} 
    \leq  -\frac{\xi\tau}{8} \norm{\Fnor(\sz^m)}^2 + \frac{8\scs^2}{\tau} ,
\end{align*}
where the last inequality is again due to $\tau \leq \frac{\sqrt{\xi}}{8\Llambda}$.
\end{proof}

\section{Complexity bound and global convergence}\label{sec:complex and convergence}
In this section, we derive complexity bounds and prove global convergence of $\NSGD$. Our analysis is inspired by the stochastic approximation literature \cite{ljung1977analysis,KusCla78,benaim1999dynamics,borkar2009stochastic, tadic2015convergence} and uses the concept of alternative time scales $\tau_{t_k,t_{k+1}}$ (for some carefully designed subsequence $\revise{\{t_k\}_{k\in\N}}$) introduced by Tadi{\'c} \cite{tadic2015convergence}. In the following, we relate the aggregated step sizes $\tau_{t_k, t_{k+1}}$ to a given ``time window'' parameter $\sT$ and specify the choice of $\revise{\{t_k\}_{k\in\N}}$.  

The recursive construction of $\revise{\{t_k\}_{k\in\N}}$ closely follows the discussion in \cite{tadic2015convergence}.  

\begin{definition}[Time indices]\label{def:time win}
Let $\sT>0$ and $\delta\in(0,1)$ be given. We introduce the mapping $\varpi:\N \to \N\cup \{\pm\infty\}$, $\varpi(k):=\sup\{n > k: \tau_{k,n} \leq \sT\}$, where $\tau_{m,n} = {\sum}_{k=m}^{n-1} \alpha_k$. The \emph{time indices} $\revise{\{t_k\}_{k\in\N}}$, associated with the \emph{time window} $\sT$ \revise{(and the step sizes $\{\alpha_k\}_{k \in \N}$)}, are defined as 
\[
t_0 := \inf\{i: \alpha_j \leq (1-\delta)\sT \text{ for all $j\geq i$}\}\quad \text{and} \quad t_{k+1} := \varpi(t_k).
\]
\end{definition}

To show convergence of $\NSGD$, we typically first analyze the behavior of the iterates and (merit) function values restricted to \revise{the time indices $\{t_k\}_{k\in\N}$}. Subsequently, full convergence can be deduced by establishing a relationship between a generic index $i$ and the subsequence $\revise{\{t_k\}_{k\in\N}}$. To this end, we define
\begin{equation}\label{index-def}
\fk(i) = \max\{k: t_k < i\}\quad \text{for all $i>t_0$}.
\end{equation}
It is evident that we have $t_{\fk(i)}<i\leq t_{\fk(i)+1}$ for all $i>t_0$.

\begin{lemma}\label{lem:exist time}
Let $\sT > 0$, $\delta \in (0,1)$ be given. The following assertions are valid.
 \begin{enumerate}[label=\textup{\textrm{(\alph*)}},leftmargin=3em,topsep=0pt,itemsep=.5ex,partopsep=0ex]
 \item Assume $\alpha_k \in (0,(1-\delta) \sT]$, $k< K$, and $\alpha_k = \alpha \in(0,(1-\delta)\sT]$ for all $k \geq K$ and some $K\in\N$. We then have $t_k \to \infty$ and $\delta\sT \leq \tau_{t_k,t_{k+1}} \leq \sT$ for all $k \in\N$.
    
 \item
    If $\lim_{k\to\infty}\alpha_k = 0$ and $\sum_{k=0}^\infty \alpha_k = \infty$, it holds that $t_k \to \infty$ and $\delta \sT \leq \tau_{t_k,t_{k+1}} \leq \sT$ for all $k\in\N$.
    \end{enumerate}
\end{lemma}

\begin{proof} 
Here, we use the convention $\sup \emptyset = -\infty$. According to the construction of $\revise{\{t_k\}_{k\in\N}}$ in \cref{def:time win} and $\alpha_k \leq (1-\delta) \sT$, we have $t_0 = 0$ and $t_{k+1}=\varpi(t_k) \neq -\infty$ for all $k\in\N$. Furthermore, we have $t_{k+1} \neq \infty$ since $\sum_{k=0}^{\infty} \alpha_k = \infty$ implies $\tau_{t_k,n}\rightarrow \infty$ as $n\rightarrow \infty$. This also yields $t_{k} \to \infty$, i.e., $\revise{\{t_k\}_{k\in\N}}$ is an infinite subsequence.  By the definition of $\varpi(k)$, it follows $\tau_{t_k,t_{k+1}} =  \tau_{t_k,\varpi(t_{k})} \leq \sT$. In addition, by the optimality of $\varpi(t_k)$, we obtain $\tau_{t_k,t_{k+1}} =  \tau_{t_k,\varpi(t_{k})} \geq \sT - \alpha_{\varpi(t_k)} \geq \delta \sT$ for all $k$.
We continue with part (b). Due to $\alpha_k \to 0$ and the definition of $t_0$, we have $\alpha_k \leq (1-\delta)\sT$ for all $k\geq  t_0$. Hence, all previous derivations are applicable in this case.
\end{proof}

To effectively utilize the approximate descent property shown in \cref{lem:merit-descent-1} and to simplify the notation, we introduce the following constants and terms:
\begin{equation}
    \label{eq:def delta and T}
    \sr_k := \scs_{t_{k},t_{k+1}}, \;\; \delta := \frac45, \;\;  \sT := \frac{\sqrt{\xi}}{8\Llambda},\quad \text{where} \quad \Llambda =\sL+\frac{2}{\lambda}, \;\; \xi =\frac{1}{2+2\lambda^2 \sL^2}.
\end{equation}

\subsection{Complexity bound} We now establish complexity of $\NSGD$.
\begin{theorem}[Complexity bound]
\label{thm:complexity}
Let $K\in\N$ be given and let the assumptions \ref{A1}--\ref{A2} and \ref{B1}--\ref{B2} hold. Let the iterates $\{\sz^k\}_{k=0}^{K}$ be generated by $\NSGD$ using step sizes $\alpha_k \in (0,\tfrac{\sT}{5}]$, $k=0,\dots,K-1$. Then, we have
    \begin{equation}\label{eq:complexity bound}
        \min_{i=0,\ldots,K-1} \Exp[{\|\Fnor(\sz^i)\|^2}] \leq \frac{10[H(z^0)-\bar{\psi}]}{\xi \sum_{i=0}^{K-1}\alpha_i} + \frac{400 \sum_{i=0}^{K-1}\alpha_i^2\sigma_i^2 }{\sT\xi \sum_{i=0}^{K-1}\alpha_i }.
    \end{equation}
\end{theorem}

\begin{proof} First, for all $k \geq K$, we may set $\alpha_k := \hat\alpha \in (0,\sT/5]$. This ensures that \cref{lem:exist time} (a) is applicable with $\delta = \frac45$ and we can infer $t_0 = 0$ and $4\sT/5 \leq \tau_{t_k,t_{k+1}} \leq \sT$ for all $k$. Furthermore, to take advantage of the time window technique, we introduce surrogate iterates for the remaining iterations within the last $\kappa(K)$-th time window. Specifically, for all $k=K,\dots,t_{\fk(K)+1}-1$, we define $\hat \sz^{K} = \sz^{K}$, $\hat \sx^{K} = \sx^{K}$, 
\[
\hat \sz^{k+1} = \hat \sz^k - \alpha_k(\nabla f(\hat \sx^k) + \lambda^{-1}(\hat \sz^k - \hat \sx^k)), \quad \text{and} \quad
    \hat \sx^{k+1} = \proxl(\hat \sz^{k+1}),
\]
where $\alpha_k = \hat\alpha \in(0,\sT/5]$. These surrogate iterates use noiseless gradients, $\sg^k = \nabla f(\hat\sx^k)$, for $k = K, \dots, t_{\kappa(K)+1} - 1$ and are introduced for technical purposes. (They are not actually computed by $\NSGD$). By construction, the variances associated with these steps are zero, i.e., we have $\sigma_k^2 = 0$ for all $k = K, \dots, t_{\kappa(K)+1} - 1$.

For all $0\leq j \leq \kappa(K)$, setting $m = t_{j}$ and $n = t_{j+1}$ in \cref{lem:merit-descent-1} and using $\frac{4}{5} \sT\leq \tau_{t_j,t_{j+1}}$ and $\sr_{j} = \scs_{t_{j},t_{j+1}}$, this yields
\begin{equation}
    \label{eq:tmp-merit-extra}
     \mer(\sz^{t_{j+1}}) - \mer(\sz^{t_{j}}) \leq -\frac{\xi\sT}{{10}}\|\Fnor(\sz^{t_{j}})\|^2 + \frac{10}{\sT}\cdot\sr_{j}^2,\quad \forall
    \; 0\leq j \leq \kappa(K)-1
\end{equation}
and $\mer(\hat \sz^{t_{\fk(K)+1}}) - \mer(\sz^{t_{\fk(K)}}) \leq -\frac{\xi\sT}{{10}}\|\Fnor(\sz^{t_{\fk(K)}})\|^2 + \frac{10}{\sT}\cdot\sr_{\fk(K)}^2$.

Summing these estimates for $j=0,\ldots,\fk(K)-1$ and $j=\fk(K)$, using $t_0 = 0$ and the bound $\mer(z)\geq \bar \psi$, and taking expectation on both sides, we obtain
\begin{equation*}
\frac{\xi\sT}{10}  {\sum}_{j=0}^{\fk(K)} \Exp[\|\Fnor(\sz^{t_j})\|^2] \leq \mer(z^{0}) - \bar \psi + \frac{10}{\sT}{\sum}_{j=0}^{\fk(K)} \Exp[\sr_j^2].
\end{equation*}
Since $\tau_{m,n}=\sum_{i=m}^{n-1}\alpha_i$, we have ${\sum}_{i=0}^{K-1} \alpha_i = \tau_{t_{\fk(K)},K} +  {\sum}_{i=0}^{\fk(K)-1}\tau_{t_i,t_{i+1}} \leq (\fk(K)+1)\sT$. Thus, due to $\{t_0,\dots,t_{\fk(K)}\} \subseteq \{0,\dots,K-1\}$, it follows  
\begin{equation}\label{LHS}
    \frac{\xi}{10}\Big({\sum}_{i=0}^{K-1} \alpha_i \Big) \cdot \min_{i = 0,1,\ldots,K-1} \Exp[\|\Fnor(\sz^{i})\|^2] \leq \mer(z^{0}) - \bar \psi + \frac{10}{\sT}{\sum}_{j=0}^{\fk(K)} \Exp[\sr_j^2].
\end{equation}
Next, applying \Cref{lem:err_estimate}, we have $\Exp[\sr_j^2] \leq 4 \sum_{i=t_j}^{t_{j+1}-1} \alpha_i^2 \sigma_i^2$ for all $j=0,\ldots,\fk(K)-1$ and 
$\Exp[\sr_{\fk(K)}^2] \leq   4 \sum_{i=t_{\fk(K)}}^{t_{\fk(K)+1}-1} \alpha_i^2 \sigma_i^2 = 4\sum_{i=t_{\fk(K)}}^{K-1} \alpha_i^2 \sigma_i^2$. Combining these estimates with \cref{LHS}, we can infer
\[
    \frac{\xi}{10}\Big({\sum}_{i=0}^{K-1} \alpha_i \Big) \cdot \min_{i = 0,1,\ldots,K-1} \Exp[\|\Fnor(\sz^{i})\|^2] \leq \mer(z^{0}) - \bar \psi + \frac{40}{\sT}{\sum}_{i=0}^{K-1}\alpha_i^2 \sigma_i^2.
\]
The proof concludes after dividing both sides of this inequality by $\frac{\xi}{10}({\sum}_{i=0}^{K-1} \alpha_i)$.
\end{proof}
\begin{remark}
    Setting $\lambda = \Theta (1/\sL)$, we have $\xi = \Theta(1)$ and $\sT = \Theta(1/\sL)$ and the complexity \cref{eq:complexity bound} for $\NSGD$  matches the known bound for $\PSGD$ established in \cite[Corollary 3.6]{davdru19} up to constant factors. 
    \revise{Moreover, as discussed in \Cref{sec:stationary}, the different stationarity measures satisfy the following relation:
\[
\|\Fnat(x)\| \leq \|\partial \psi(x)\|_{-} \leq \|\Fnor(z)\|, \quad \text{where} \quad x = \prox{\lambda \vp}(z).
\]
Consequently, the bound in \cref{eq:complexity bound} can also be directly expressed in terms of the stationarity measure $\|\partial \psi(\sx^k)\|_{-}$. This does not seem to be possible for related complexity results, \cite{GhaLanZha16,davdru19}, that rely on the natural residual (or the gradient of the Moreau envelope). In particular, $\|\Fnat(x)\| < \varepsilon$ does not necessarily imply $\|\partial \psi(x)\|_{-} < \varepsilon$.}
\end{remark}
%
\subsection{Global convergence}
Under assumption \ref{B3}, the complexity bound in \cref{thm:complexity} only ensures $\liminf_{k\to\infty}\Exp[\|\Fnor(\sz^{k})\|^2] = 0$ (when taking the limit $K \to \infty$). We address this limitation in \cref{thm:global_convergence-2}. In particular, we prove that both the stationarity measure and function values converge in expectation and in an a$.$s$.$ sense.
 
In the following, we first establish quantitative relationships for both normal map and function values at iterations $i$ and $\kappa(i)$. The proof is relegated to \cref{proof:lem:index-time bound}. 

\begin{lemma}\label{lem:index-time bound}
    Let \ref{A1}--\ref{A2} and \ref{B1}--\ref{B2} hold and let $\revise{\{\sx^k\}_{k\in\N}}$ and $\revise{\{\sz^k\}_{k\in\N}}$ be generated by $\NSGD$. Recalling $\Llambda := \sL+2\lambda^{-1}$, we have for all $i>t_0$:
    \begin{subequations}
    \begin{align}
    \label{eq:tech1}\|\Fnor(\sz^i)\|^2 &\leq 2\prt{2\Llambda\sT +1}^2\| \Fnor(\sz^{t_{\fk(i)}})\|^2 + 8\Llambda^2 \sr_{\fk(i)}^2, \\
    \label{eq:tech2}\psi(\sx^i) &\leq \psi(\sx^{t_{\fk(i)}})+(4\Llambda\sT^2+{\lambda}/{4})\norm{\Fnor(\sz^{t_{\fk(i)}})}^2 + 4\Llambda \sr_{\fk(i)}^2,\\
    \label{eq:tech3}\psi(\sx^i) &\geq \psi(\sx^{t_{\fk(i)+1}})-(4\Llambda\sT^2+{\lambda}/{4})\norm{\Fnor(\sz^{t_{\fk(i)}})}^2 - 16\Llambda \sr_{\fk(i)}^2.
\end{align}
 \end{subequations}
\end{lemma}

\revise{\Cref{lem:index-time bound} will be used in the proof of \Cref{thm:global_convergence-2} to extend certain convergence properties from the subsequence $\{t_k\}_{k \in \N}$ to the entire sequence of indices. Specifically, inequality \eqref{eq:tech1} provides an upper bound on $\|\Fnor(\sz^i)\|^2$ in terms of $\|\Fnor(\sz^{t_{\fk(i)}})\|^2$ and the error term $\sr_{\fk(i)}^2$ (which converges to zero). Upon showing $\|\Fnor(\sz^{t_{\fk(i)}})\| \to 0$, it then directly follows that $\|\Fnor(\sz^i)\| \to 0$. The convergence of the function values $\{\psi(\sx^k)\}_{k \in \N}$ can be established in a similar way by invoking \eqref{eq:tech2}--\eqref{eq:tech3}.}

\begin{theorem}[Global convergence]
\label{thm:global_convergence-2}
Let \ref{A1}--\ref{A2}, \ref{B1}--\ref{B3} hold and let $\revise{\{\sx^k\}_{k\in\N}}$, $\revise{\{\sz^k\}_{k\in\N}}$ be generated by $\NSGD$. Then, there is $\spsi^* : \Omega \to \R$ such that: 
 \begin{enumerate}[label=\textup{\textrm{(\alph*)}},leftmargin=3em,topsep=0.5ex,itemsep=.5ex,partopsep=0ex]
        \item $\Fnor(\sz^k)\to 0$, $\|\partial \psi(\sx^k)\|_{-} \to 0$, and $\psi(\sx^k)\to \spsi^*$ a.s..
	\item $\Exp[\norm{\Fnor(\sz^k)}^2]\to 0$, $\Exp[\|\partial \psi(\sx^k)\|_{-}^2]\to0$, and $\Exp[\psi(\sx^k)] \to \Exp[\spsi^*]$.
\end{enumerate}
\end{theorem}
\begin{remark}  
\revise{$\PSGD$ is known to converge in the sense that $\|\Fnat(\sx^k)\| \to 0$, \cite{LiMil22}. As mentioned, this does not directly imply convergence of the true subgradient norm, i.e., this does not ensure $\|\partial \psi(\sx^k)\|_{-} \to 0$. Moreover, the existing asymptotic convergence results for $\PSGD$ require \emph{global Lipschitz} continuity of the possibly nonsmooth term $\vp$, see \cite[Corollary 3.6]{LiMil22}. Our approach allows us to circumvent this assumption while establishing stronger convergence guarantees.}
\end{remark}
\begin{proof}
By \cref{lem:merit-descent-1}, \cref{lem:exist time} (b) (with $\delta = 4/5$; cf$.$ \cref{eq:def delta and T}), we have 
\begin{equation}
    \label{eq:thm:global-1}
    \mer(\sz^{t_{k+1}}) - \mer(\sz^{t_{k}}) \leq -\frac{\xi\sT}{{10}}\|\Fnor(\sz^{t_k})\|^2 + \frac{10}{\sT}\cdot\sr_{k}^2\quad \forall \,k\in\N.
\end{equation}
The condition $\sum_{k=0}^{\infty}\alpha_k^2\sigma_k^2<\infty$ and \cref{lem:err_estimate} imply that ${\sum}_{k=0}^\infty \sr_{k}^2<\infty$ a$.$s$.$ and ${\sum}_{k=0}^\infty\Exp[{\sr_{k}^2}]<\infty$. Setting $\su_k:=\frac{10}{\sT}\, \sum_{i=k}^\infty \sr_{i}^2$, it follows
\begin{equation}
	\label{eq:strong-0}
	\mer(\sz^{t_{k+1}}) + \su_{k+1} \leq  \mer(\sz^{t_{k}}) + \su_{k}- ({\xi \sT}/{10}) \|\Fnor(\sz^{t_k})\|^2.
\end{equation}	
Summing \cref{eq:thm:global-1} for $k=0,1,\ldots$ and using the bound $\mer(z)\geq \bar \psi$, this yields
\begin{equation}\label{eq:norm-sum}
{\xi\sT}\,{\sum}_{k=0}^\infty \|\Fnor(\sz^{t_k})\|^2 \leq 10[\mer(\sz^{t_{0}}) + \su_{0} - \bar \psi].
\end{equation}
Hence, due to ${\sum}_{k=0}^\infty \sr_{k}^2<\infty$ a$.$s$.$, we can conclude $\|\Fnor(\sz^{t_k})\|\to 0$, $\sr_{k}\to0$, and $\su_{k}\to0$ a$.$s$.$. Moreover, based on the definition \cref{index-def} of $\fk(i)$, we can infer $\|\Fnor(\sz^{t_{\fk(i)}})\|\to 0$ and $\sr_{\fk(i)}\to0$ as $i\to\infty$ a$.$s$.$. This proves $\|\Fnor(\sz^{i})\|\to 0$ a$.$s$.$ owing to \cref{eq:tech1}. 

Taking expectation in \cref{eq:norm-sum} and using $\sum_{k=0}^\infty\Exp{[\sr_{k}^2]}<\infty$, we further obtain $\Exp[\sr_{k}^2]\to0$, $\Exp[\su_{k}^2]\to0$, and $\Exp[{\|\Fnor(\sz^{t_k})\|^2}]\to 0$. Thus, we can infer $\Exp[\|\Fnor(\sz^{i})\|^2]\to 0$ by taking expectation in \cref{eq:tech1} and passing $i\to\infty$. By \cref{eq:stat-am}, we have $\|\partial \psi(\proxl(z))\|_{-} \leq \|\Fnor(z)\|$, which implies $\|\partial \psi(\sx^k)\|_{-} \to 0$ a$.$s$.$ and $\Exp[{\|\partial \psi(\sx^k)\|^2_{-}}] \to 0$.

Next, we show $\psi(\sx^{t_k}) \to \spsi^*$ a$.$s$.$. By \cref{eq:strong-0}, $\{\mer(\sz^{t_k})+\su_k\}_k$ is monotonically non-increasing, and thus, we conclude that $\mer(\sz^{t_k})+\su_k \to \spsi^*$ a$.$s$.$ for some random variable $\spsi^*$. Using $\|\Fnor(\sz^{t_k})\|\to0$, $\su_k\to 0$, and $\psi(\sx^{t_k}) = \mer(\sz^{t_k}) + \su_k - (\frac{\xi\lambda}{2}\|\Fnor(\sz^{t_k})\|^2+\su_k)$, we obtain $\psi(\sx^{t_k}) \to \spsi^*$ a$.$s$.$. Taking $i\to \infty$ in \cref{eq:tech2} and \cref{eq:tech3}, it follows %
\begin{align*}
     \psi(\sx^i) &\leq \psi(\sx^{t_{\fk(i)}})+(4\Llambda\sT^2+\lambda/4)\norm{\Fnor(\sz^{t_{\fk(i)}})}^2 + 4\Llambda \sr_{\fk(i)}^2 \to \spsi^*,\\
    \psi(\sx^i) &\geq \psi(\sx^{t_{\fk(i)+1}})-(4\Llambda\sT^2+\lambda/4)\norm{\Fnor(\sz^{t_{\fk(i)}})}^2 - 16\Llambda \sr_{\fk(i)}^2 \to \spsi^*, \quad \text{a.s.},
 \end{align*}
and hence $\psi(\sx^i) \to \spsi^*$. Finally, we verify $\Exp[\psi(\sx^{i})] \to \Exp[\spsi^*]$. Since $\{\mer(\sz^{t_k}) + \su_k\}_k$ is bounded below by $\bar\psi$, non-increasing, and converges to $\spsi^*$, we can obtain $\Exp[\mer(\sz^{t_k}) + \su_k]\to \Exp[\spsi^*]$ by the monotone convergence theorem. Using similar arguments as before, this yields $\Exp[\psi(\sx^{t_k})] \to \Exp[\spsi^*]$. Let us now define $\ssigma_i := \psi(\sx^{t_{\fk(i)}})+(4\Llambda\sT^2+\lambda/4)\norm{\Fnor(\sz^{t_{\fk(i)}})}^2 + 4\Llambda \sr_{\fk(i)}^2$.
We then have $\psi(\sx^i) \leq \ssigma_i$ (see \cref{eq:tech2}), $\ssigma_i\to \spsi^*$ a$.$s$.$, and $\Exp[\ssigma_i]\to \Exp[\spsi^*]$. Applying Fatou's Lemma for $\psi(\sx^i)-\bar{\psi}$ and $\ssigma_i-\psi(\sx^i)$, it follows 
 \begin{align*}
      \Exp[\spsi^*]-\bar{\psi} &= \Exp[{{\lim}_{i\to\infty}({\psi(\sx^i)-\bar{\psi}}})] \leq {\liminf}_{i\to\infty}{\Exp{[\psi(\sx^i)]}}-\bar{\psi}\quad 
 \end{align*}
and $0 = \Exp[{{\lim}_{i\to\infty} (\ssigma_i-\psi(\sx^i))}] \leq \Exp[{\spsi^*}] - {\limsup}_{i\to\infty} {\Exp{[\psi(\sx^i)]}}$. Consequently, we can infer $\Exp{[\psi(\sx^i)]}\to \Exp[{\spsi^*}]$. This finishes the proof. 
\end{proof}

\section{Iterate convergence and manifold identification}\label{sec:iter conv and mani}
We now study iterate convergence and manifold identification properties of $\NSGD$.
\subsection{\texorpdfstring{A Kurdyka-{\L}ojasiewicz-type property for $H$}{A Kurdyka-Lojasiewicz-type property for H}}
Our analysis is based on the concept of definable functions (in an o-minimal structure)\cite{van1998tame}; see \cite{kur98,coste2000introduction,AttBolRedSou10}. 
\begin{enumerate}[label=\textup{\textrm{(H.3)}},leftmargin=3em,topsep=2mm,itemsep=1ex,partopsep=0ex]
        \item \label{A3} The objective function $\psi : \Rn \to \Rex$ is definable. 
\end{enumerate}
   
Definable functions encompass a broad range of important function classes, including real semialgebraic and globally subanalytic functions \cite{kur98,loj63,gabrielov1996complements}, and functions in log-exp structures \cite{van1998tame}. This underlines the generality and broad applicability of \ref{A3}, as discussed in \cite{AttBolRedSou10,boldanlew06,boldanlewshi07}. We refer to \cite{AttBolRedSou10,li2018calculus,davis2020stochastic} for concrete examples.

Leveraging \ref{A3}, we now state a variant of the Kurdyka-{\L}ojasiewicz (KL) inequality \cite[Lemma 4.11]{josz2024proximal} proposed by Josz et al$.$. The result in \cite{josz2024proximal} is the basis to establish iterate convergence of $\NSGD$ when $\psi$ is merely definable. The desingularization function $\varrho$ in \cref{lem:KL-mer} differs slightly from that used in \cite{josz2024proximal}. In addition, since our approximate descent property is stated for $\mer$, we directly provide a KL-type inequality for $\mer$. This also generalizes \cite[Lemma 5.3]{ouyang2021trust} to the o-minimal setting. 

\begin{lemma}[Kurdyka-{\L}ojasiewicz-type inequality]
\label{lem:KL-mer}
    Let \ref{A3} hold and let $\bar z$ be given with $\bar{x}=\proxl(\bar{z}) \in \crit(\psi)$. For all $\theta\in[\frac12,1)$, there then are $\varsigma, \rho > 0, \eta\in(0,1]$ and a continuous, concave function $\varrho:[0,\eta) \to \R_+$, which is $C^1$ on $(0,\eta)$, satisfying 
      \begin{enumerate}[label=\textup{\textrm{(\alph*)}},leftmargin=3em,topsep=1pt,itemsep=.5ex,partopsep=0ex]
         \item $1/\varrho^\prime(s+t) \leq 1/\varrho^\prime(s) + \varsigma t^{\theta}$ for all $s,t >0$ with $s+t<\eta$,
         \item $\varrho(0) = 0$ and $\varrho^\prime(t) \geq t^{-\theta}$ for all $t\in(0,\eta)$,
     \end{enumerate}
     such that, for all $z\in\cB_\rho(\bar z) \cap\{z:0<|\mer(z)-\mer(\bar z)| + \frac{\xi\lambda}{2}\|\Fnor(z)\|^2<\eta\}$, we have
     \begin{equation}
         \label{eq:mer-kl}
         \varrho^\prime(|\mer(z)-\mer(\bar z)|) \cdot  \|\Fnor(z)\|\geq 1.
     \end{equation} 
\end{lemma}
\begin{proof}
The proof is divided into two steps. First, we show that $\psi$ satisfies the KL property, with the desingularization function $\varrho$ being concave, differentiable, and fulfilling properties (a) and (b). This step closely follows the proof of \cite[Lemma 4.11]{josz2024proximal}. Although the differences are minor, we include a complete proof for self-containedness. In the second step, we show that these properties can be transferred to $\mer$. 

\noindent\textbf{Step 1.}  
Applying \cite[Corollary 15]{boldanlewshi07} and \cite[Theorem 4.1]{AttBolRedSou10}, there are $\rho,\eta_1 >0$ and a definable, concave, continuous function $\varrho_1:[0,\eta_1)\to\R_+$, $\varrho_1 \in C^1((0,\eta_1))$, such that 
\begin{equation}
\varrho_1^\prime(|\psi(x)-\psi(\bar x)|) \cdot \|\partial \psi(x)\|_{-}\geq 1, \quad \forall\; x\in\cB_{\rho}(\bar x) \cap\{x:0<|\psi(x)-\psi(\bar x)|<\eta_1\}. \label{eq:kl-1}
\end{equation}
We now set $\varrho(t)=\int_{0}^t \varrho^\prime(\tau)\,\rmd \tau$, $\varrho^\prime(t):=\max\{\varrho_1^\prime(t),t^{-\theta}\}$, and $\bar \eta:=\min\{\eta_1,1\}$. Then, $\varrho(0)=0$ and \cref{eq:kl-1} remains to hold with $\varrho_1$ and $\eta_1$ being replaced by $\varrho$ and $\bar \eta$, respectively.
The function $\varrho$ is continuous, $C^1$ on $(0,\bar\eta)$, and concave as $\varrho^\prime$ is non-increasing. 
Moreover, $\varrho^\prime$ is definable because $t \mapsto t^{-\theta}$ and $\varrho_1^\prime$ are definable \cite[Lemma 6.1]{coste2000introduction}, so is $t\mapsto 1/\varrho^\prime(t)$ \cite[Lemma 2.3 p.14]{van1998tame} (see also \cite[Section 4.3]{AttBolRedSou10}). According to \cite[Theorem 6.4]{coste2000introduction}, the mapping $t\mapsto 1/\varrho^\prime(t)$ is differentiable for all but finitely many points in $(0,\bar \eta)$. This result, along with the monotonicity theorem \cite[(1.2) p.43]{van1998tame}, indicates that $t\mapsto 1/\varrho^\prime(t)$ is differentiable and monotone on the interval $(0,\eta)$ for some $\eta\in(0,\bar \eta)$.  

If $(1/\varrho^\prime)^\prime(t)\to \infty$ as $t\searrow0$, then $(1/\varrho^\prime)^\prime$ is decreasing on $(0,\eta)$ and $1/\varrho^\prime$ is concave. By \cite[Lemma 3.5]{josz2023convergence}, this yields $1/\varrho^\prime(s+t) \leq 1/\varrho^\prime(s) + 1/\varrho^\prime(t) \leq 1/\varrho^\prime(s) + t^\theta$ for all $s,t>0$ and $s+t<\eta$. Otherwise, we have $\bar\varsigma:=\max_{\tau\in(0,\eta)}(1/\varrho^\prime)^\prime(\tau) < \infty$ and by the mean value theorem, it follows $1/\varrho^\prime(s+t) - 1/\varrho^\prime(s)= (1/\varrho^\prime)^\prime(\tau) t \leq \bar\varsigma t \leq \bar\varsigma t^{\theta}$, $\tau \in (0,\eta)$, where we used $\theta < 1$ and $t \in (0,1)$. Note that $\bar\varsigma$ is nonnegative since $(1/\varrho^\prime)^\prime(\tau) = -\varrho^{\prime\prime}(\tau)/(\varrho^\prime(\tau))^2 \geq 0$ (by the concavity of $\varrho$). We may now set $\varsigma := \max\{\bar\varsigma,1\}$. 

\noindent\textbf{Step 2:} \emph{Property \cref{eq:mer-kl} for $\mer$.} Let $z\in\cB_\rho(\bar z) \cap\{z:0<|\mer(z)-\mer(\bar z)|+ \frac{\xi\lambda}{2}\|\Fnor(z)\|^2<\eta\}$ be arbitrary. Setting $x=\proxl(z)$ and using the nonexpansiveness of $\proxl$, we have $\|x-\bar x\| \leq \|z-\bar z\| \leq \rho$ and $|\psi(x)-\psi(\bar x)| \leq |\mer(z)-\mer(\bar z)|+ \frac{\xi\lambda}{2}\|\Fnor(z)\|^2<\eta$.

Thus, if $x$ satisfies $\psi(x) \neq \psi(\bar x)$, then \textbf{Step 1} implies $\varrho^\prime(|\psi(x)-\psi(\bar x)|) \|\partial \psi(x)\|_{-}\geq 1$. Furthermore, since $1/\varrho^\prime$ is non-decreasing and applying the bound (a), we can infer
\begin{align*}
          1/\varrho^\prime(|\mer(z)-\mer(\bar z)|) 
          &\leq  1/\varrho^\prime(|\psi(x)-\psi(\bar x)| + \tfrac{\xi\lambda}{2}\|\Fnor(z)\|^2)\\
          &\leq 1/\varrho^\prime(|\psi(x)-\psi(\bar x)|) + \varsigma (\tfrac{\xi\lambda}{2}\|\Fnor(z)\|^2)^\theta \\ &\leq \|\partial \psi(x)\|_{-} + \varsigma \sqrt{\xi\lambda/2}\|\Fnor(z)\| \leq (1+\varsigma \sqrt{\xi\lambda/2})\|\Fnor(z)\|.
\end{align*}
Here, the third inequality is due to $\tfrac{\xi\lambda}{2}\|\Fnor(z)\|^2<1$ and $\theta\in[\frac12,1)$, and the last step follows from $\Fnor(z)\in\partial \psi(x)$. If $\psi(x) = \psi(\bar x)$, then we obtain $1/\varrho^\prime(|H(z)-H(\bar z)|) \leq (\tfrac{\xi\lambda}{2}\|\Fnor(z)\|^2)^\theta \leq \sqrt{\xi\lambda/2}\|\Fnor(z)\|$ by (b). We can conclude the proof after rescaling $\varrho$ to $(1+\max\{1,\varsigma\} \sqrt{\xi\lambda/2}) \varrho$ (this does not affect the properties (a)--(b)).
\end{proof}

\subsection{Iterate convergence}
Let us introduce the event 
\[ \cX := \{\omega\in\Omega: {\liminf}_{k\to\infty}\; \|\sx^k(\omega)\| < \infty\}. 
\]
Clearly, if $\omega\in\cX$, then the realization $\revise{\{x^k\}_{k\in\N}} = \revise{\{\sx^k(\omega)\}_{k\in\N}}$ has at least one bounded subsequence and accumulation point. We will work with the following assumption: 
 \begin{enumerate}[label=\textup{\textrm{(S.4)}},leftmargin=3em,topsep=1.5mm,itemsep=1ex,partopsep=0ex]
        \item \label{B4} The event $\cX \in \mathcal F$ occurs with probability $1$, i.e., $\Prob(\cX)=1$.
    \end{enumerate}
%
The condition \ref{B4} is less stringent than the ubiquitous bounded iterates assumption used in KL-based analyses of stochastic algorithms, cf$.$ \cite{tadic2015convergence,li2021convergence,dereich2021convergence}.

We now present one of our main results. The iterate convergence shown in \cref{thm:iter converge} crucially relies on the aforementioned KL property for definable functions. 

\begin{theorem}
	\label{thm:iter converge}
	 Let \ref{A1}--\ref{A3}, \ref{B1}--\ref{B2}, and \ref{B4} hold and let $\revise{\{\sx^k\}_{k\in\N}}$ be generated by $\NSGD$ with step sizes $\revise{\{\alpha_k\}_{k\in\N}}\subset \R_{++}$ satisfying: 
	\begin{equation}\label{eq:kl-step}  
		\lim_{k\to\infty}\alpha_k=0,\quad {\sum}_{k=0}^\infty \alpha_k = \infty \quad \text{and} \quad {\sum}_{k=0}^{\infty} \alpha_k^2 \sigma_k^2 \Big( {\sum}_{i=0}^k\alpha_i \Big)^{\mu}< \infty,
	\end{equation}
for some $\mu>1$. Then, we have $\sx^k \to \sx^*$ a$.$s$.$ for some $\sx^*:\Omega\to\crit(\psi)$.
\end{theorem}

\begin{remark} We continue with several remarks on \cref{thm:iter converge}.
\begin{enumerate}[label=$\bullet$,leftmargin=2em,topsep=0pt,itemsep=.5ex,partopsep=0ex]
    \item \emph{(Step sizes).} 
    The step size requirements in \cref{eq:kl-step} are fairly mild and have appeared in similar form in \cite{tadic2015convergence,qiu2024convergence}. Polynomial step sizes of the form $\alpha_k = \Theta(k^{-\gamma})$ satisfy \cref{eq:kl-step} if $\gamma \in (\frac{2}{3},1]$ and if the variance parameters are constant ($\sigma_k = \sigma$ for all $k$). 
    \item \emph{(Basis for manifold identification).} Previous works on manifold identification of proximal stochastic gradient methods typically postulate iterate convergence, $\sx^k \to \sx^*$ a$.$s$.$, as an assumption, \cite{pooliasch18,huang2022training}. To our knowledge, such convergence guarantees have not been shown so far for $\PSGD$ and its variance-reduced variants in the nonconvex setting. By contrast, invoking \ref{A3}, \cref{thm:iter converge} provides a comprehensive foundation for identification properties of $\NSGD$.
\end{enumerate}
\end{remark}

\subsubsection{\texorpdfstring{Proof of \cref{thm:iter converge}}{Proof of Theorem 4.2}}

We start the proof by deriving certain summability results for the error terms $\su_k$ and $\sr_k$. 

\noindent\textbf{Step 1:} \emph{Summability of error terms.} Due to condition \cref{eq:kl-step} (i.e., $\sum_{k=0}^\infty \alpha_k^2 \tau_{0,k}^{\mu}\sigma_k^2  <\infty$) and \cref{lem:err_estimate}, we have $\bar \sr:=\sum_{k=0}^\infty \tau_{0,t_k}^{\mu}\sr_k^2<\infty$ a$.$s$.$. This yields \vspace{-0.5ex}
\[
\sT\su_k = {10}\, {\sum}_{i=k}^\infty \sr_{i}^2 =  {10}\, {\sum}_{i=k}^\infty \tau_{0,t_i}^{-\mu} \tau_{0,t_i}^{\mu} \sr_{i}^2 \leq {10} \tau_{0,t_k}^{-\mu} \cdot {\sum}_{i=k}^\infty  \tau_{0,t_i}^{\mu} \sr_{i}^2 \leq {10\bar \sr}  \cdot \tau_{0,t_k}^{-\mu}.
\]
Using $\delta \sT \leq \tau_{t_k,t_{k+1}}$, it follows $\tau_{0,t_k} \geq  \tau_{0,t_0} + \delta \sT k$ and $\su_k\to 0$ a$.$s$.$, which implies 
\begin{equation}
    \label{eq:summability u_k}
    {\sum}_{k=0}^\infty \su_k^\theta \leq (10\bar\sr/\sT)^\theta  \cdot {\sum}_{k=0}^\infty (\tau_{0,t_0} + \delta \sT k)^{-\mu\theta}<\infty \quad \text{for all}\quad \theta > \tfrac{1}{\mu}.
\end{equation}
Hence, we have $\sum_{k=j}^\infty \su_k^\theta \to 0$ as $j\to\infty$ for all $\theta > {\mu}^{-1}$. Moreover, it holds that \vspace{-1ex}
\begin{equation}\label{eq:rk-sum}
{\sum}_{k=0}^\infty \sr_k = {\sum}_{k=0}^\infty \tau_{0,t_k}^{-\frac{\mu}{2}}\cdot \tau_{0,t_k}^{\frac{\mu}{2}} \sr_k \leq \Big({ \bar{\sr} \cdot {\sum}_{k=0}^\infty (\tau_{0,t_0} + \delta \sT  k)^{-\mu} }\Big)^{1/2} < \infty \quad \text{a$.$s$.$.}
\end{equation}
Let us fix $\theta\in(\max\{\frac{1}{\mu},\frac12\},1)$ and let us introduce the events $\mathcal E := \{\omega : \sum_{k=0}^\infty \su_k(\omega)^\theta < \infty, \, \sum_{k=0}^\infty \sr_k(\omega) < \infty\}$ and $\mathcal C := \mathcal X \cap \{\omega : \Fnor(\sz^k(\omega)) \to 0, \, \psi(\sx^k(\omega))\to\spsi^*(\omega)\}$. Condition \cref{eq:kl-step} implies that \ref{B3} is satisfied. Thus, by \cref{thm:global_convergence-2}, \cref{eq:summability u_k}, \cref{eq:rk-sum}, and \ref{B4}, we can infer $\Prob(\mathcal C) = \Prob(\mathcal E) = \Prob(\mathcal D)=1$, where $\mathcal D = \mathcal C \cap \mathcal E$. 

\noindent\textbf{Step 2:} \emph{Subsequence convergence and roadmap.} We first present an auxiliary result.
 
\begin{lemma}[Subsequence convergence]\label{lem:subsequence converge}
Let \ref{A1}--\ref{A2} and \ref{B1}--\ref{B4} hold and let $\omega\in\mathcal C$ be arbitrary. Then, there exist $z^*\in\Rn$ and a subsequence $\revise{\{q_k\}_{k\in\N}}\subseteq \N$ such that $x^*:=\proxl(z^*)\in\crit(\psi)$, $x^{q_k} \to x^*$, $z^{q_k} \to z^*$, $\psi(x^k) \to \psi(x^*) = \psi^*$, and $\mer(z^{k}) \to \mer(z^*) = \psi^*$ as $k\to\infty$, where  
$\revise{\{x^k\}_{k\in\N}} = \revise{\{\sx^k(\omega)\}_{k\in\N}}$ and $\revise{\{z^k\}_{k\in\N}}= \revise{\{\sz^k(\omega)\}_{k\in\N}}$. 
\end{lemma}

The proof of \cref{lem:subsequence converge} can be found in \cref{proof:lem:subsequence converge}. To simplify the notation, we set
$\sd_k := {\max}_{t_k<i\leq t_{k+1}} \max\{\|\sx^i-\sx^{t_k}\|,\|\sz^i-\sz^{t_k}\|\}$.

\cref{lem:est-err} and $\tau_{t_k,t_{k+1}} \leq \sT$ then imply 
\begin{equation}\label{dk-bound}
    \sd_k \leq 2 \sT \|\Fnor(\sz^{t_k})\| + 2 \sr_k.
\end{equation}
Our goal is to show $\sx^k \to \sx^*$ a$.$s$.$. Thanks to $\Prob(\mathcal D)=1$, it suffices to prove that $\revise{\{\sx^k(\omega)\}_{k\in\N}}$ is a Cauchy sequence for all $\omega \in \mathcal D$. Let $\omega \in \mathcal D$ be arbitrary and let $\revise{\{x^k\}_{k\in\N}} = \revise{\{\sx^k(\omega)\}_{k\in\N}}$,  $\revise{\{z^k\}_{k\in\N}} = \revise{\{\sz^k(\omega)\}_{k\in\N}}$, $\revise{\{r_k\}_{k\in\N}} = \revise{\{\sr_k(\omega)\}_{k\in\N}}$, etc. denote the corresponding realizations. For given $\varepsilon>0$, we need to show that there is $k$ such that $\|x^m - x^n\|<\varepsilon$ for all $k\leq m<n$. 

We claim that it suffices to establish ${\sum}_{k=0}^{\infty}\,d_k< \infty$. In this case, there is $k^\prime\in\N$ such that ${\sum}_{j=k}^{\infty}\,d_j<\varepsilon/2$ for all $k \geq k^\prime$. For any $n > m > t_{k^\prime}$, we have $k^\prime \leq \fk(m)\leq \fk(n)$. Therefore, using the triangle inequality and the above bound, we obtain 
\begin{equation*}
\|x^m - x^n\| \leq d_{\fk(m)} + d_{\fk(n)} + {\sum}_{j=\fk(m)}^{{\fk(n)-1}}\|x^{t_{j+1}}-x^{t_j}\| \leq 2 {\sum}_{j=\fk(m)}^{{\fk(n)}} d_j< \varepsilon, 
\end{equation*}
which proves our claim. In the remaining steps, we establish summability of $\revise{\{d_k\}_{k\in\N}}$.

\noindent\textbf{Step 3:} \emph{Entering the local region and a KL-based bound for $\revise{\{d_k\}_{k\in\N}}$.} By \cref{lem:subsequence converge}, there exists a subsequence $\revise{\{q_k\}_{k\in\N}}$ such that $x^{q_k} \to x^*$, $z^{q_k} \to z^*$ for some $z^*$ with $x^* = \proxl(z^*) \in\crit(\psi)$, and $\mer(z^{k}) \to \mer(z^*) = \psi^*$. Based on the definition of $\fk$ and $d_k$, it holds that $\max\{\|x^{q_k} - x^{t_{\fk(q_k)}}\|,\|z^{q_k} - z^{t_{\fk(q_k)}}\|\} \leq d_{\fk(q_k)}$. Due to $\|\Fnor(z^{t_{{\fk(q_k)}}})\|\to0$, $r_{{\fk(q_k)}}\to0$, and \cref{dk-bound}, we then have $d_{{\fk(q_k)}} \to 0$, $x^{t_{\fk(q_k)}} \to x^*$, and $z^{t_{\fk(q_k)}} \to z^*$.

Applying \cref{lem:KL-mer}, there exist $\rho>0$, $\eta\in(0,1]$, and a concave function $\varrho:[0,\eta)\to\R_+$, $\varrho \in C^1((0,\eta))$ (satisfying properties (a)--(b) in \cref{lem:KL-mer}) such that
\begin{equation}
         \label{eq:super-proof-kl}
         \varrho^\prime(|\mer(z)-\psi^*|) \cdot  \|\Fnor(z)\|\geq 1\quad 
\end{equation}
for all $z\in\cB_\rho(z^*) \cap\{z \in \Rn:0<|\mer(z)-\psi^*| + \frac{\xi\lambda}{2}\|\Fnor(z)\|^2<\eta\}$.
    
Next, we define $\Delta_k:= \varrho(\mer(z^{t_k}) -\psi^* + u_k)$.
Since \revise{$\{\mer(z^{t_k}) -\psi^* + u_k\}_{k \in \N}$} is monotonically decreasing (see \cref{eq:strong-0}) with $\mer(z^{t_k}) + u_k \to \psi^*$, $\Delta_k$ is well-defined. By $\omega \in \mathcal D$, $\varrho(0) = 0$, and by the previous discussion, there then is $\hat k \in \N$ such that
\begin{subequations}
\label{eq:choose rho eta}
    \begin{align}
        \label{eq:choose rho eta-1} &\|z^{t_{\hat k}} - z^*\| + 20\xi^{-1}\Delta_{\hat{k}}  + 2\varsigma\sT {\sum}_{i=\hat k}^\infty\, u_i^\theta + 2{\sum}_{i=\hat k}^\infty\, r_i< \rho \quad \text{and} \quad \\ \label{eq:choose rho eta-2} &|\mer(z^{t_k})-\psi^*| + \max\{\tfrac{\xi\lambda}{2}\|\Fnor(z^{t_k})\|^2, u_k\} \leq  \eta\quad\text{for all} \quad k \geq \hat k. 
    \end{align}
\end{subequations}

Let $k \geq \hat k$ be now given with $z^{t_k} \in \mathcal B_\rho(z^*)$. (In particular, this is true for $k = \hat k$). Hence, by \cref{eq:choose rho eta-2}, the KL inequality \cref{eq:super-proof-kl} is applicable and we have
\begin{subequations}
    \begin{align}
	\label{eq:long-equation-1}\Delta_k - \Delta_{k+1}
    & \geq \varrho^\prime(\mer(z^{t_k}) - \psi^*+u_k)(\mer(z^{t_k}) +u_k - [\mer(z^{t_{k+1}})+u_{k+1}]) \\
	\label{eq:long-equation-2}& \geq \varrho^\prime(|\mer(z^{t_k}) - \psi^*|+u_k) (\xi \sT \|\Fnor(z^{t_k})\|^2/10) \\
	\label{eq:long-equation-3}&\geq (1/\varrho^\prime({|\mer(z^{t_k}) - \psi^*|})  + \varsigma u_k^\theta)^{-1} (\xi \sT \|\Fnor(z^{t_k})\|^2/10) \\
    \label{eq:long-equation-4}&\geq \frac{\xi\sT}{10}\cdot \frac{\|\Fnor(z^{t_k})\|^2}{\|\Fnor(z^{t_k})\| +  \varsigma u_k^\theta}\geq \frac{\xi \sT}{10}\srt{\|\Fnor(z^{t_k})\| -  \varsigma 
 u_k^\theta}.
\end{align}
\end{subequations}
The inequality \cref{eq:long-equation-1} is due to the concavity of $\varrho$, \cref{eq:long-equation-2} follows from the descent-type property \cref{eq:strong-0} and the fact that $\varrho^\prime$ is non-increasing, \cref{eq:long-equation-3} utilizes $\varrho^\prime(s+t) \geq [1/\varrho^\prime(s) + \varsigma t^{\theta}]^{-1}$ (see \cref{lem:KL-mer} (a)), and \cref{eq:long-equation-4} is based on \cref{eq:super-proof-kl} and ${a^2}\geq a^2-b^2$, $a,b\in \R$. Rearranging yields $\sT\|\Fnor(z^{t_k})\| \leq 10\xi^{-1}(\Delta_k - \Delta_{k+1}) + \varsigma \sT u_k^\theta$. Combining this estimate with $d_k \leq 2 \sT \|\Fnor(z^{t_k})\| + 2 r_k$ (see \cref{dk-bound}), we can infer
\begin{equation}
    \label{eq:d bound}
d_k \leq 20\xi^{-1}(\Delta_k - \Delta_{k+1}) + 2\varsigma\sT u_k^\theta + 2r_k.
\end{equation}

\noindent\textbf{Step 4:} \emph{Attracted by the local region.}  We now show that the following statements hold for all $k \geq \hat k$: 
	\begin{enumerate}[label=(\alph*),topsep=4pt,itemsep=0.5ex,partopsep=0ex]
	\item $z^{t_k} \in \mathcal{B}_{\rho}(z^*)$ and $|\mer(z^{t_k}) - \psi^*| + \max\{\frac{\xi\lambda}{2}\|\Fnor(z^{t_k})\|^2,u_k\}< \eta$.
	\item $\sum_{i=\hat k}^{k}d_i  \leq 20\xi^{-1}[\Delta_{\hat k} - \Delta_{k+1}] + 2\varsigma\sT \sum_{i=\hat k}^k u_i^\theta + 2\sum_{i=\hat k}^k r_i.$
	\end{enumerate}
We will prove these statements by induction. First, (a) and (b) are satisfied for $k=\hat k$ thanks to \cref{eq:choose rho eta} and \cref{eq:d bound}. To conduct the induction step, let us assume that (a)--(b) hold for some $k =n \geq \hat k$. We first verify (a) for $k=n+1$. Due to \cref{eq:choose rho eta-2}, the bound $|\mer(z^{t_{n+1}}) - \psi^*| + \max\{\frac{\xi\lambda}{2}\|\Fnor(z^{t_{n+1}})\|^2,u_{n+1}\} < \eta$ is immediate. Moreover, we have
\begin{align*}
    \|z^{t_{n+1}} - z^*\| &\leq \|z^{t_{\hat k}} - z^*\| + {\sum}_{i=\hat{k}}^n \|z^{t_{i+1}} - z^{t_{i}}\| \leq  \|z^{t_{\hat k}} - z^*\| + {\sum}_{i=\hat k}^n \,d_i\\
    &\leq \|z^{t_{\hat k}} - z^*\| +  20\xi^{-1}[\Delta_{\hat k} - \Delta_{n+1}] + 2\varsigma\sT {\sum}_{i=\hat k}^\infty\, u_i^\theta + 2{\sum}_{i=\hat k}^\infty\, r_i < \rho,
\end{align*}
where the second inequality follows from the definition of $d_k$, the third step uses (b) for $k=n$, and the last step holds due to \cref{eq:choose rho eta-1} and $\Delta_{n+1}\geq 0$. Next, we prove that (b) holds for $k=n+1$. Invoking \cref{eq:d bound} with $k=n+1$ and (b) with $k=n$, we obtain
\[
{\sum}_{i=\hat k}^{n+1}\,d_i = d_{n+1} + {\sum}_{i=\hat k}^{n}\,d_i \leq 20\xi^{-1}[\Delta_{\hat k} - \Delta_{n+2}] + 2\varsigma\sT {\sum}_{i=\hat k}^{n+1}\, u_i^\theta + 2{\sum}_{i=\hat k}^{n+1}\, r_i.
\]
This completes the induction. Thus, (a), (b), and \cref{eq:d bound} hold for all $k \geq \hat k$. 

\noindent\textbf{Step 5:} \emph{Summability of $\revise{\{d_k\}_{k\in\N}}$.} Taking $k\to\infty$ in \textbf{Step 4} (b) and due to $\Delta_k \geq 0$, we have \vspace{-1.5ex}
\begin{equation*} 
{\sum}_{i=\hat k}^{\infty}d_i  \leq 20\xi^{-1}\Delta_{\hat k} + 2\varsigma\sT {\sum}_{i=\hat k}^{\infty} u_i^\theta + 2{\sum}_{i=\hat k}^{\infty} r_i < \rho < \infty.
\end{equation*}
This yields $\sum_{k=0}^{\infty}d_k < \infty$ and, following \textbf{Step 2}, this completes the proof. \hfill \proofbox

\subsubsection{Existing KL-based results for stochastic methods} 
The KL property has been used extensively in the past decades to study the limiting behavior of deterministic optimization algorithms; see, e.g., \cite{absil2005convergence,AttBol09,AttBolRedSou10,BolSabTeb14,fragarpey15}. 
In contrast to the deterministic case, applications of KL-based techniques to stochastic methods are rare, with most results being limited to the smooth setting. To our knowledge, Tadi\'{c}, \cite{tadic2015convergence}, provides one of the first KL-based iterate convergence results for $\SGD$.
Recently, in \cite{dereich2021convergence,qiu2024convergence} a$.$s$.$ convergence of $\SGD$-type methods with momentum is studied under the stronger {\L}ojasiewicz (\L) inequality (the desingularization function must have the form $\varrho(t) = ct^{1-\theta}$). 
Li et al., \cite{li2021convergence}, derive KL-based iterate convergence for random reshuffling ($\RR$). $\RR$ is a without-replacement variant of $\SGD$ for finite-sum problems. Similar to $\SGD$, it does not satisfy a sufficient descent condition, but the stochastic errors in $\RR$ can be controlled in a deterministic fashion. 
Chouzenoux et al., \cite{chouzenoux2023kurdyka}, study KL-based iterate convergence of $\PSGD$, though under restrictive conditions: $\varphi$ has to be smooth, $\inf_k \alpha_k > 0$, and the variance $\sigma_k$ has to be summable. Thus, the analysis in \cite{chouzenoux2023kurdyka} aligns more with the deterministic framework by Frankel et al$.$, \cite{fragarpey15}. 
In \cite{latafat2022block,latafat2022finito}, iterate convergence is derived for a family of block-coordinate (incremental-type) algorithms for \cref{eq:SO}. The results in \cite{latafat2022block,latafat2022finito} are based on the {\L}-inequality and leverage a sure descent property. Recent KL-based analyses of proximal $\RR$ methods can be found in \cite{qiu2023new,josz2024proximal}. Driggs et al$.$, \cite{driggs2021stochastic}, show iterate convergence in expectation for a proximal alternating method with variance reduction via the KL property.

\subsection{Manifold identification} \label{sec:manifold}
We first introduce the notion of partial smoothness following \cite{lewis2002active,lewis2013partial,liang2017activity}. For simplicity, we specialize the definition to \revise{the nonsmooth part $\vp$ of} the objective function $\psi$ in \cref{eq:SO}.  A set $\cM_{\bar x}$ is a ${\cal C}^2$-\emph{manifold} around a point $\bar x\in\cM_{\bar x}$ if (locally) $\cM_{\bar x}$ is the solution set of a collection of ${\cal C}^2$-equations with linearly independent gradients. 

\begin{definition}[Partial smoothness] The function \revise{$\vp$} is \emph{partly smooth} at $\bar x$ relative to a set $\cM_{\bar x}$ containing $\bar x$ if $\partial \vp(\bar x) \neq \emptyset$, and
\begin{enumerate}[label=\textup{\textrm{(\alph*)}},leftmargin=3em,topsep=1pt,itemsep=.5ex,partopsep=0ex]
    \item \emph{(smoothness)} $\cM_{\bar x}$ is a ${\cal C}^2$-manifold around $\bar x$ and $\revise{\vp}|_{\cM_{\bar x}}$ is ${\cal C}^2$ near $\bar x$;
    \item \emph{(sharpness)} the affine span of $\partial \revise{\vp}(\bar x)$ is a translate of the normal space ${\cal N}_{\cM_{\bar x}}(\bar x)$;
    \item \emph{(continuity)} the set-valued mapping $\partial \revise{\vp}$ restricted to $\cM_{\bar x}$ is continuous at $\bar x$.
\end{enumerate}
Here, we refer to the set $\cM_{\bar x}$ as the \emph{active manifold}. 
\end{definition}

The class of partly smooth functions enjoys broad applicability in optimization, machine learning, and statistics \revise{\cite{lewis2013partial,vaiter2017model,lee2023accelerating}}; we also refer to the discussions in \cite[Section 3]{lewis2002active}, \cite[Section 2.4]{lee2012manifold}, \cite[Example 3.7]{liang2017activity}. \revise{Our identification result is based on the following lemma by Lee \cite[Lemma 2]{lee2023accelerating}, which is established for composite-type functions with a convex nonsmooth component. A key advantage of this lemma is that it does not require the function $ f $ in \cref{eq:SO} to be $ {\cal C}^2 $ on the manifold}.



\begin{theorem}[Manifold identification \revise{\cite{lee2023accelerating}}] \label{lem:identification}
Let \ref{A1}--\ref{A2} hold and let \revise{$\vp$} be partly smooth at $\bar x$ relative to the manifold $\cM_{\bar x}$ with $0\in \mathrm{ri}(\partial \psi(\bar x))$. \revise{Suppose $x^k \to \bar x$ and $\|\partial \psi(x^k)\|_{-} \to0$, then we have $x^k \in \cM_{\bar x}$ for all large $k\in\N$}.
\end{theorem}

\begin{theorem}\label{thm:identification}
    Let the assumptions in \cref{thm:iter converge} hold and let $\revise{\{\sx^k\}_{k\in\N}}$ be generated by $\NSGD$. Let $\sx^*$ denote the (a$.$s$.$) limit of $\revise{\{\sx^k\}_{k\in\N}}$ and let \revise{the function $\vp$} be {partly smooth} at $\sx^*(\omega)$ relative to $\cM_{\sx^*(\omega)}$ with $0\in \mathrm{ri}(\partial \psi(\sx^*(\omega))$ for almost every $\omega\in\Omega$. Then, we have $\sx^k \in \cM_{\sx^*}$ for all sufficiently large $k \geq \sk$ a$.$s$.$.  
\end{theorem}

\begin{proof}
\revise{By \cref{thm:global_convergence-2,thm:iter converge}, it holds that $\sx^k\to \sx^*$ and $\|\partial \psi(\sx^k)\|_{-}\to 0$ a$.$s$.$}. Consequently, invoking \cref{lem:identification}, we conclude that $\sx^k(\omega)\in\cM_{\sx^*(\omega)}$ for all $k \geq \sk(\omega) \in \N$ and for almost every $\omega \in \Omega$. 
\end{proof}

\begin{remark}
Manifold identification properties of stochastic proximal-type methods have been a topic of significant interest. As mentioned, previous works successfully established identification results, but primarily in the context of convex optimization \cite{lee2012manifold,duchi2021asymptotic,dai2023}. The nonconvex case poses additional challenges and existing identification properties require the generated iterates to converge a$.$s$.$ to some stationary point, see \cite{pooliasch18,huang2022training}. So far, this assumption has remained an unverified abstract prerequisite. Our work allows addressing this gap. For $\NSGD$, we show $\|\partial \psi(\sx^k)\|_{-}\to 0$ a$.$s$.$ (\cref{thm:global_convergence-2}) and we prove a$.$s$.$ iterate convergence when $\psi$ is definable (\cref{thm:iter converge}). Combining these results enables us to conclusively demonstrate manifold identification under mild assumptions, \revise{following the framework provided in \cite{lee2023accelerating}}.
\end{remark}


\section{Numerical illustrations}\label{sec:num_exp}
In this section, we present numerical comparisons of $\NSGD$, $\PSGD$\revise{, and $\RDA$}. Our primary aim is to illustrate that $\NSGD$ generally has better identification properties than $\PSGD$. All experiments are performed using \texttt{MATLAB} (R2024b) on an iMac, Apple M1, 2021 with 16GB memory. 

\subsection{Nonconvex classification}\label{subsec:exp-2}
We first \revise{compare $\NSGD$ and $\PSGD$ on} a sparse nonconvex binary classification problem, \cite{wang2017stochastic,milzarek2019stochastic}, 
\begin{equation}
	\label{eq:binary-clas}
	{\min}_{x\in \Rn}~\psi(x) := \frac1{N}\,{\sum}_{i=1}^{N}\,[1-\tanh(b_i\cdot a_i^\top x)] + \nu\|{x}\|_1,
\end{equation}
where $(a_i,b_i) \in \Rn \times \{-1,1\}$, $i \in [N] := \{1,\dots,N\}$, denotes the given data and $\tanh$ is the hyperbolic tangent function. In the comparison, we use the datasets\footnote{All datasets are available at \url{www.csie.ntu.edu.tw/~cjlin/libsvmtools/datasets}} \texttt{news20} ($N = 19\,996$, $d = 1\,355\,191$), \texttt{rcv1} ($N=20\,242$, $d=47\,236$), and \texttt{gisette} ($N=6\,000$, $d=5\,000$) and the parameter $\nu$ is set to $1/N$ for all datasets. 
The $\ell_1$-norm is known to be partly smooth at $\bar x$ relative to the manifold $\mathcal M_{\bar x} := \{x: \mathrm{supp}(x) \subseteq \mathrm{supp}(\bar x)\}$, \cite{vaiter2017model,liang2017activity}. Hence, the identification results in \cref{thm:identification} are generally applicable when solving problem \cref{eq:binary-clas}. \revise{At iteration $k$, we construct a mini-batch $S_k \subseteq [N]$ by sampling $|S_k|$ elements uniformly at random without replacement from $[N]$. We then use $g^k:=\frac{1}{|S_k|}\sum_{i\in S_k} \nabla f_i(x^k)$ as a standard stochastic approximation of $\nabla f(x^k)$.}

\textbf{Implementational details.} We use the initial point $x^0=\mathds{1}/d\in\Rn$ and mini-batches $S_k$ with \revise{several} fixed \revise{sizes} \revise{$|S_k| \in \{16,64,256\}$}.   We apply step sizes of the form $\alpha_k=\alpha/(\sL+k)$ and choose $\lambda=\alpha$ in $\NSGD$, where \revise{$\alpha \in [10^{0},10^5]$} and $\sL  = 4\|A\|^2/(5N)$ is the Lipschitz constant of $\nabla f$ and $A=[a_1^\top,\dots,a_N^\top]^\top$. Here, the proximity operator is the well-known shrinkage operator $\prox{\lambda \nu \|\cdot\|_1}(z) = \mathrm{sgn}(z) \odot \max\{0,|z|-\lambda\nu\}$. Note that $\PSGD$ computes $\prox{\alpha_k \nu \|\cdot\|_1}$, $k\in\N$, while we work with $\prox{\lambda \nu \|\cdot\|_1}$ in $\NSGD$. 

\begin{figure}[t]
\centering
	\setlength{\abovecaptionskip}{-3pt plus 3pt minus 0pt}
	\setlength{\belowcaptionskip}{-10pt plus 3pt minus 0pt}
	\centering
	\begin{tikzpicture}[scale=1]
    \node[right] at (3.44,2) {\footnotesize $\NSGD$};
    \node[right] at (9.24,2) {\footnotesize $\PSGD$};
    \draw[draw=MyGray] (0.54,1.75) rectangle (12.14,2.25);
    \draw[draw=DeepBlue,fill=DeepBlue,line width=1pt] (2.34,2) -- (3.24,2);
    \node[circle,draw=DeepBlue,fill=DeepBlue,minimum size=4pt,inner sep=0pt, outer sep=0pt] at (2.79,2) {};
    \draw[draw=DeepRed,fill=DeepRed,line width=1pt] (8.14,2) -- (9.04,2);
    \node[draw=DeepRed,fill=DeepRed,minimum size=4pt,inner sep=0pt, outer sep=0pt] at (8.59,2) {};

	\node[right] at (0.0,0) {\includegraphics[height=3.2cm,trim=20 0 20 0,clip]{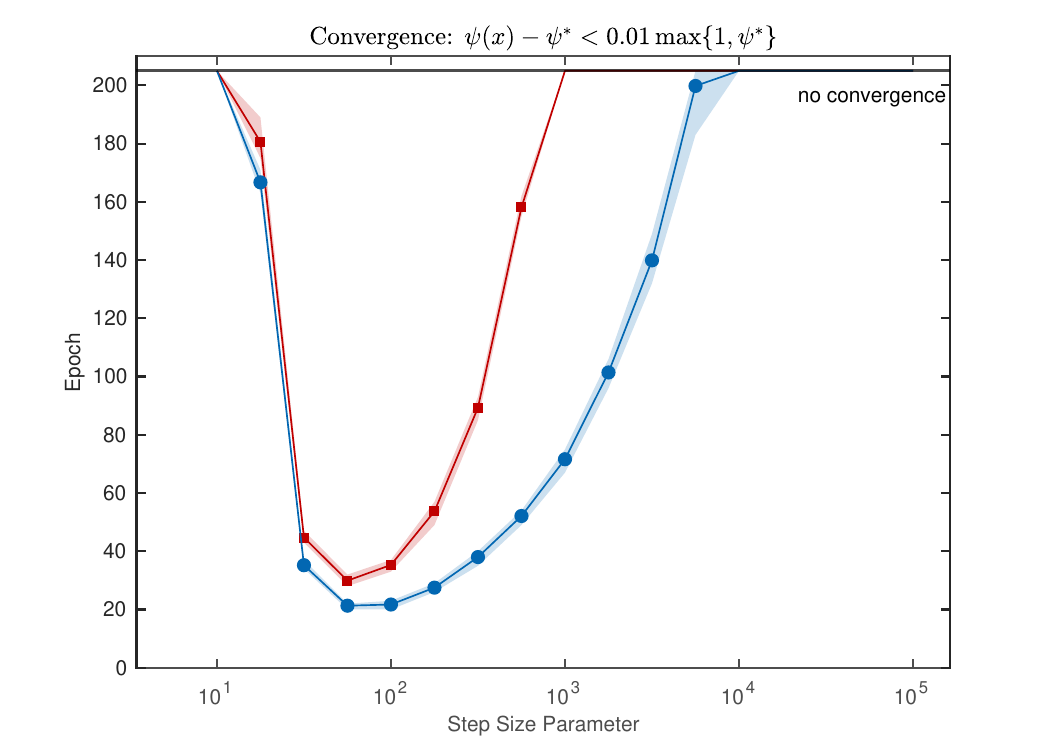}};
	\node[right] at (4.25,0) {\includegraphics[height=3.2cm,trim=40 0 20 0,clip]{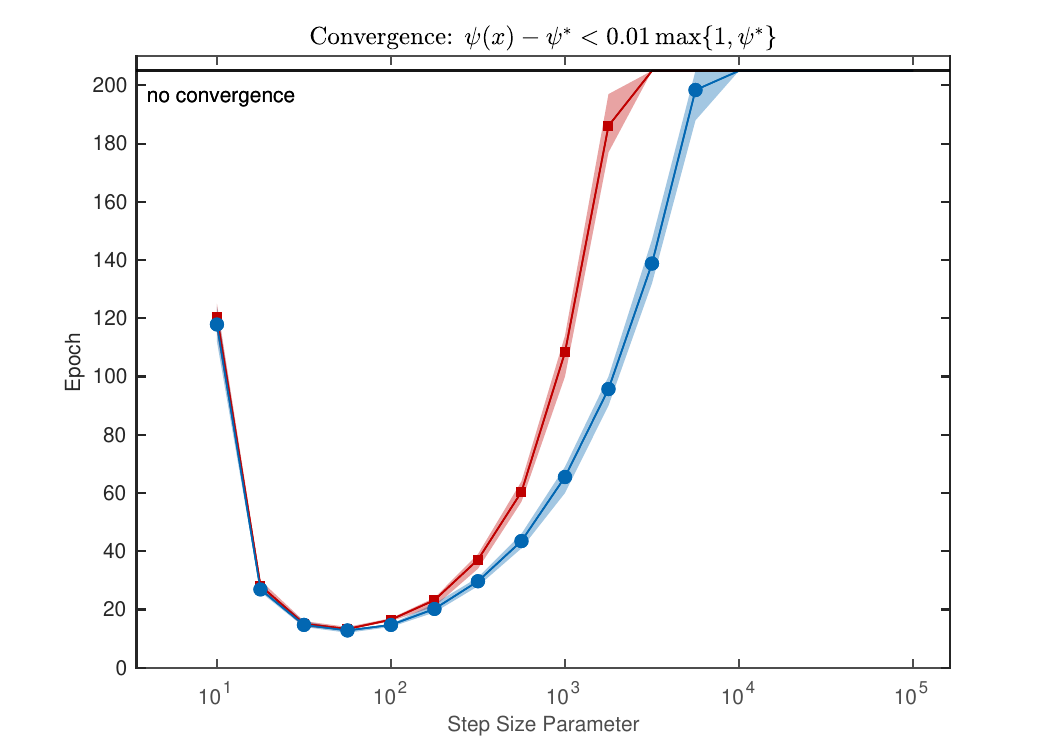}};
	\node[right] at (8.3,0) {\includegraphics[height=3.2cm,trim=40 0 20 0,clip]{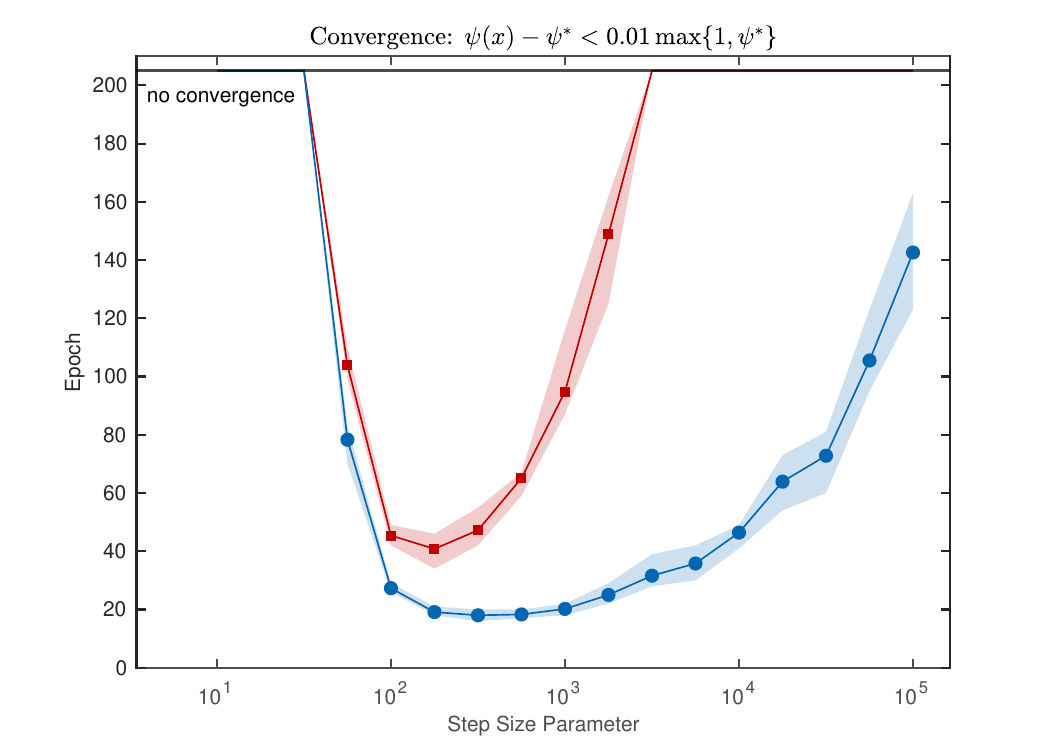}};
	\node[right] at (0.0,-3.2) {\includegraphics[height=3.2cm,trim=20 0 20 0,clip]{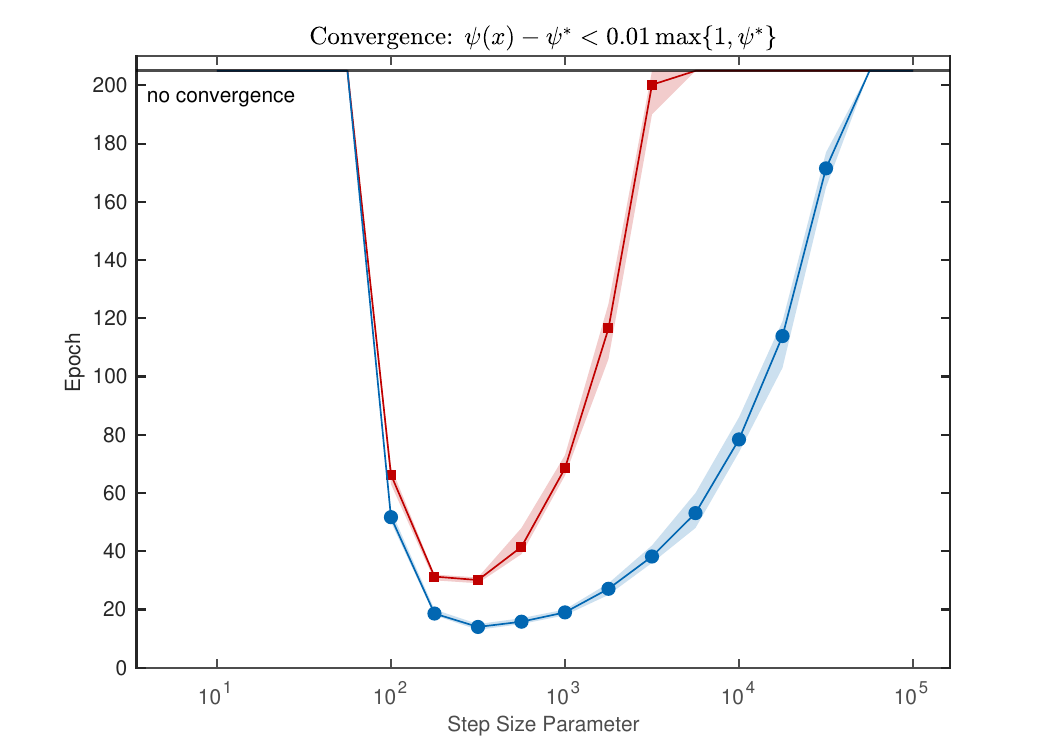}};
	\node[right] at (4.25,-3.2) {\includegraphics[height=3.2cm,trim=40 0 20 0,clip]{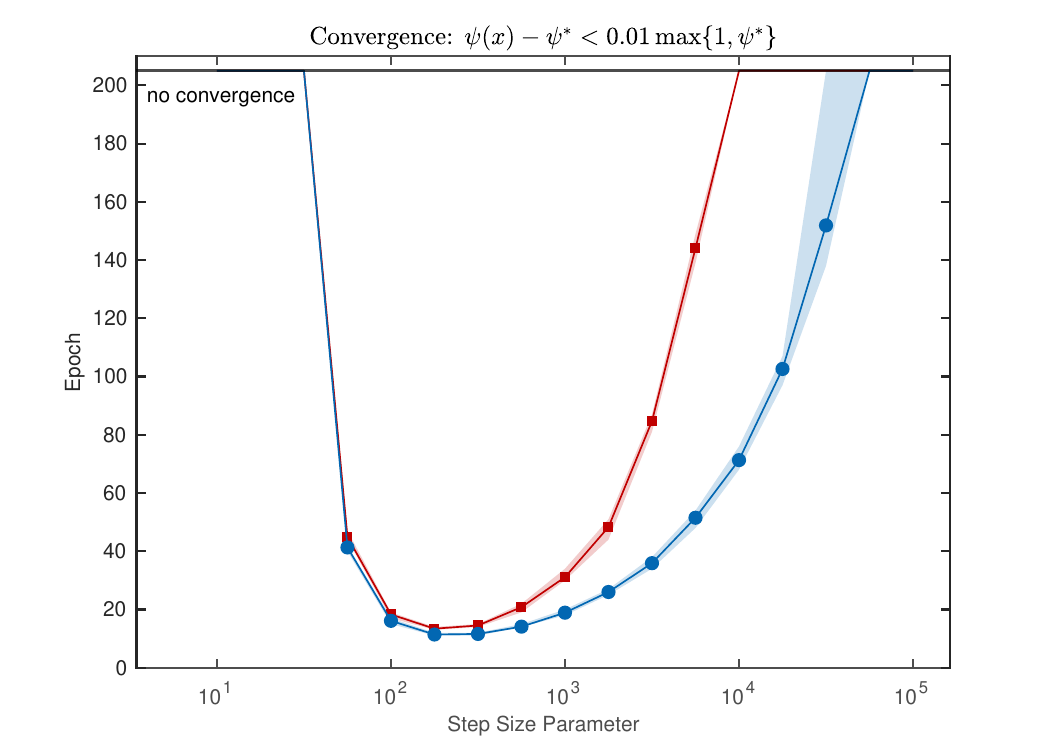}};
	\node[right] at (8.3,-3.2) {\includegraphics[height=3.2cm,trim=40 0 20 0,clip]{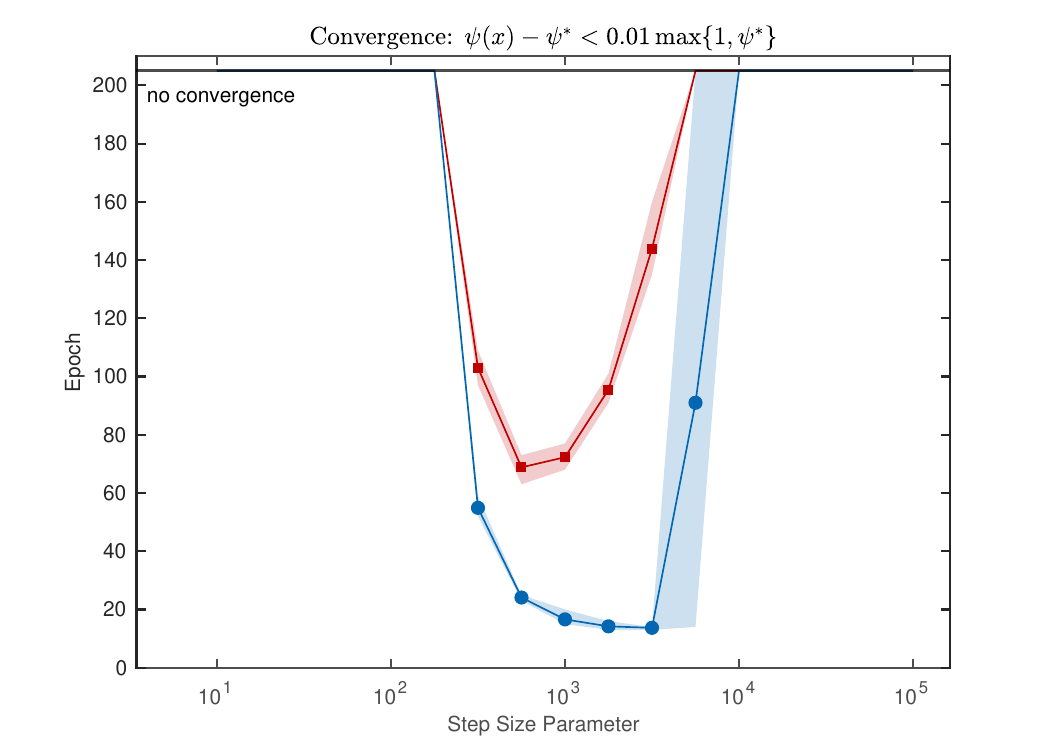}};
	\node at (12.55,0) {\rotatebox{-90}{{\scriptsize $|S_k|=16$}}};
	\node at (12.55,-3.3) {\rotatebox{-90}{{\scriptsize $|S_k|=64$}}};
    \node at (12.55,-6.6) {\rotatebox{-90}{{\scriptsize $|S_k|=256$}}};
    
	\node[right] at (0.0,-6.4) {\includegraphics[height=3.2cm,trim=20 0 20 0,clip]{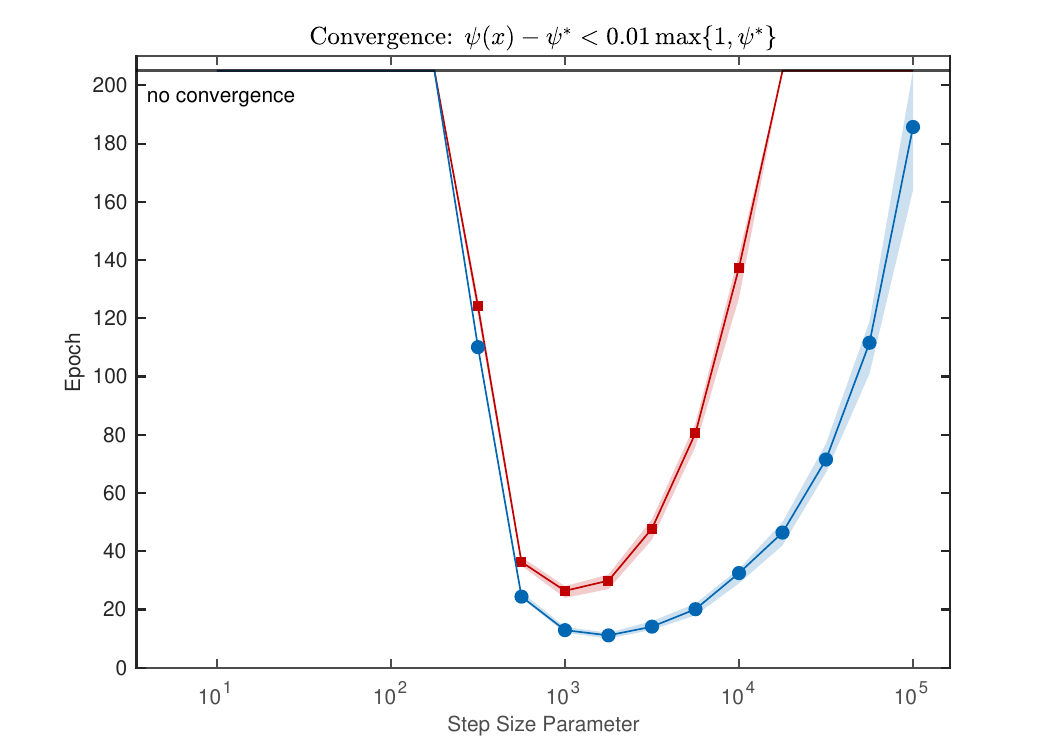}};
	\node[right] at (4.25,-6.4) {\includegraphics[height=3.2cm,trim=40 0 20 0,clip]{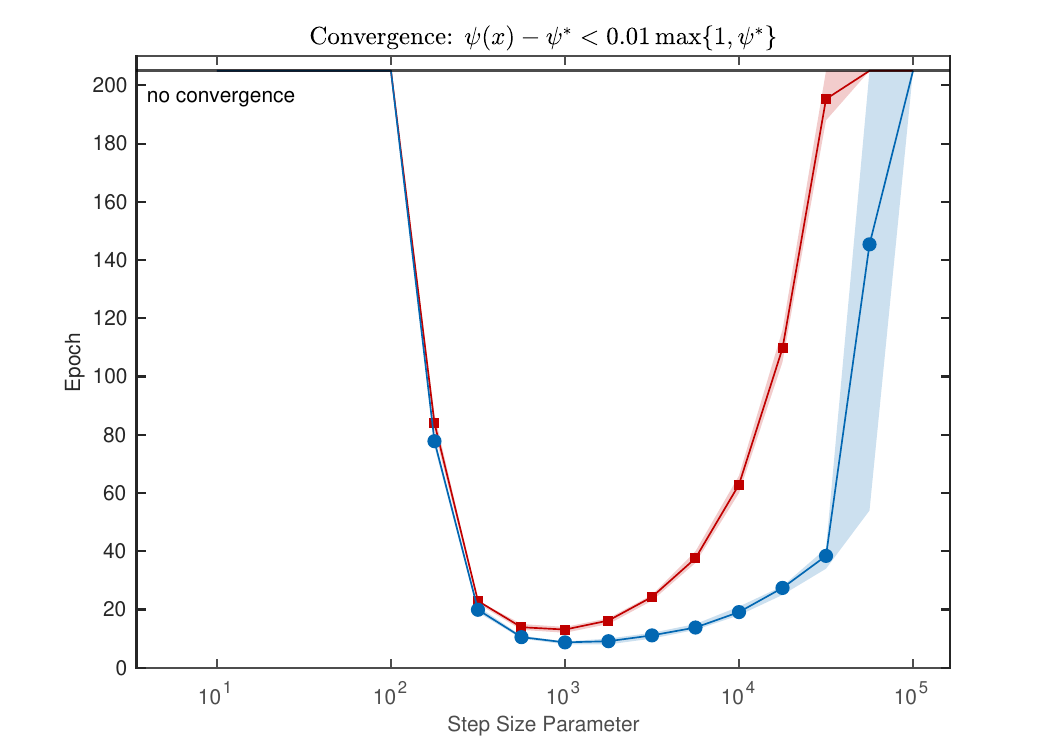}};
	\node[right] at (8.3,-6.4) {\includegraphics[height=3.2cm,trim=40 0 20 0,clip]{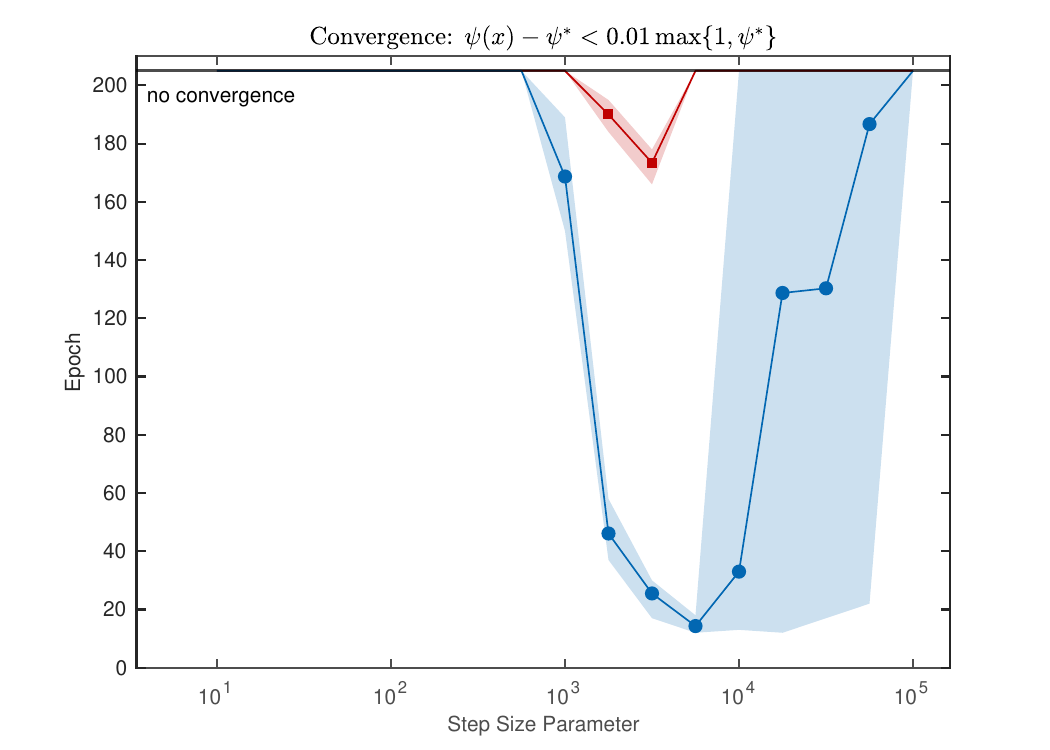}};

    \node[right] at (1.4,-8.4) {{\footnotesize(a)~\texttt{news20}}};
	\node[right] at (5.5,-8.4) {{\footnotesize(b)~\texttt{rcv1}}};
	\node[right] at (9.5,-8.4) {{\footnotesize(c)~\texttt{gisette}}};
\end{tikzpicture}
\vspace{1ex}
	\caption{
    \revise{Convergence behavior of $\NSGD$ and $\PSGD$ for the classification problem \cref{eq:binary-clas} using different step size ranges and mini-batches. (Averaged over $10$ runs).}
    }
	\label{fig:exp2}
\end{figure}

In \cref{fig:exp2}, we depict the number of epochs required to satisfy the relative accuracy criterion ``$\psi(x^k)-\psi^* \leq 0.01\max\{1,\psi^*\}$'' vs$.$ different step size parameters \revise{$\alpha \in [10^0,10^5]$.} \revise{We report ``no convergence'' if this criterion is not met within $200$ epochs.} The average performance (over $10$ independent runs) is shown using a darker color and thicker line style. Here, one single epoch corresponds to \revise{$\lceil N/|S_k| \rceil$} iterations of $\PSGD$ or $\NSGD$ and the approximate optimal value $\psi^*$ is obtained by running the deterministic proximal gradient method starting from three different initial points until $\|F_{\mathrm{nat}}^1(x^k)\|<10^{-12}$. \revise{The results shown in \cref{fig:exp2}} indicate that $\NSGD$ is more robust with respect to the choice of the parameter $\alpha$ and that it can generally achieve faster convergence.

\begin{figure}[t]
\centering
	\setlength{\abovecaptionskip}{-3pt plus 3pt minus 0pt}
	\setlength{\belowcaptionskip}{-10pt plus 3pt minus 0pt}
	\centering
	\begin{tikzpicture}[scale=1]
    \node[right] at (3.44,2) {\footnotesize $\NSGD$};
    \node[right] at (9.24,2) {\footnotesize $\PSGD$};
    \draw[draw=MyGray] (0.54,1.75) rectangle (12.14,2.25);
    \draw[draw=DeepBlue,fill=DeepBlue,line width=1pt] (2.34,2) -- (3.24,2);
    \node[circle,draw=DeepBlue,fill=DeepBlue,minimum size=4pt,inner sep=0pt, outer sep=0pt] at (2.79,2) {};
    \draw[draw=DeepRed,fill=DeepRed,line width=1pt] (8.14,2) -- (9.04,2);
    \node[draw=DeepRed,fill=DeepRed,minimum size=4pt,inner sep=0pt, outer sep=0pt] at (8.59,2) {};

	\node[right] at (0.0,0) {\includegraphics[height=3.2cm,trim=20 0 20 0,clip]{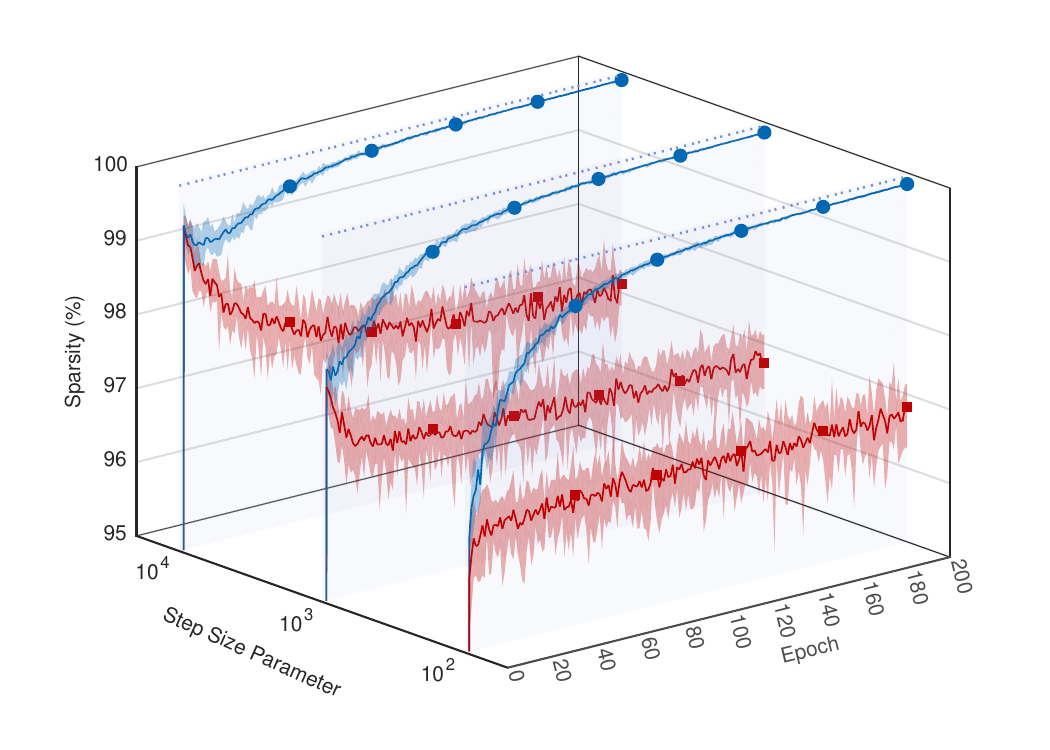}};
	\node[right] at (4.25,0) {\includegraphics[height=3.2cm,trim=40 0 20 0,clip]{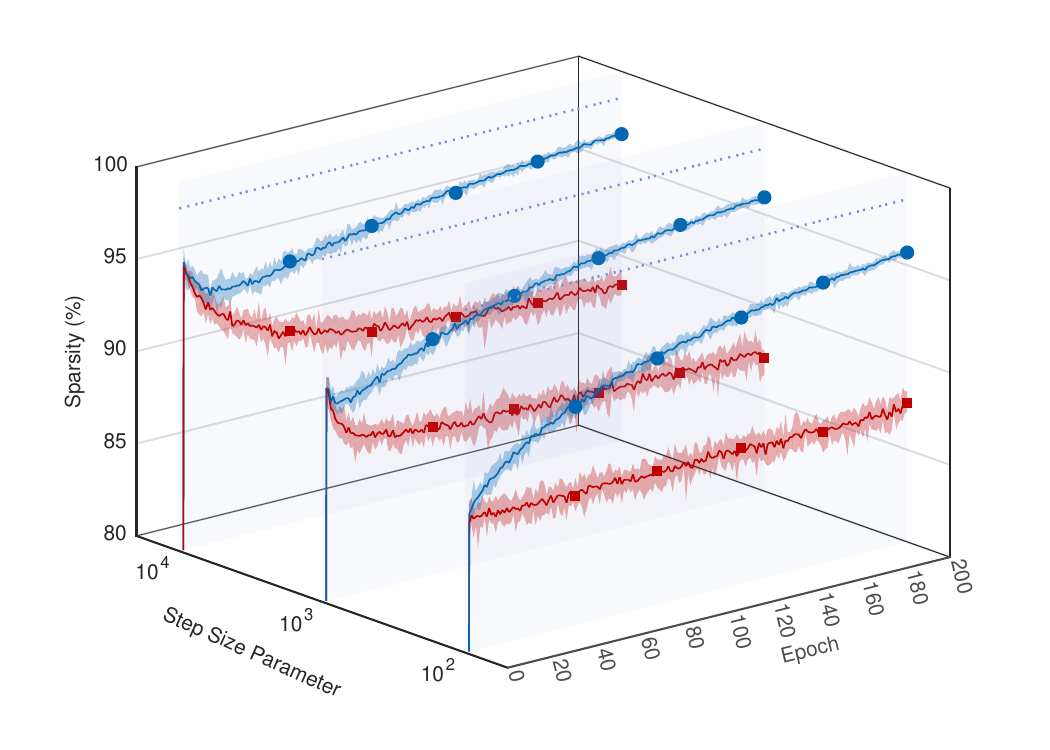}};
	\node[right] at (8.3,0) {\includegraphics[height=3.2cm,trim=40 0 20 0,clip]{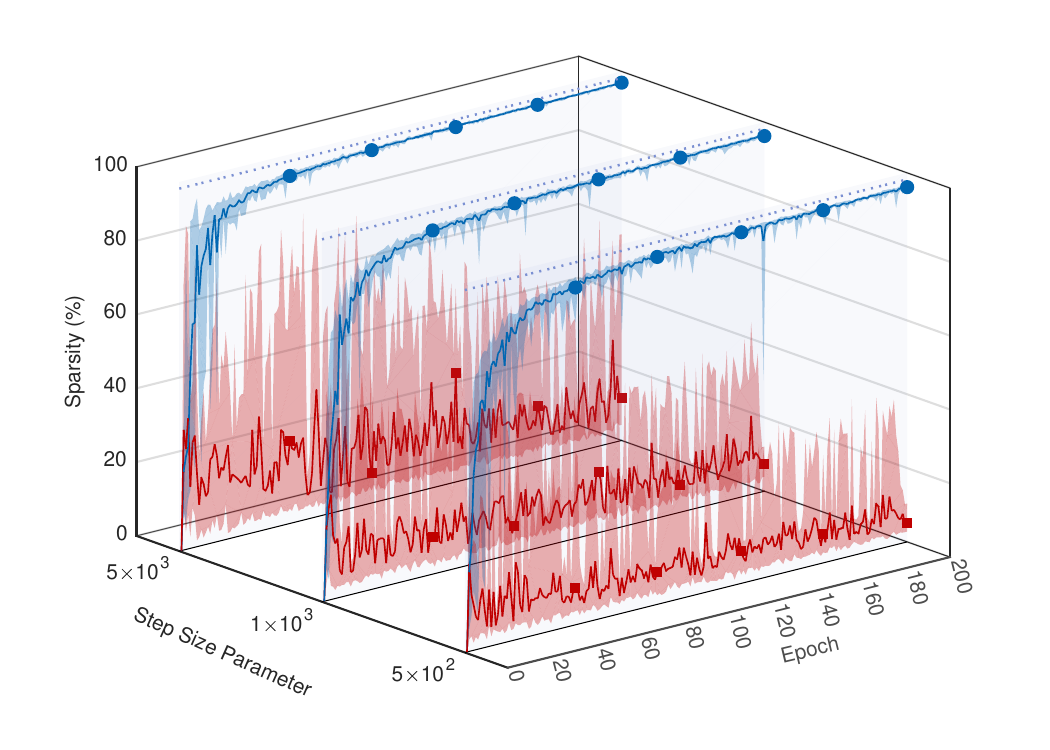}};
	\node[right] at (0.0,-3.2) {\includegraphics[height=3.2cm,trim=20 0 20 0,clip]{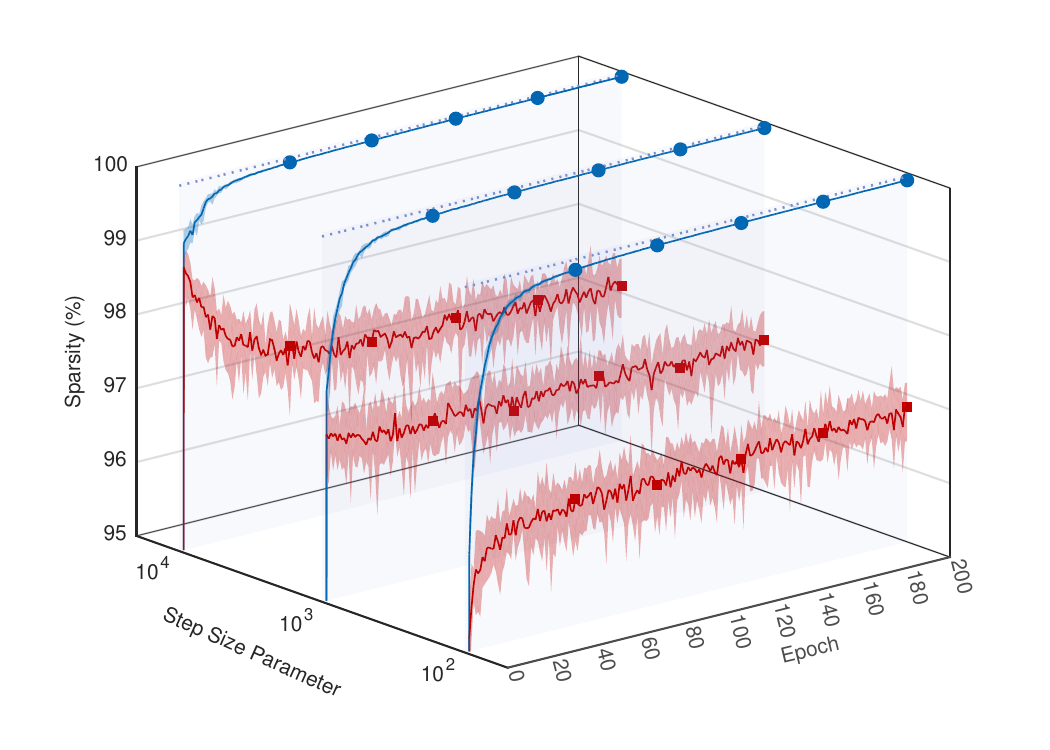}};
	\node[right] at (4.25,-3.2) {\includegraphics[height=3.2cm,trim=40 0 20 0,clip]{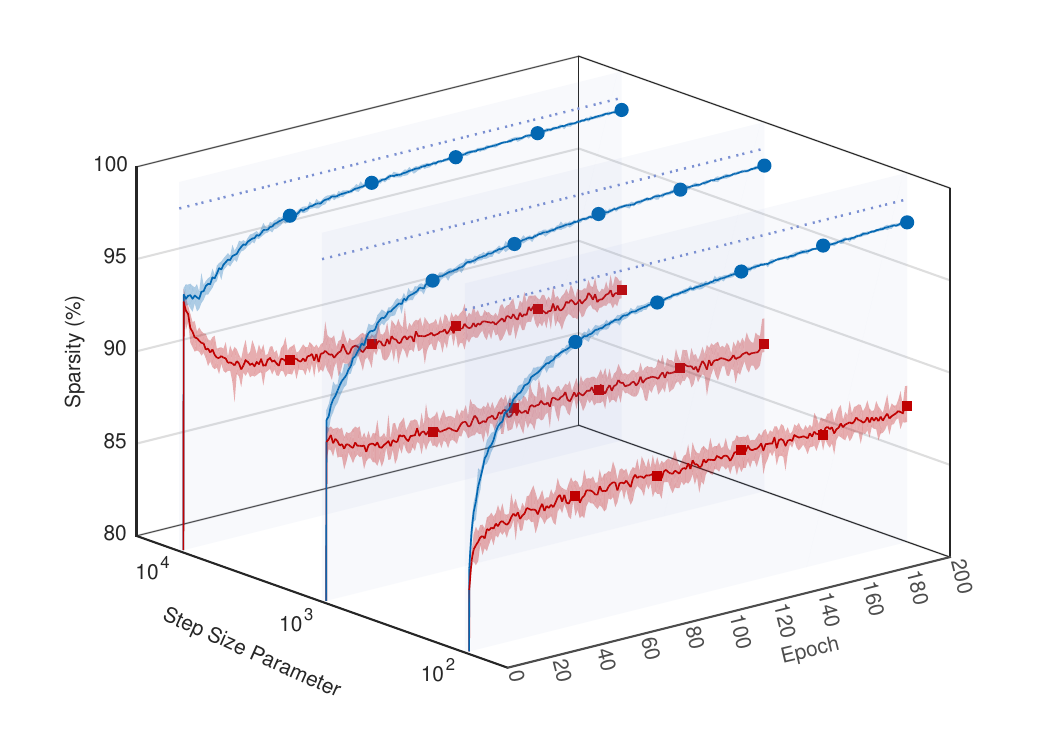}};
	\node[right] at (8.3,-3.2) {\includegraphics[height=3.2cm,trim=40 0 20 0,clip]{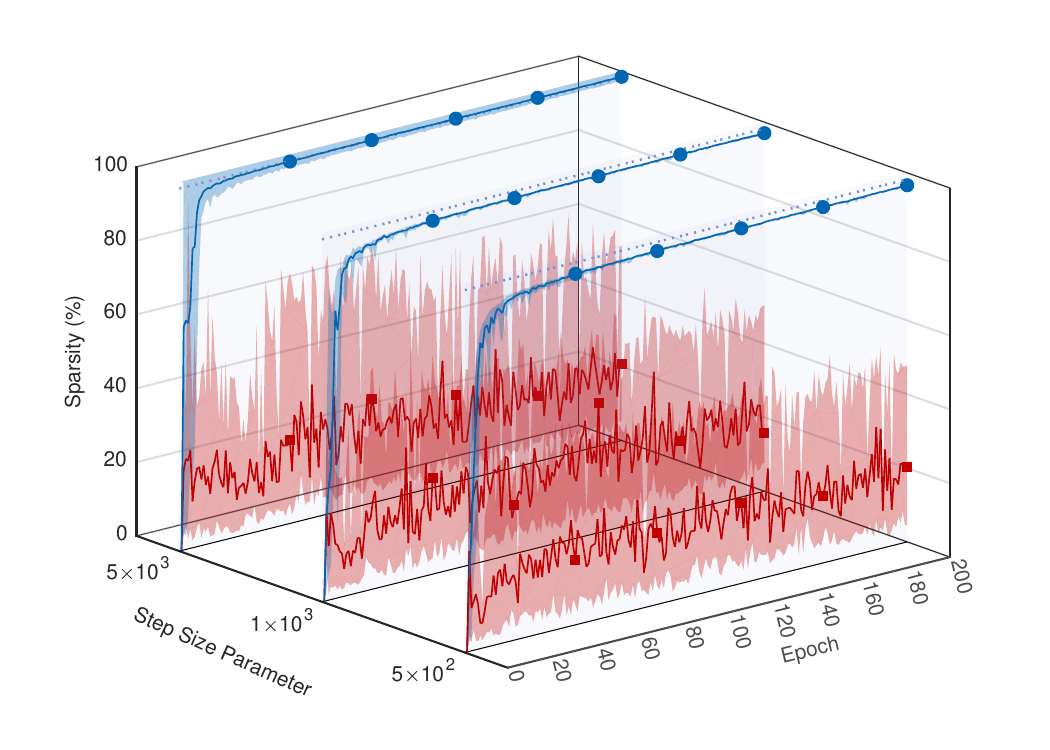}};
	\node at (12.55,0) {\rotatebox{-90}{{\scriptsize $|S_k|=16$}}};
	\node at (12.55,-3.3) {\rotatebox{-90}{{\scriptsize $|S_k|=64$}}};
    \node at (12.55,-6.6) {\rotatebox{-90}{{\scriptsize $|S_k|=256$}}};
    
	\node[right] at (0.0,-6.4) {\includegraphics[height=3.2cm,trim=20 0 20 0,clip]{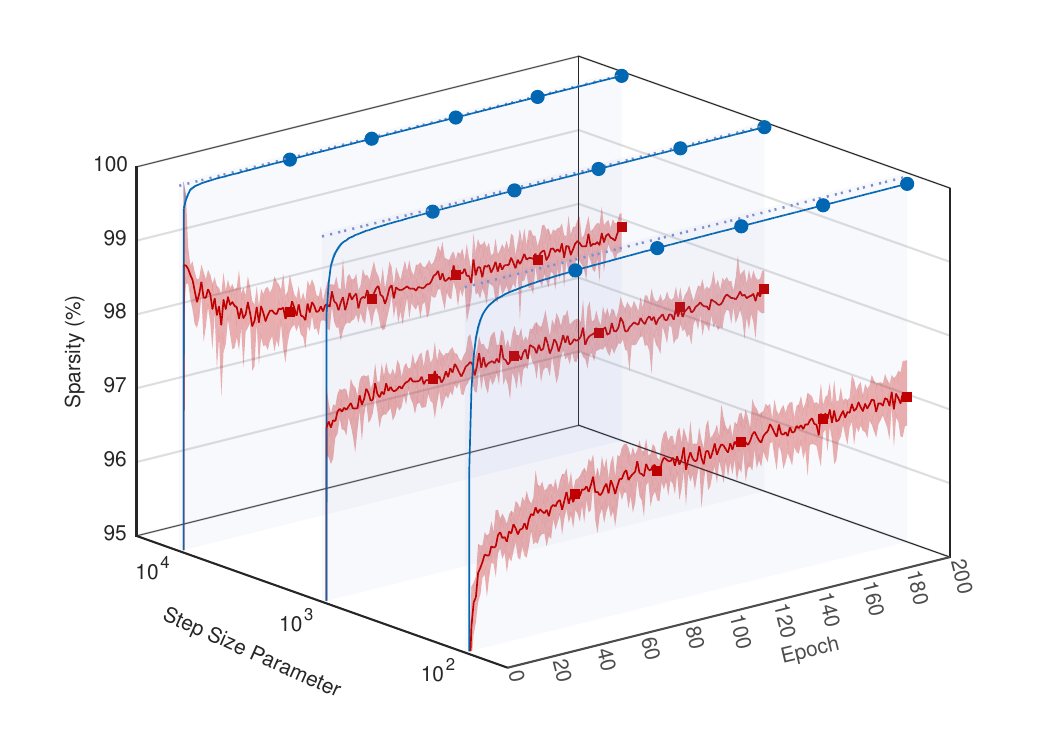}};
	\node[right] at (4.25,-6.4) {\includegraphics[height=3.2cm,trim=40 0 20 0,clip]{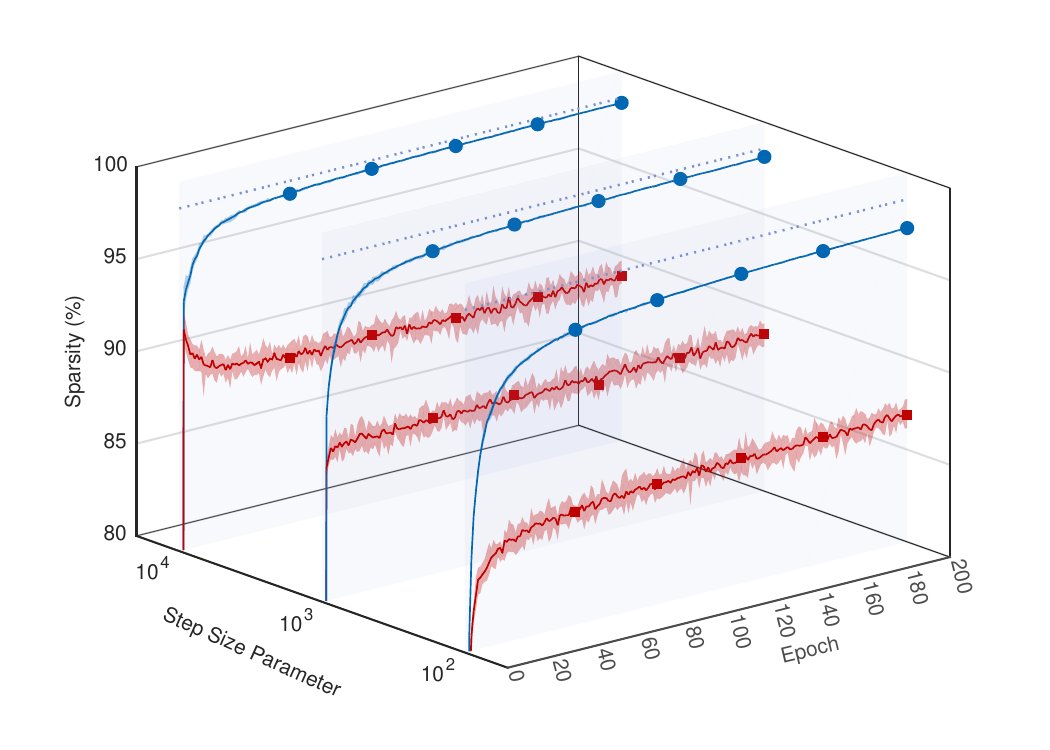}};
	\node[right] at (8.3,-6.4) {\includegraphics[height=3.2cm,trim=40 0 20 0,clip]{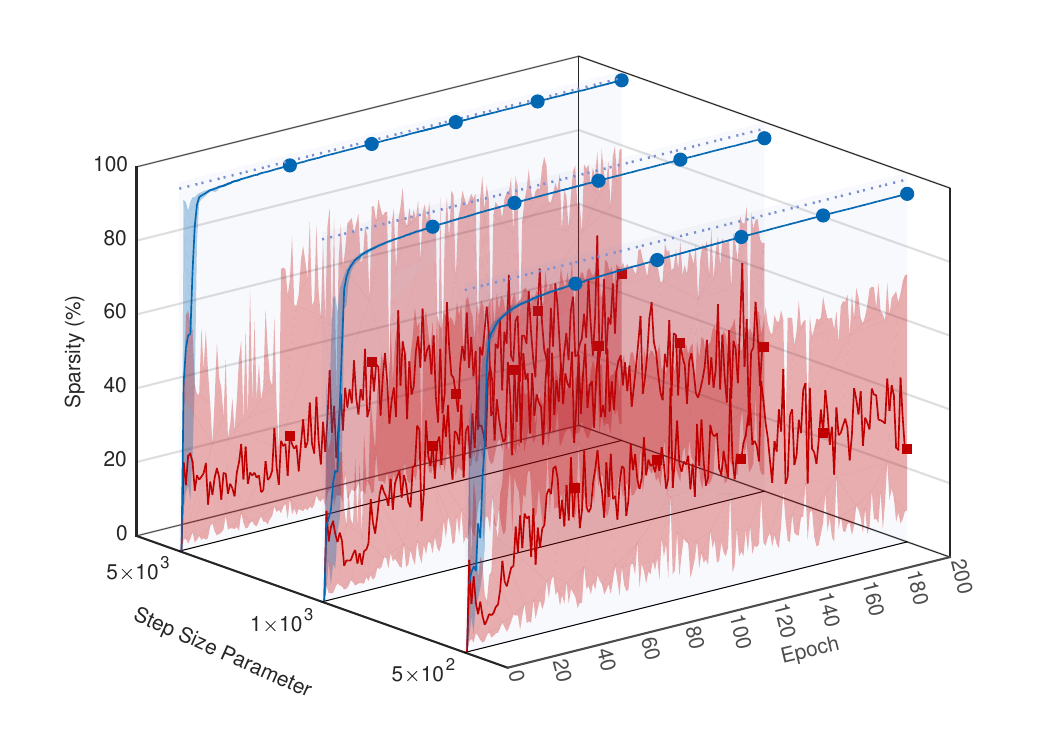}};

    \node[right] at (1.4,-8.4) {{\footnotesize(a)~\texttt{news20}}};
	\node[right] at (5.5,-8.4) {{\footnotesize(b)~\texttt{rcv1}}};
	\node[right] at (9.5,-8.4) {{\footnotesize(c)~\texttt{gisette}}};
\end{tikzpicture}
\vspace{1ex}
	\caption{
    \revise{Sparsity information of $\NSGD$ and $\PSGD$ for the classification problem \cref{eq:binary-clas}. (Averaged over $10$ runs). The sparsity is measured via $100\% \cdot |\{i: x_i^k=0\}|/d$.}
    }
	\label{fig:exp2-sparse}
\end{figure}

\revise{In \cref{fig:exp2-sparse},} we plot the \revise{sparsity level $100\% \cdot |\{i: x_i^k=0\}|/d$} of the iterates $\revise{\{x^k\}_{k\in\N}}$ vs$.$ the number of epochs (again averaged over 10 runs). The reference sparsity levels (obtained by $\PGD$) for \texttt{news20}, \texttt{rcv1}, and \texttt{gisette} are $99.95\%$, $98.57\%$, and $98.23\%$\revise{, respectively}. Overall, $\NSGD$ tends to recover sparser solutions than $\PSGD$ which reflects the different scalings of the $\ell_1$-norm in the proximity operators and the better identification properties of $\NSGD$.

\subsection{\texorpdfstring{Sparse\,$+$\,low-rank matrix decomposition}{Sparse + low-rank matrix decomposition}}
We now study the performance of $\NSGD$, $\PSGD$, \revise{and $\RDA$} on a principal component pursuit task, \cite{candes2011pcp}. In particular, given a matrix $M \in \Rmn$, we want to recover $X,Y\in\Rmn$ that minimize the following ``sparse\,$+$\,low-rank'' model
\begin{equation} \label{eq:prob-pcp}
{\min}_{X,Y\in\Rmn} \; f(X,Y) + \nu_1 \|X\|_* + \nu_2\|Y\|_1, \quad \nu_1, \nu_2 > 0,
\end{equation}
where $f(X,Y) := \frac{1}{2} \|X+Y-M\|^2_F$ and $\|\cdot\|_*$ denotes the nuclear norm. 
According to \cite{vaiter2017model}, the nuclear norm is partly smooth at $\bar X \in \R^{m \times n}$ relative to the constant rank manifold $\mathcal M_{\bar X} := \{X: \mathrm{rank}(X) = \mathrm{rank}(\bar X)\}$. Thus, we can again expect $\NSGD$ to identify the underlying low-rank structure and sparsity when solving \cref{eq:prob-pcp}.

We consider the special case where $M$ is a video clip with each column $M_i$ of $M$ corresponding to a vectorized frame. The model \cref{eq:prob-pcp} aims to decompose $M$ into a low-rank background $X$ and into a sparse component $Y$ representing movements or moving objects. This application is known as video background subtraction, \cite{candes2011pcp,ji2011video}. We use the stochastic gradients $g^k := \frac{1}{|S_k|}\sum_{i \in S_k} \nabla f_i(X^k,Y^k)$, $f_i(X,Y) := \frac{n}2\|X_i+Y_i-M_i\|^2$, $i\in[n]$, and $S_k$ is a mini-batch selected \revise{uniformly at random} without replacement from $[n]$. Hence, each stochastic gradient only accesses a total of $|S_k|$ frames of $M$ at iteration $k$. 

\begin{figure}[t]
\centering
	\setlength{\abovecaptionskip}{-3pt plus 3pt minus 0pt}
	\setlength{\belowcaptionskip}{-10pt plus 3pt minus 0pt}
	\centering
\tikzmath{\h1=-2.9;}
\begin{tikzpicture}[scale=1]
    \node[right] at (2.5,1.9) {\footnotesize $\NSGD$, $\lambda = 1$};
    \node[right] at (6.48,1.9) {\footnotesize $\NSGD$, $\lambda = 2$};
    \node[right] at (10.46,1.9) {\footnotesize $\PSGD$};
    \draw[draw=MyGray] (0.51,1.65) rectangle (12.47,2.15);
    \draw[draw=DeepBlue,fill=DeepBlue,line width=1pt] (1.5,1.9) -- (2.4,1.9);
    \node[circle,draw=DeepBlue,fill=DeepBlue,minimum size=4pt,inner sep=0pt, outer sep=0pt] at (1.95,1.9) {};
    \draw[draw=MyBlue,fill=MyBlue,line width=1pt] (5.48,1.9) -- (6.38,1.9);
    \node[diamond,draw=MyBlue,fill=MyBlue,minimum size=4.5pt,inner sep=0pt, outer sep=0pt] at (5.93,1.9) {};
    \draw[draw=DeepRed,fill=DeepRed,line width=1pt] (9.46,1.9) -- (10.36,1.9);
    \node[draw=DeepRed,fill=DeepRed,minimum size=4pt,inner sep=0pt, outer sep=0pt] at (9.91,1.9) {};
    \node[right] at (0,0) {\includegraphics[height=2.6cm]{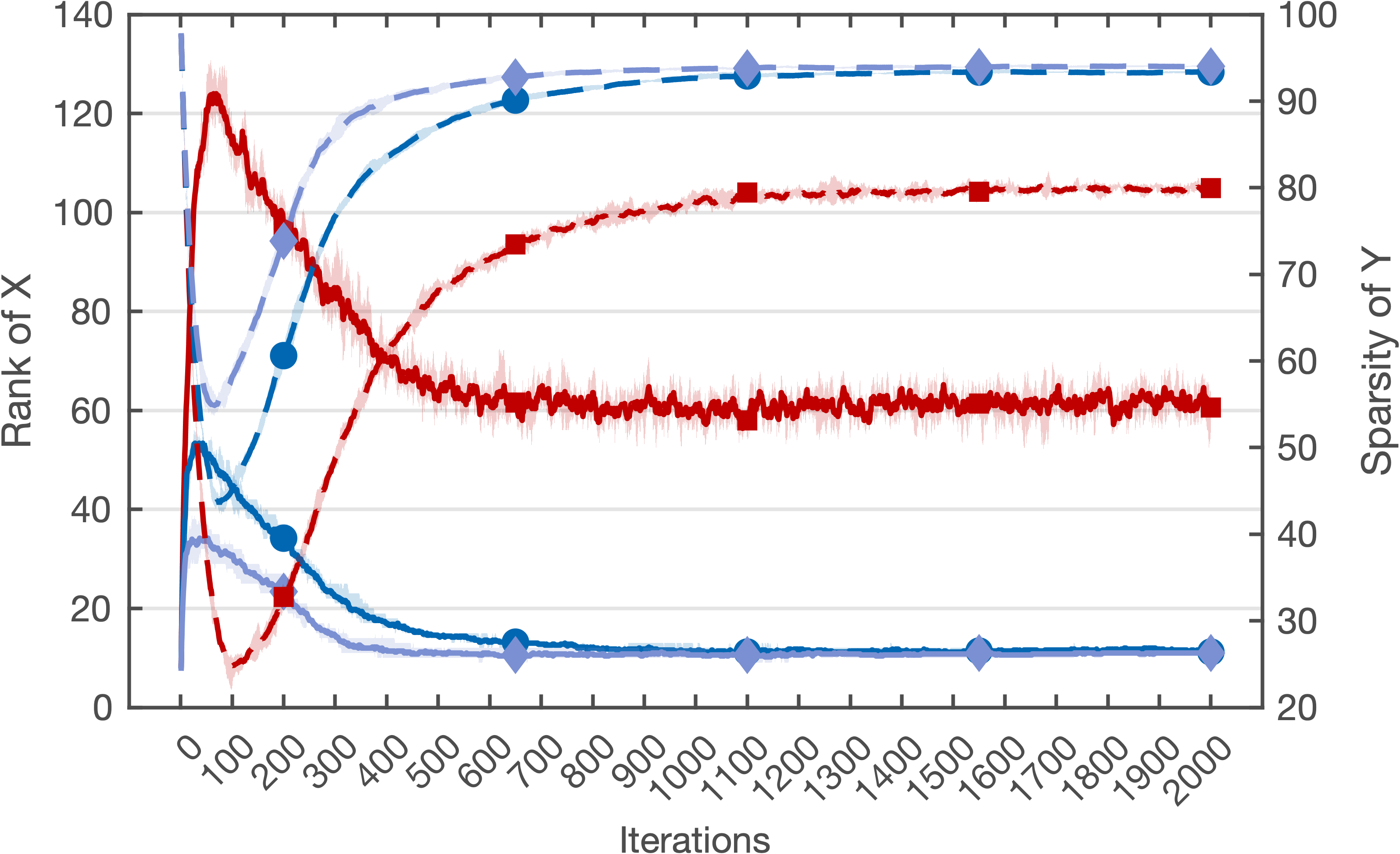}};
    \node[right] at (4.45,0) {\includegraphics[height=2.6cm]{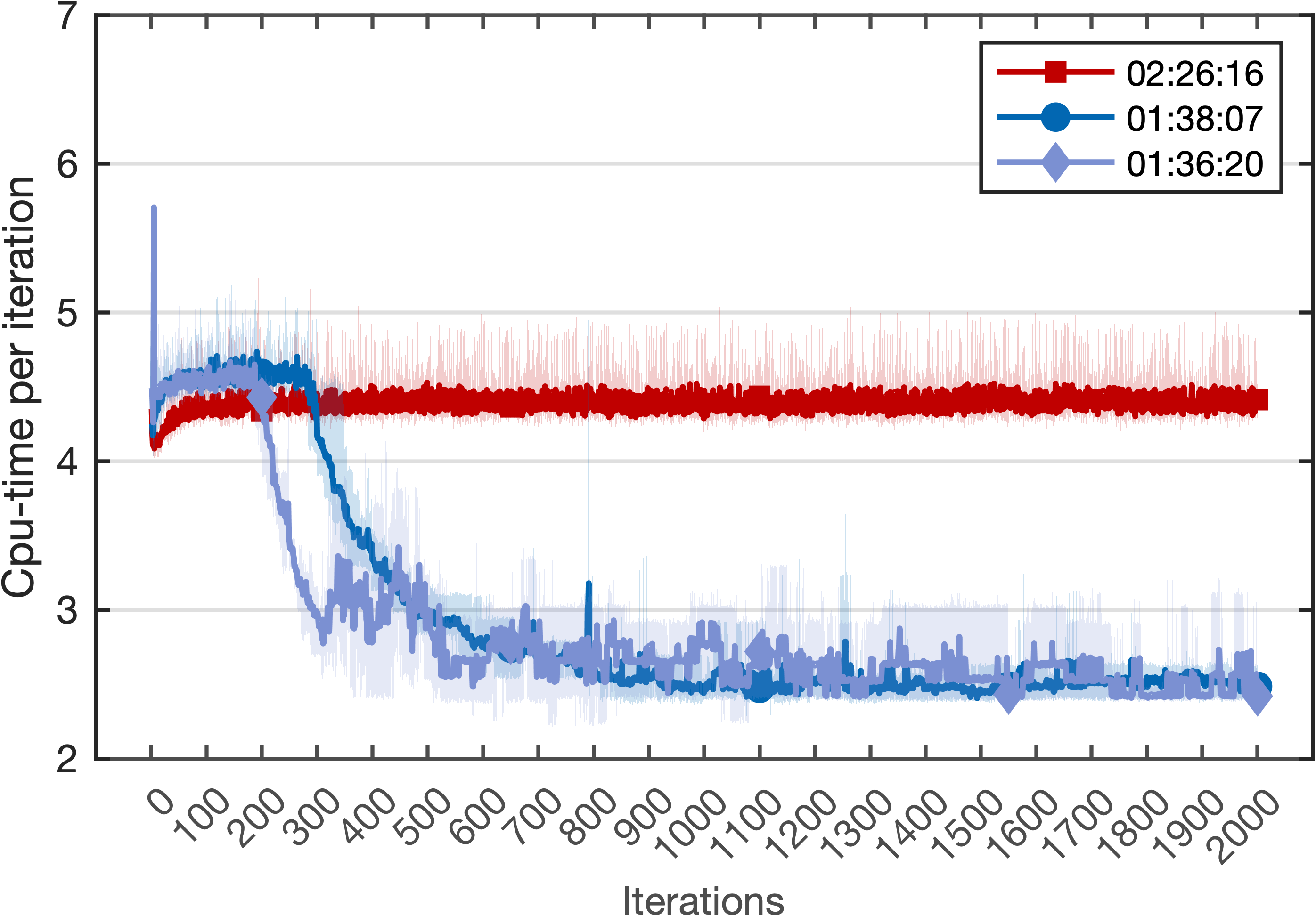}};
    \node[right] at (8.4,0) {\includegraphics[height=2.6cm]{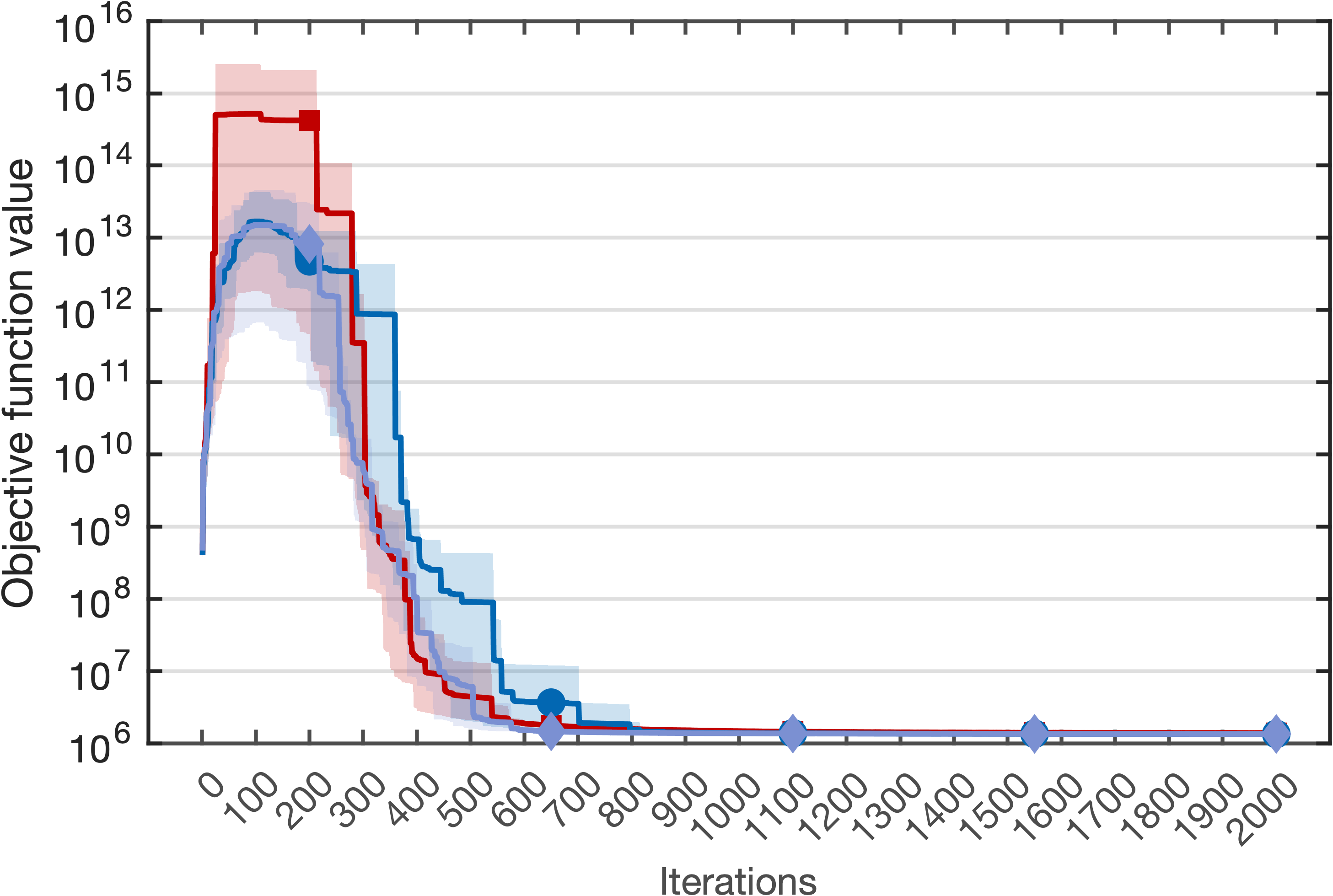}};
    \draw[draw=gray,densely dotted] (11.76,-0.73) -- (10.3,0);
    \draw[draw=gray,densely dotted] (12.2,-0.73) -- (12.3,0);
    \draw[draw=gray] (11.43,-0.89) rectangle (12.29,-0.73); 
    \node[right,draw=gray,fill=white,inner sep=0.5mm] at (10,0.6) {\includegraphics[height=1.6cm]{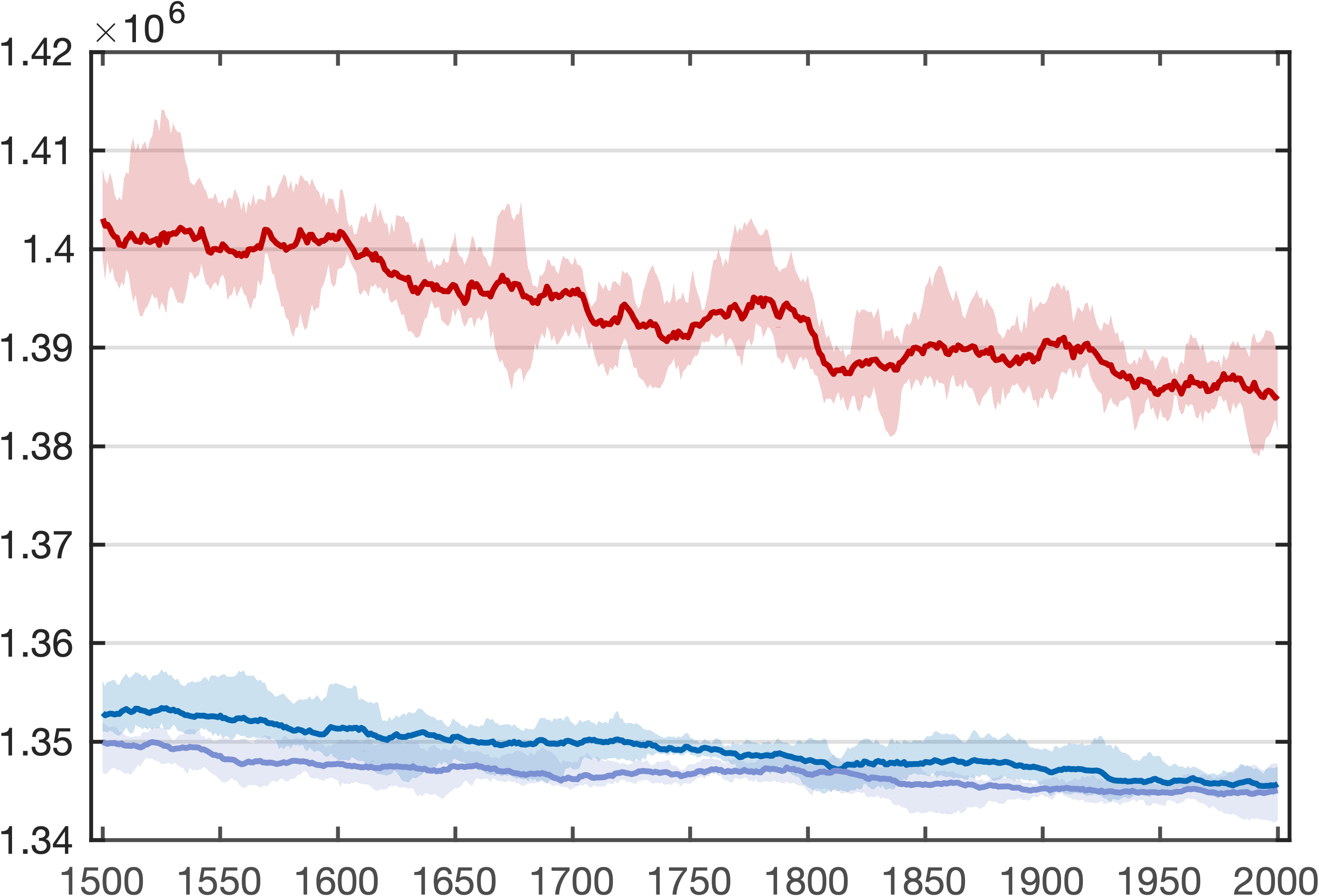}};
    \node at (12.8,0) {\rotatebox{-90}{{\scriptsize $\gamma = 2/3$}}};
    \node[right] at (0,\h1) {\includegraphics[height=2.6cm]{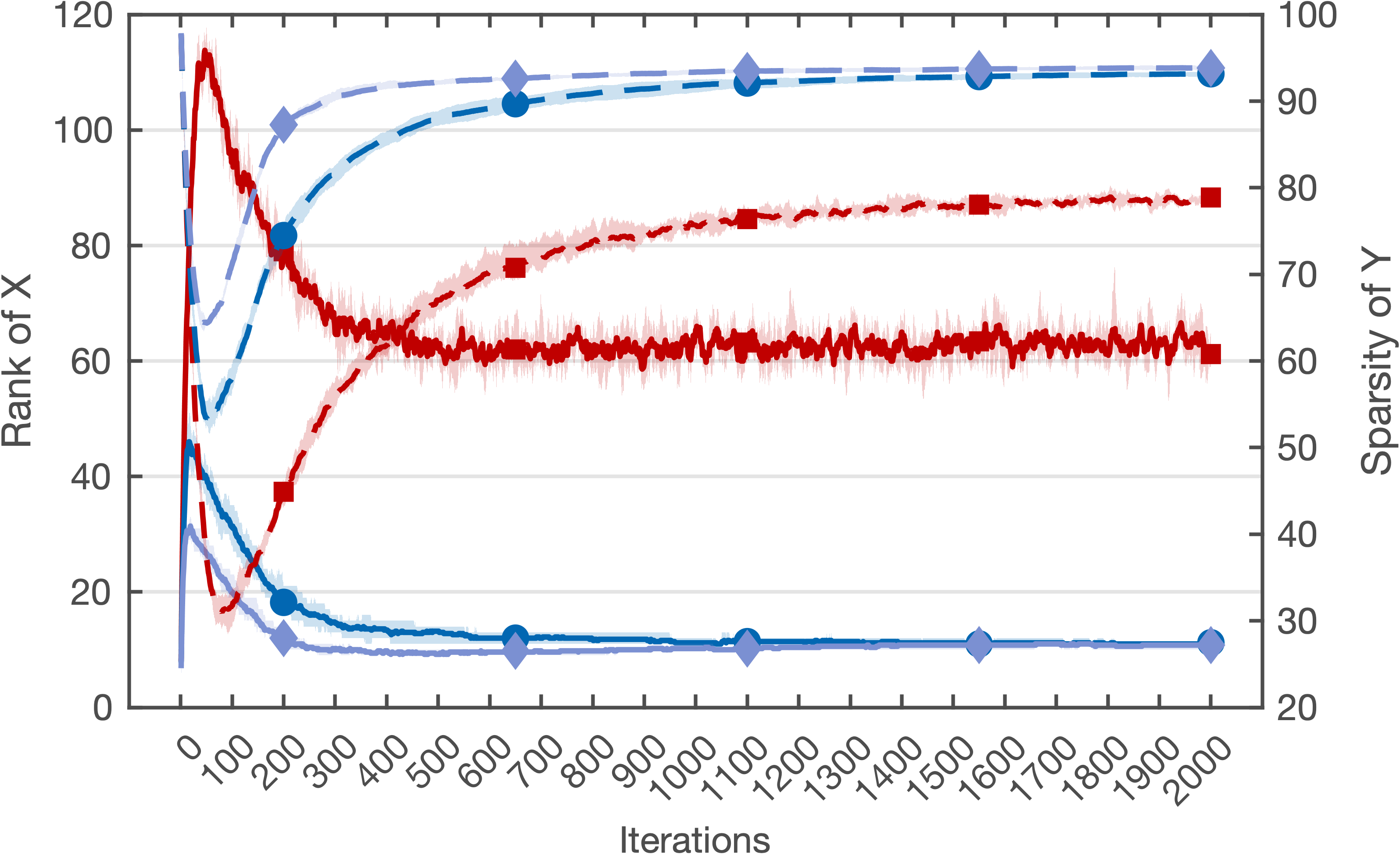}};
    \node[right] at (4.45,\h1) {\includegraphics[height=2.6cm]{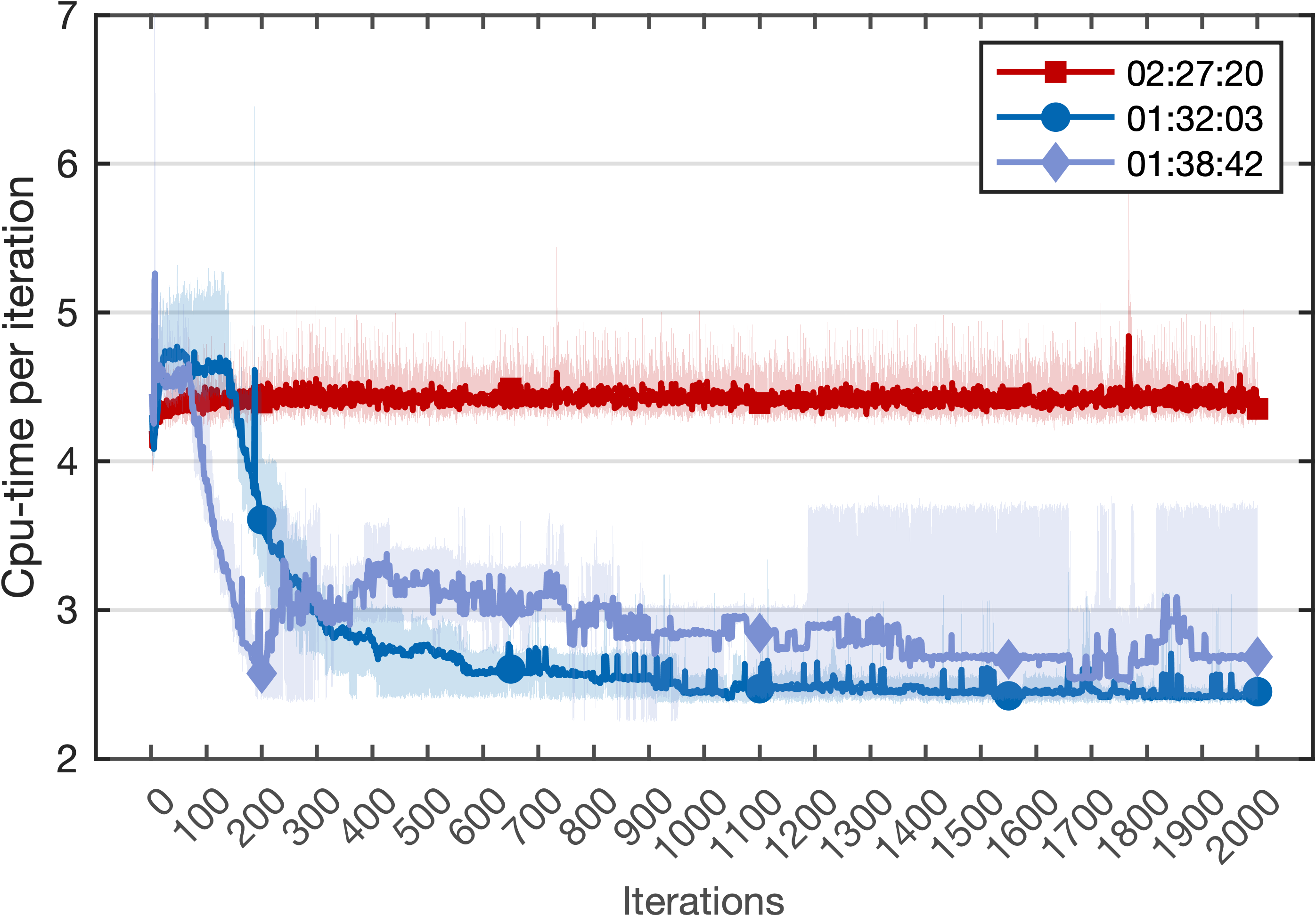}};
    \node[right] at (8.4,\h1) {\includegraphics[height=2.6cm]{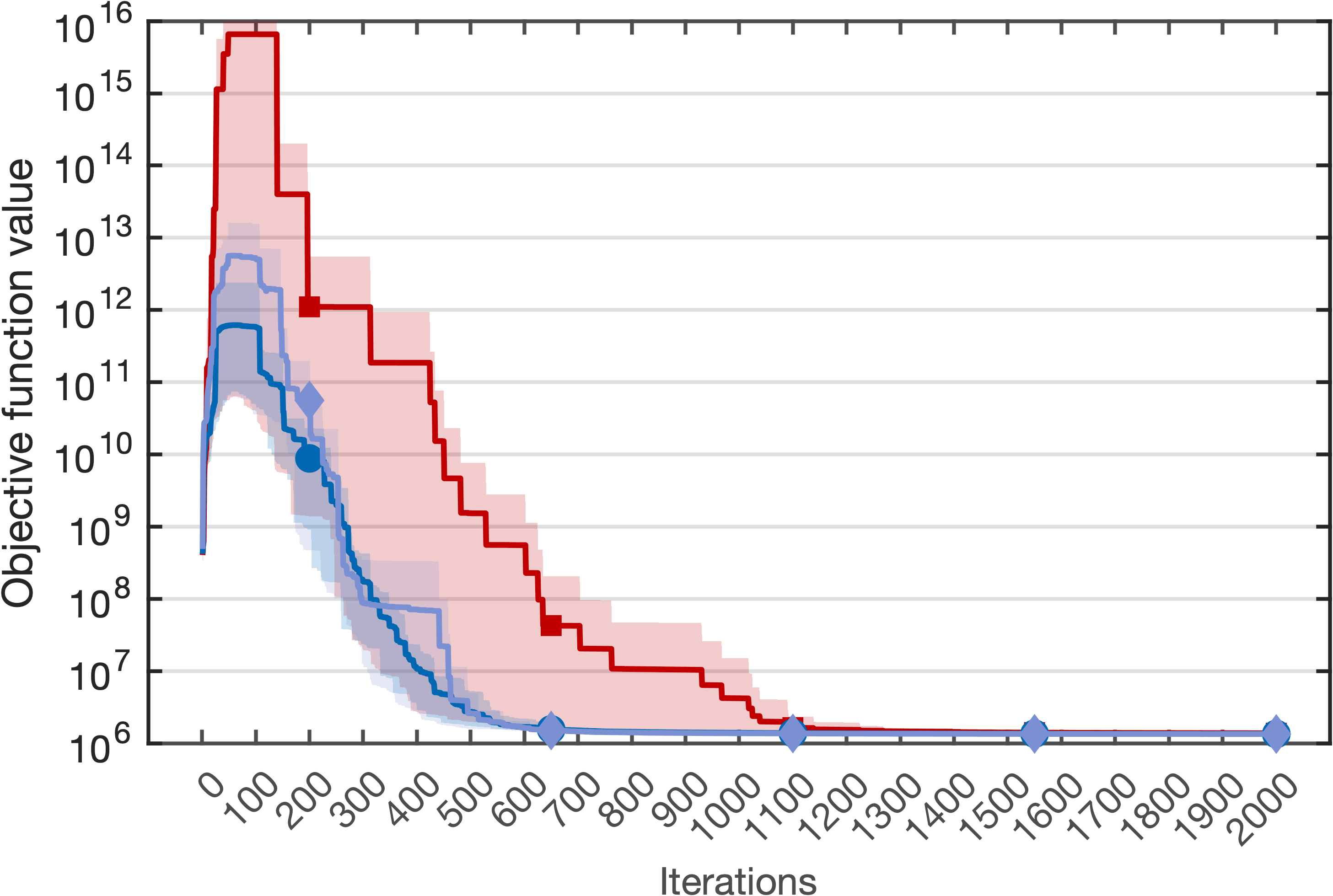}};
    \draw[draw=gray,densely dotted] (11.76,-0.73+\h1) -- (10.3,\h1);
    \draw[draw=gray,densely dotted] (12.2,-0.73+\h1) -- (12.3,\h1);
    \draw[draw=gray] (11.43,-0.89+\h1) rectangle (12.29,-0.73+\h1); 
    \node[right,draw=gray,fill=white,inner sep=0.5mm] at (10,0.6+\h1) {\includegraphics[height=1.6cm]{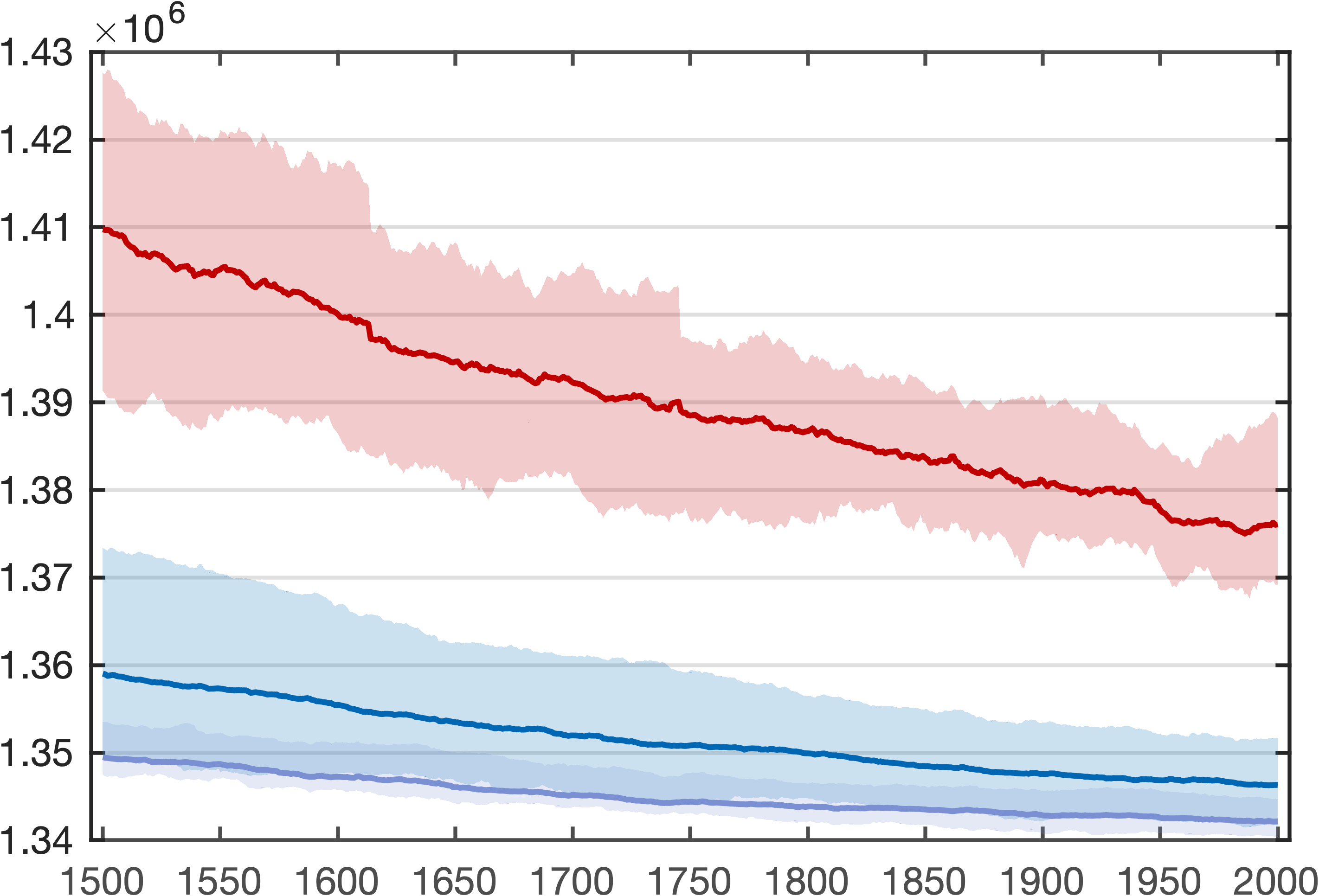}};
    \node at (12.8,\h1) {\rotatebox{-90}{{\scriptsize $\gamma = 3/4$}}};
    \node[right] at (0,2*\h1) {\includegraphics[height=2.6cm]{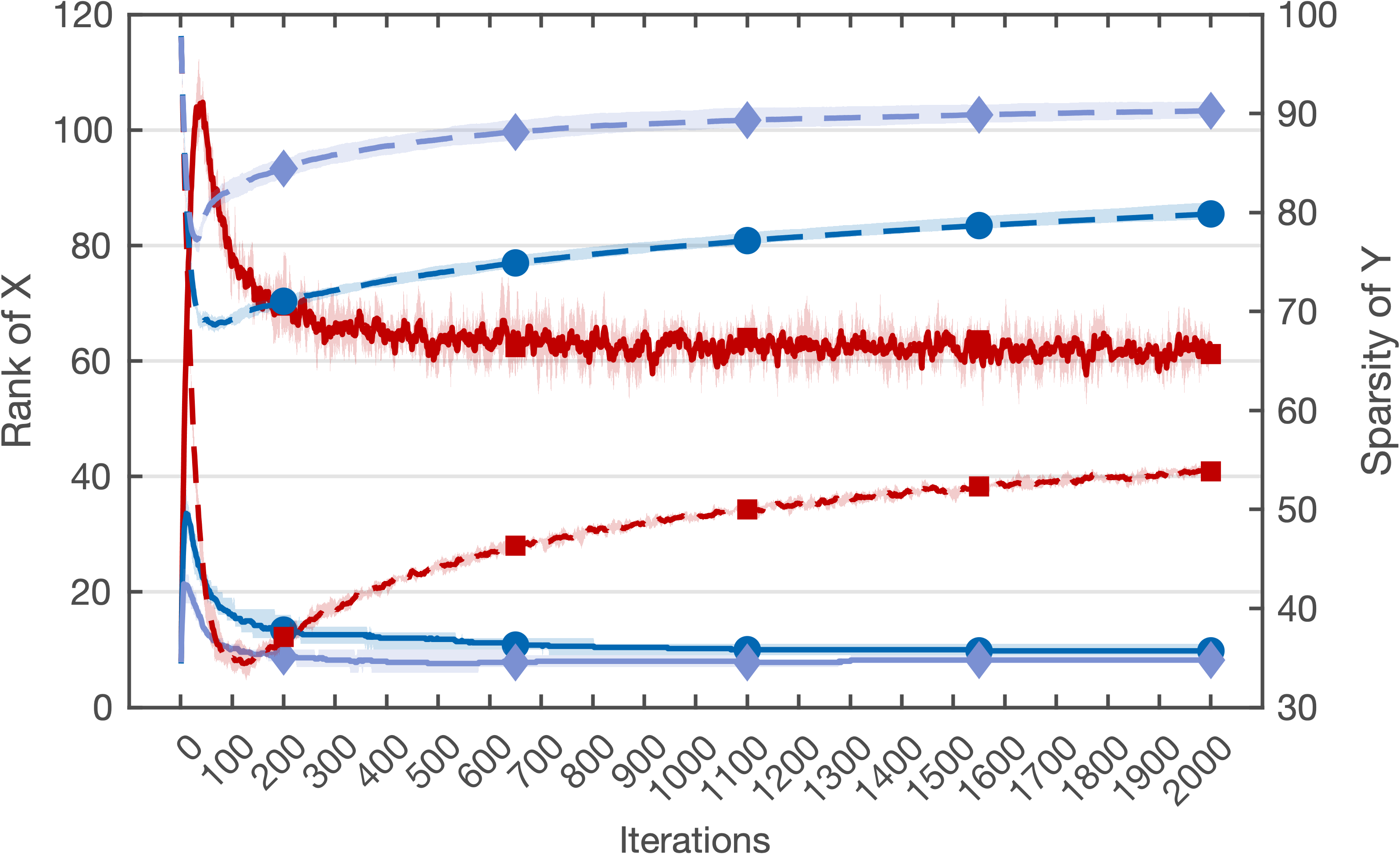}};
	\node[right] at (4.45,2*\h1) {\includegraphics[height=2.6cm]{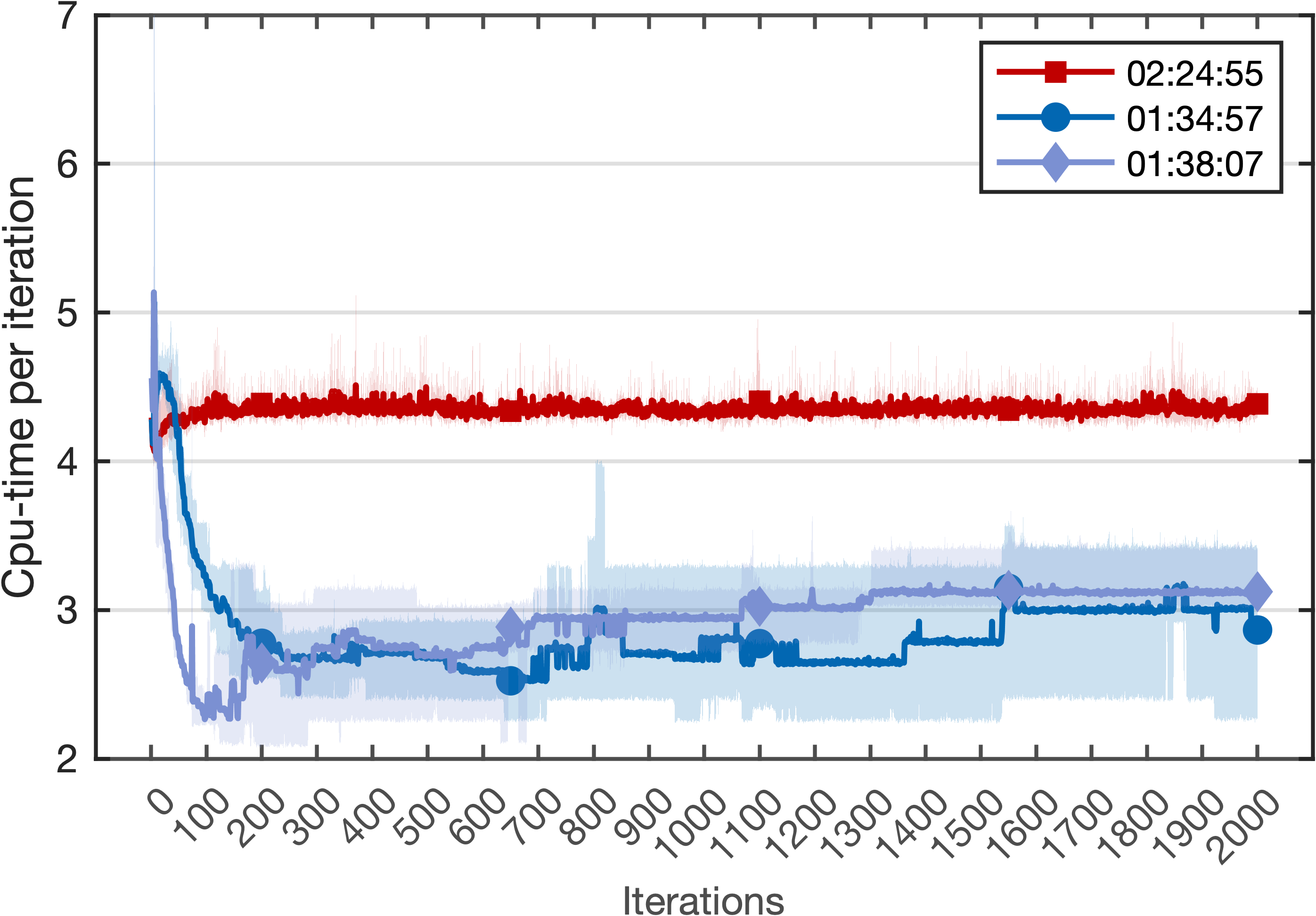}};
    \node[right] at (8.4,2*\h1) {\includegraphics[height=2.6cm]{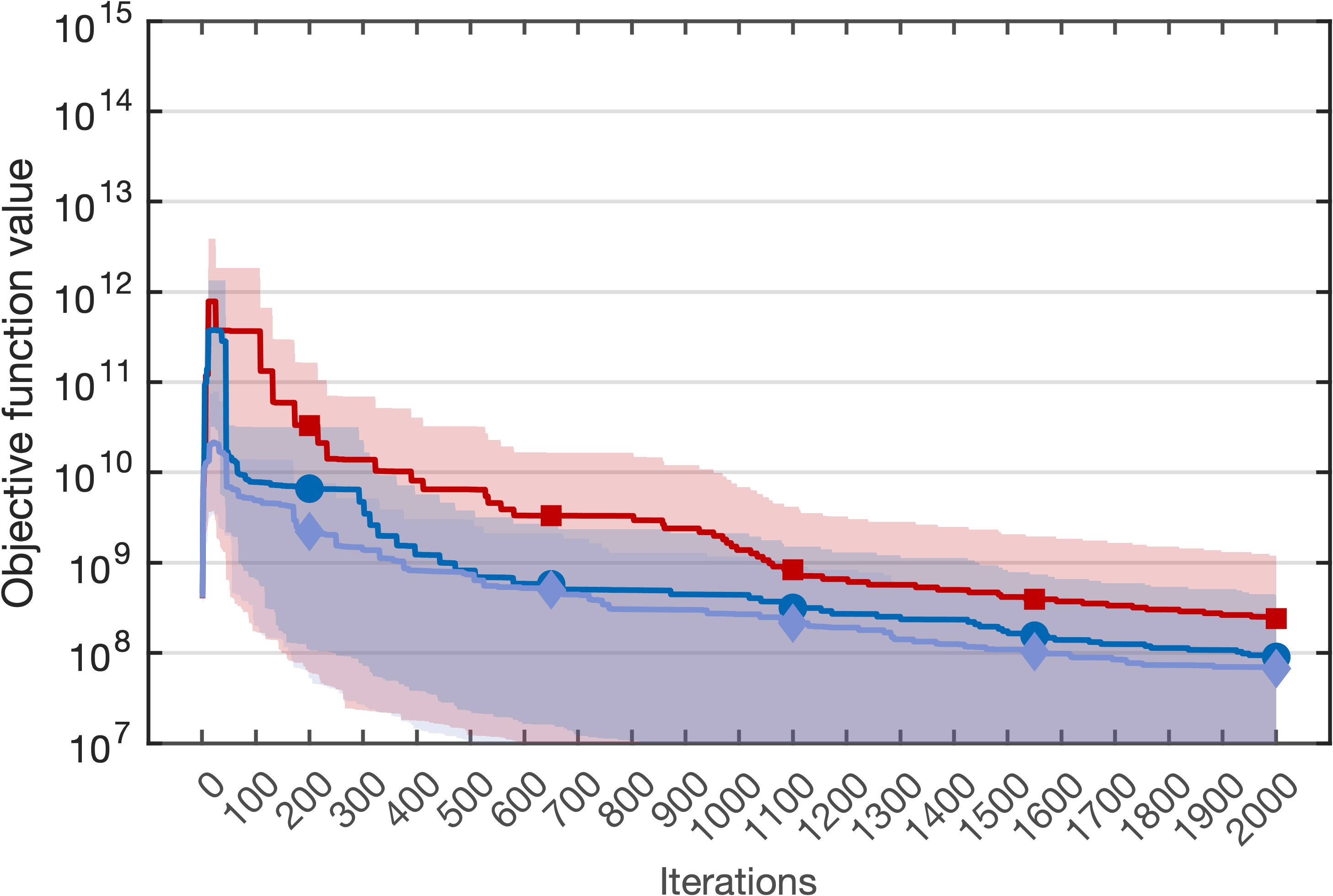}};
    \draw[draw=gray,densely dotted] (11.76,-0.39+2*\h1) -- (10.3,2*\h1);
    \draw[draw=gray,densely dotted] (12.2,-0.39+2*\h1) -- (12.3,2*\h1);
    \draw[draw=gray] (11.43,-0.7+2*\h1) rectangle (12.29,-0.39+2*\h1); 
    \node[right,draw=gray,fill=white,inner sep=0.5mm] at (10.07,0.6+2*\h1) {\includegraphics[height=1.6cm]{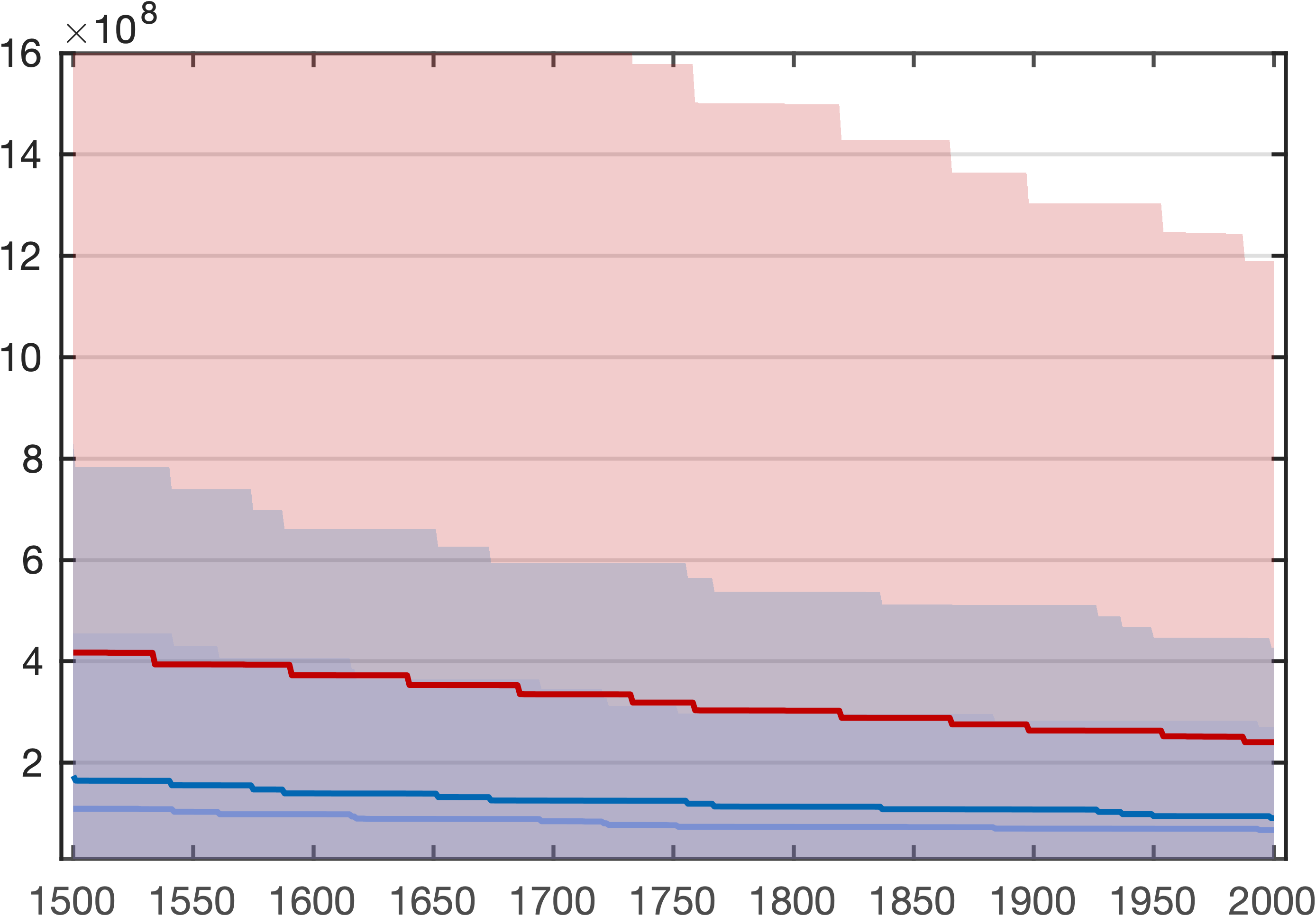}};
    \node at (12.8,2*\h1) {\rotatebox{-90}{{\scriptsize $\gamma = 1$}}};
\end{tikzpicture}
%
\caption{Performance of $\NSGD$ and $\PSGD$ on the sparse\,$+$\,low-rank task \cref{eq:prob-pcp}. Each row depicts the results for different step sizes $\alpha_k = \frac12(k+1)^{-\gamma}$, $\gamma \in \{\frac23,\frac34,1\}$. In the first column, we plot the rank of \revise{$\{X^k\}_{k\in\N}$} (the number of \revise{singular values larger than $0$;} using solid lines) and the \revise{sparsity level $100\%\cdot|\{i,j:Y_{ij}^k = 0\}|/(mn)$} of $\revise{\{Y^k\}_{k\in\N}}$ (using dashed lines). The cpu-time per iteration is shown in the middle column; the legend depicts the total running time. The right column illustrates the change in the objective function values. (Averaged over $5$ runs).}
\label{fig:exp-3}
\end{figure}

\begin{figure}[t]
\centering
	\setlength{\abovecaptionskip}{-3pt plus 3pt minus 0pt}
	\setlength{\belowcaptionskip}{-10pt plus 3pt minus 0pt}
	\centering
\tikzmath{\h1=-2.9;}
\begin{tikzpicture}[scale=1]
    \node[right] at (1.8,1.9) {\footnotesize $\RDA$, $\gamma = 5$};
    \node[right] at (5.78,1.9) {\footnotesize $\RDA$, $\gamma = 25$};
    \node[right] at (1.8,2.4) {\footnotesize $\NSGD$, $\lambda = 1$};
    \node[right] at (5.78,2.4) {\footnotesize $\NSGD$, $\lambda = 2$};
    \node[right] at (9.76,2.4) {\footnotesize $\NSGD$, $\lambda = 5$};
    \draw[draw=MyGray] (0.51,1.65) rectangle (12.47,2.65);
    \draw[draw=MyOrange,fill=MyOrange,line width=1pt] (0.8,1.9) -- (1.7,1.9);
    \node[isosceles triangle,rotate=90,isosceles triangle apex angle=60,draw=MyOrange,fill=MyOrange,minimum size=4pt,inner sep=0pt, outer sep=0pt] at (1.25,1.9) {};
    \draw[draw=MyTurq,fill=MyTurq,line width=1pt] (4.78,1.9) -- (5.68,1.9);
    \node[isosceles triangle,rotate=270,isosceles triangle apex angle=60,draw=MyTurq,fill=MyTurq,minimum size=4.5pt,inner sep=0pt, outer sep=0pt] at (5.23,1.9) {};
    \draw[draw=DeepBlue,fill=DeepBlue,line width=1pt] (0.8,2.4) -- (1.7,2.4);
    \node[circle,draw=DeepBlue,fill=DeepBlue,minimum size=4pt,inner sep=0pt, outer sep=0pt] at (1.25,2.4) {};
    \draw[draw=MyBlue,fill=MyBlue,line width=1pt] (4.78,2.4) -- (5.68,2.4);
    \node[diamond,draw=MyBlue,fill=MyBlue,minimum size=4.5pt,inner sep=0pt, outer sep=0pt] at (5.23,2.4) {};
    \draw[draw=MyPurple,fill=MyPurple,line width=1pt] (8.76,2.4) -- (9.66,2.4);
    \node[draw=MyPurple,fill=MyPurple,minimum size=4pt,inner sep=0pt, outer sep=0pt] at (9.21,2.4) {};
    \node[right] at (0,0) {\includegraphics[height=2.6cm]{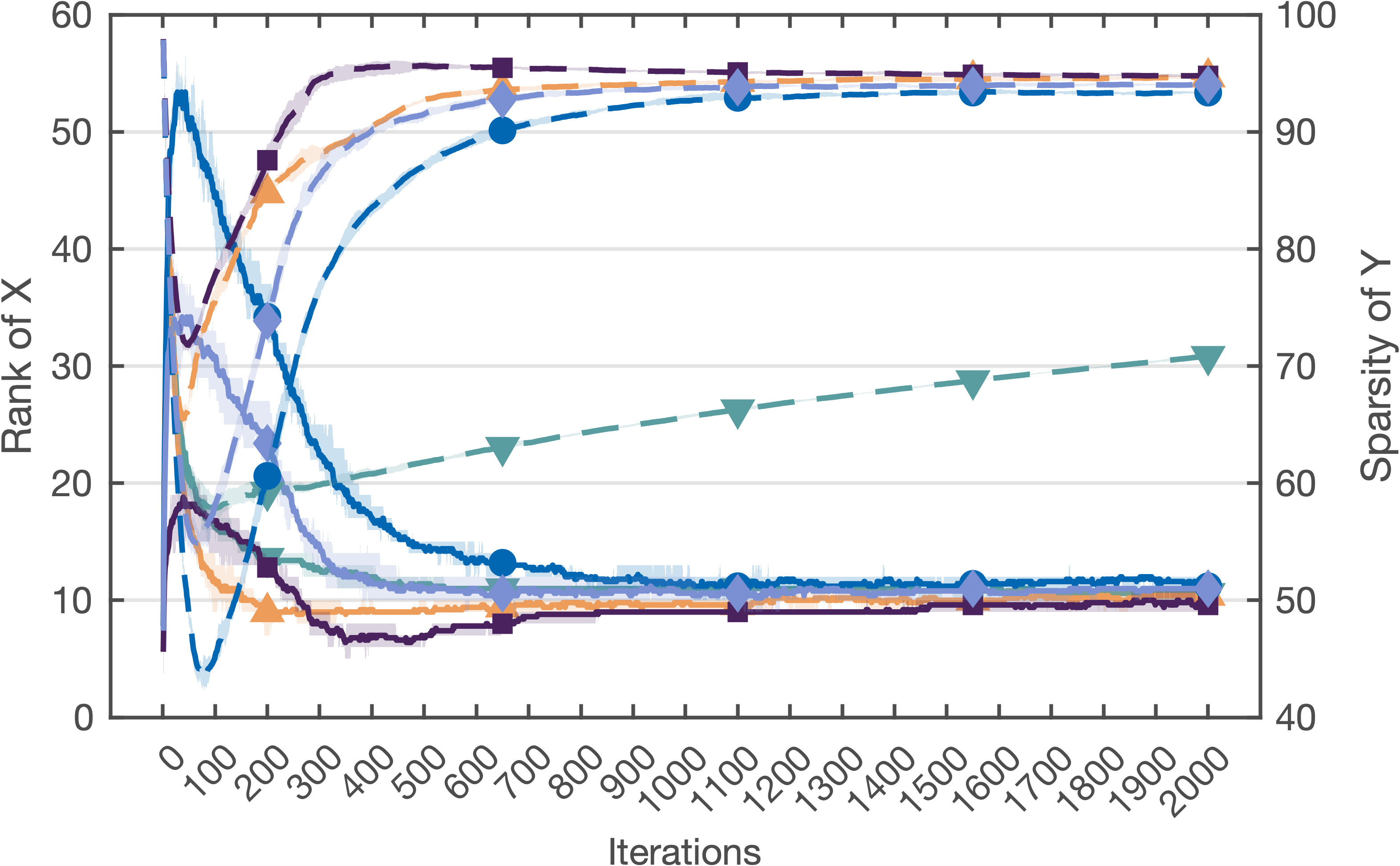}};
    \node[right] at (4.45,0) {\includegraphics[height=2.6cm]{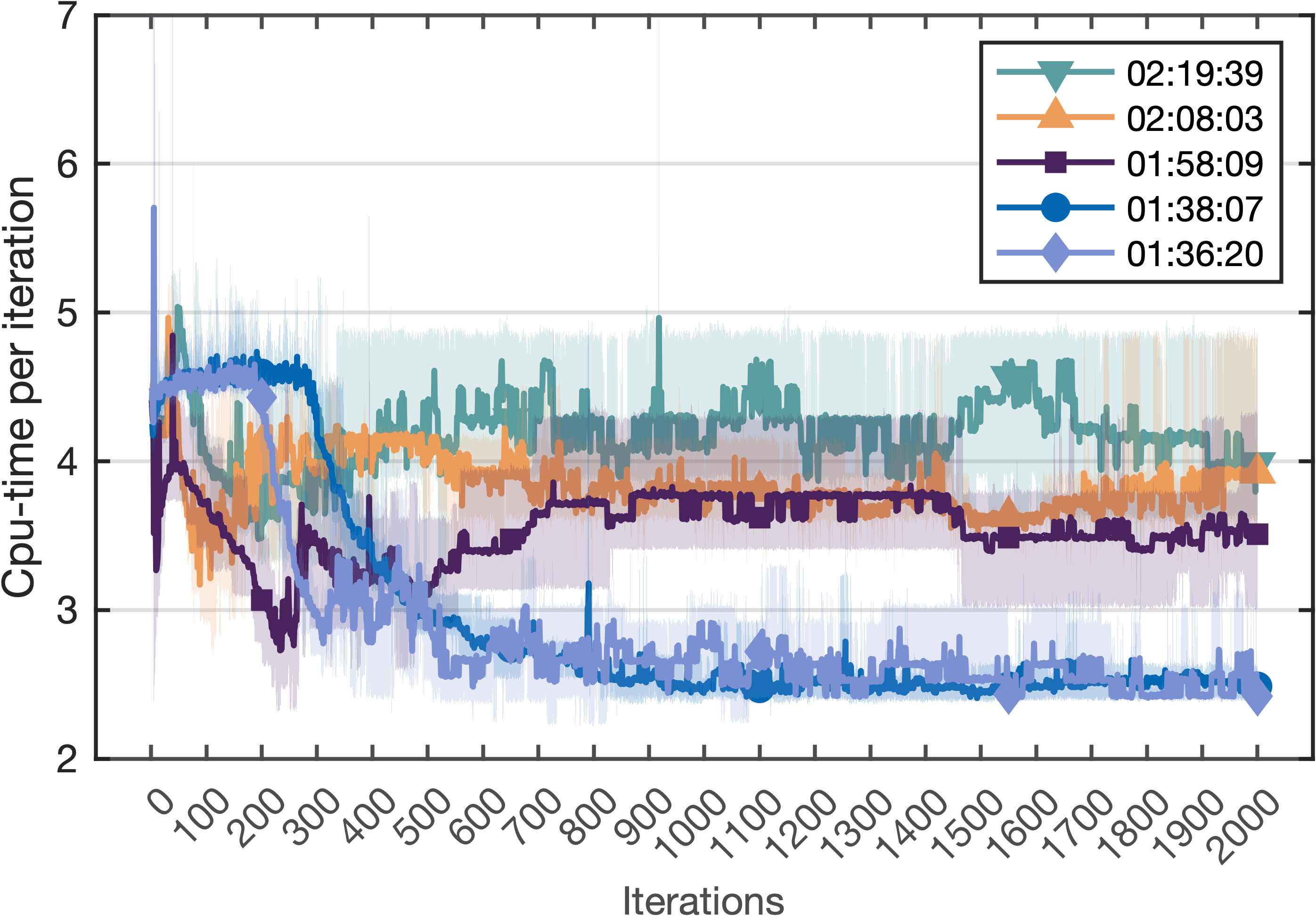}};
    \node[right] at (8.4,0) {\includegraphics[height=2.6cm]{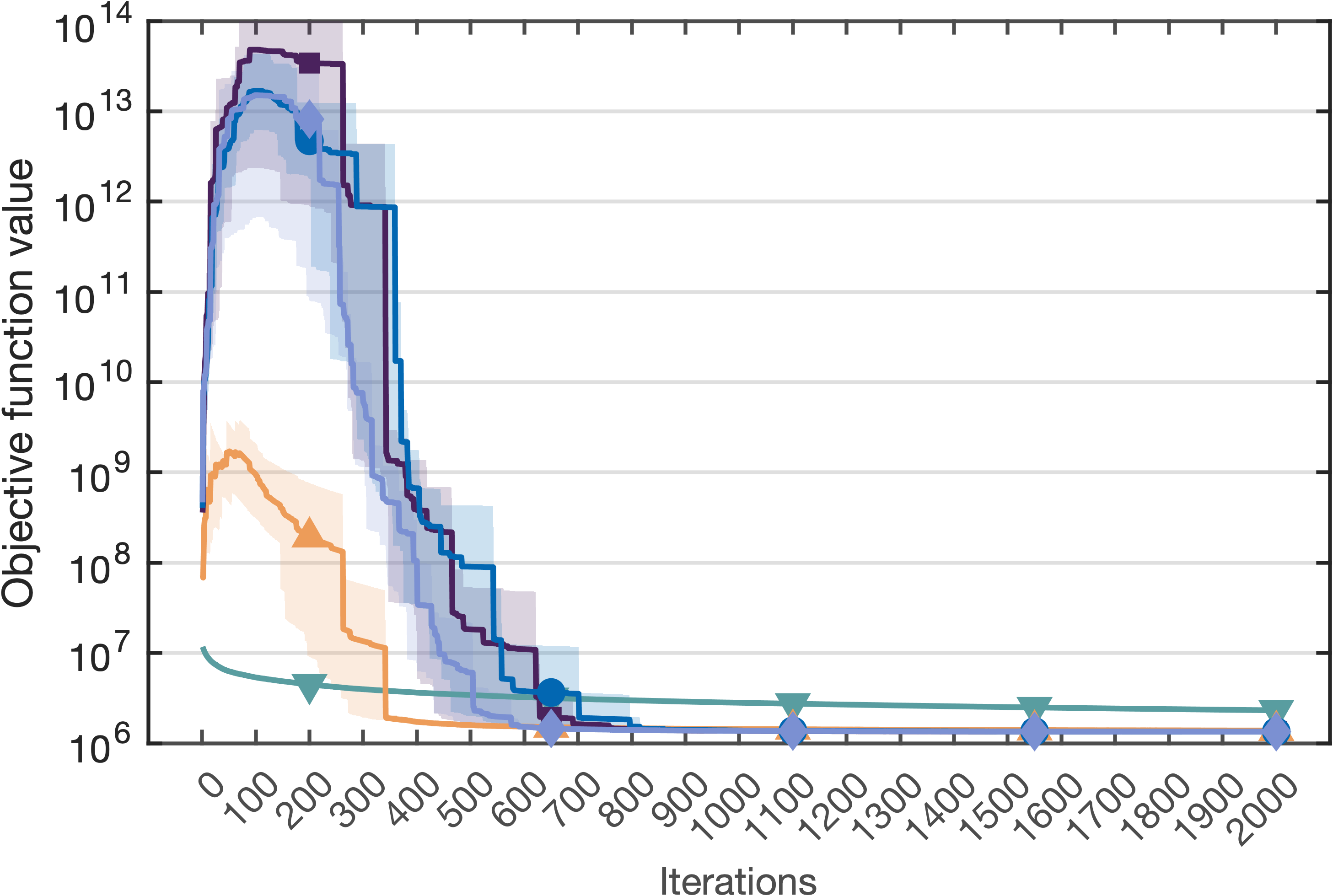}};
    \draw[draw=gray,densely dotted] (11.76,-0.73) -- (10.3,0);
    \draw[draw=gray,densely dotted] (12.2,-0.73) -- (12.3,0);
    \draw[draw=gray] (11.43,-0.89) rectangle (12.29,-0.73); 
    \node[right,draw=gray,fill=white,inner sep=0.5mm] at (10,0.6) {\includegraphics[height=1.6cm]{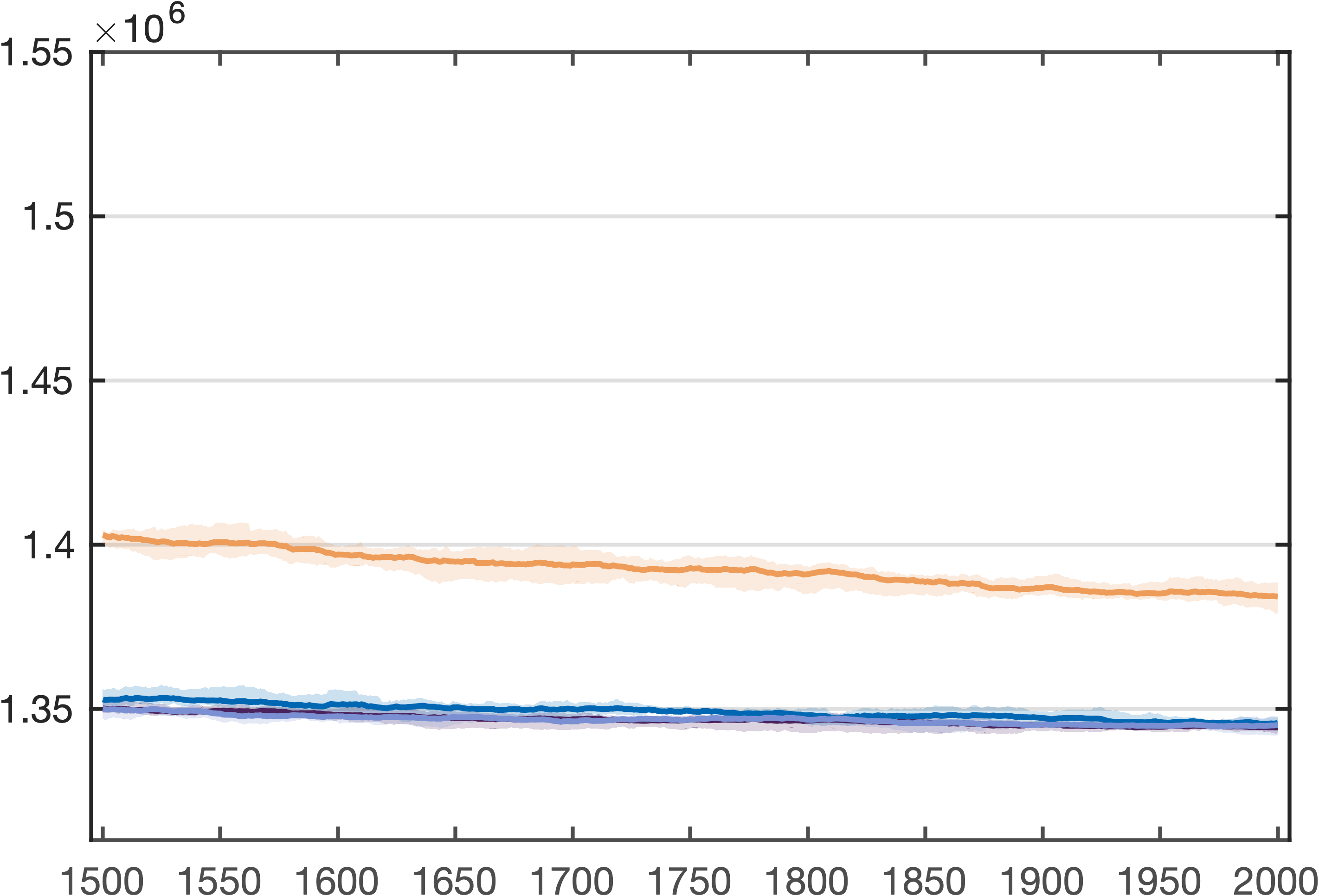}};
    %
\end{tikzpicture}
%
\caption{\revise{Performance of $\NSGD$ and $\RDA$ on the problem \cref{eq:prob-pcp}. $\NSGD$ is run with the step sizes $\alpha_k = \frac12(k+1)^{-2/3}$. In the first figure, the rank of $\{X^k\}_{k\in\N}$ (using solid lines) and the sparsity level $100\%\cdot|\{i,j:Y_{ij}^k = 0\}|/(mn)$ of $\{Y^k\}_{k\in\N}$ (using dashed lines) are plotted. The cpu-time per iteration is shown in the middle figure; the legend depicts the total running time. The right plot illustrates the change in the objective function values. (Averaged over $5$ runs).}}
\label{fig:exp-rda}
\end{figure}

\textbf{Implementational details.} In our experiment, we consider a (gray scale) video with 351 frames and resolution $360\times640$\footnote{Original video by Life-Of-Vids, \textit{San Francisco street slope down up}. Pixabay. \url{https://pixabay.com/videos/san-francisco-street-slope-down-up-3133/}. Last accessed 2025-02-26.}, i.e., $M\in[0,1]^{230400\times 351}$. 
We set $\nu_1 = 150$, $\nu_2 = 0.25$ and fix the mini-batch size to $|S_k| = 8$ (for all $k$). All initial points are set to $0$ and both $\NSGD$ and $\PSGD$ are run for $2000$ iterations with step sizes of the form $\alpha_k = 1/[\sL(k+1)^\gamma]$, where $\gamma \in \{\frac23,\frac34,1\}$ and $\sL = 2$ is the Lipschitz constant of $\nabla f$. In $\NSGD$, we select $\lambda \in \{1,2\}$. \revise{The setup for $\RDA$ generally follows \cite{xiao10a,lee2012manifold} and we run $\RDA$ with $\beta_k = \gamma \sqrt{k}$, $\gamma \in \{5,25\}$, and regularizer $h(X,Y) = \frac12 \|X\|_F^2 + \frac12 \|Y\|_F^2$; see also \Cref{app:comparison}. In this case, the updates of $\RDA$ reduce to 
\[ X^{k+1} = \mathrm{prox}_{\nu_1 \frac{\sqrt{k}}{\gamma} \|\cdot\|_*}(-\tfrac{\sqrt{k}}{\gamma}{\overline{G}^k_x}) \quad \text{and} \quad Y^{k+1} = \mathrm{prox}_{\nu_2 \frac{\sqrt{k}}{\gamma} \|\cdot\|_1}(-\tfrac{\sqrt{k}}{\gamma}{\overline{G}^k_y}), \]
where the matrices ${\overline{G}^k_x}$ and $\overline{G}^k_y$ are averaged gradients (based on the mini-batch gradients $\{g^k\}_{k \in \N}$).} The proximity operator associated with the nuclear norm is given by $\prox{\ell\nu_1\|\cdot\|_*}(X) = U\mathrm{diag}(\prox{\ell\nu_1\|\cdot\|_1}(\sigma))V^\top$, where $X = U\mathrm{diag}(\sigma)V^\top$ is the singular value decomposition (SVD) of $X$, \cite{Bec17}. (In $\NSGD$, $\ell$ is set to $\ell = \lambda$; in $\PSGD$, we have $\ell = \alpha_k$). In iteration $k$, if the rank $r_{k-1}$ of the previous iterate $X^{k-1}$ falls below $25$, we use $\texttt{svds}$ to only compute the decomposition related to the top $r_{k-1}+2$ singular values. This threshold is adaptively increased, if necessary, to ensure that $\prox{\ell\nu_1\|\cdot\|_*}$ is calculated correctly.  Otherwise, if $r_{k-1} \geq 25$, we use \texttt{MATLAB}'s inbuilt \texttt{svd} command to compute the SVD.

In our tests, we track the \revise{sparsity level $|\{i,j: Y_{ij}^k = 0\}|/(mn)$} of $\revise{\{Y^k\}_{k\in\N}}$, the \revise{rank $r_k = |\{i: \sigma_i > 0, X^k = U\mathrm{diag}(\sigma)V^\top\}|$} of $\revise{\{X^k\}_{k\in\N}}$, the required cpu-time per iteration, and the function values. We also run $\mathsf{FISTA}$ with high accuracy to generate reference solutions for \cref{eq:prob-pcp}. $\mathsf{FISTA}$ returns $X^*, Y^*$ with $\mathrm{rank}(X^*) = 11$; the sparsity level of $Y^*$ is $93.81\%$; the optimal value is $\psi^* = 1.333 \cdot 10^6$. As shown in \cref{fig:exp-3}, $\NSGD$ consistently generates iterates with significantly lower rank and higher levels of sparsity. In addition, $\NSGD$ achieves lower objective function values for $\gamma \in \{\frac23,\frac34\}$. When $\NSGD$ is run with $\lambda = 2$, it outperforms $\PSGD$ for all of the tested choices of the step size parameter $\gamma$. The low-rank structure identified by $\NSGD$ has a favorable effect on the run time and can reduce the computational time by a factor of $\approx1.5$ compared to $\PSGD$. (This can be likely further improved using more advanced singular value computations and the $\mathsf{PROPACK}$ package, \cite{larsenPropack}). The results in \cref{fig:exp-3} and \cref{table:one} further illustrate that the manifold identification property of $\NSGD$ improves when $\gamma$ is small (i.e., when larger step sizes are used). \revise{Overall, as shown in \Cref{fig:exp-rda}, $\RDA$ has convincing identification properties when $\gamma = 5$ and performs comparable to $\NSGD$. $\RDA$ slightly overestimates the sparsity level and the decrease in function values is generally slower compared to $\NSGD$. In the case $\gamma = 25$, convergence is slower and $\RDA$ will require more steps to identify the correct sparsity pattern of the solution. This indicates that a more careful tuning of $\gamma$ or of the choice of the regularization function $h$ might be necessary to further enhance $\RDA$'s performance and competitiveness in this problem.}


\begin{table}[t]
\centering
{\scriptsize
\setlength\tabcolsep{4.2pt}
\NiceMatrixOptions{cell-space-limits=2pt}
\begin{NiceTabular}{|p{1.6cm}|p{0.2cm}|p{1.425cm}p{1.425cm}|p{1.425cm}p{1.425cm}|p{1.425cm}p{1.425cm}|}%
 [ 
   code-before = 
    \rectanglecolor{MyBlue!10}{4-3}{4-8}
    \rectanglecolor{MyBlue!10}{7-3}{7-8}
 ]
\toprule
\Block[c]{2-2}{\textbf{Alg.}} & & \Block[c]{1-2}{$\gamma = \frac23$} & & \Block[c]{1-2}{$\gamma = \frac34$} & & \Block[c]{1-2}{$\gamma = 1$} & \\
& & \Block[c]{1-1}{last $500$ it.} & \Block[c]{1-1}{last $100$ it.} & \Block[c]{1-1}{last $500$ it.} & \Block[c]{1-1}{last $100$ it.} & \Block[c]{1-1}{last $500$ it.} & \Block[c]{1-1}{last $100$ it.} \\
\Hline 
$\PSGD\phantom{{}^{(0)}}$ & {\Block[c]{3-1}<\rotate>{rank}} & 61.53 (1.23) & 62.02 (1.14) & 63.07 (1.24) & 63.13 (1.48) & 61.76 (1.17) & 61.84 (1.16) \\
$\NSGD^{(1)}$ & & 11.54 (0.20) & 11.54 (0.12) & 11.09 (0.10) & 11.02 (0.06) & 9.81 (0.05) & 9.8 (0.00) \\
$\NSGD^{(2)}$ & & 10.87 (0.14) & 11 (0.00) & 10.81 (0.10) & 10.8 (0.00) & 8.2 (0.00) & 8.2 (0.00) \\ \Hline
$\PSGD\phantom{{}^{(0)}}$ & {\Block[c]{3-1}<\rotate>{sparsity}} & 79.85 (0.17) & 79.84 (0.17) & 78.28 (0.35) & 78.53 (0.20) & 53.00 (0.55) & 53.76 (0.12) \\
$\NSGD^{(1)}$ & & 93.34 (0.04) & 93.34 (0.02) & 92.96 (0.12) & 93.11 (0.02) & 79.21 (0.38) & 79.73 (0.07) \\
$\NSGD^{(2)}$ & & 94.00 (0.03) & 94.01 (0.01) & 93.75 (0.05) & 93.79 (0.00) & 90.06 (0.13) & 90.23 (0.03) \\
\bottomrule
\end{NiceTabular}
\noindent
\NiceMatrixOptions{cell-space-limits=2pt}
\begin{NiceTabular}{|p{1.6cm}|p{1.2cm}|p{2.06cm}p{2.06cm}|p{2.06cm}p{2.06cm}|}%
 [ 
   code-before = 
    \rectanglecolor{MyBlue!10}{4-3}{4-6}
 ]
\toprule
\Block[c]{2-2}{\textbf{Alg.}} & & \Block[c]{1-2}{$\beta_k = 5\sqrt{k}$} & & \Block[c]{1-2}{$\beta_k = 25\sqrt{k}$} & \\
& & \Block[c]{1-1}{last $500$ it.} & \Block[c]{1-1}{last $100$ it.} & \Block[c]{1-1}{last $500$ it.} & \Block[c]{1-1}{last $100$ it.} \\
\Hline 
\Block[c]{2-1}{$\RDA$} & rank & \Block[c]{1-1}{10.21 (0.16)} & \Block[c]{1-1}{10.4 (0.00)} & \Block[c]{1-1}{10.71 (0.14)} & \Block[c]{1-1}{10.60 (0.02)} \\
& sparsity & \Block[c]{1-1}{94.58 (0.05)} & \Block[c]{1-1}{94.65 (0.02)} & \Block[c]{1-1}{69.70 (0.67)} & \Block[c]{1-1}{70.62 (0.14)} \\
\bottomrule
\end{NiceTabular}
}
\caption{Average rank and sparsity level achieved by \revise{$\NSGD$, $\PSGD$, and $\RDA$}. Averages are taken with respect to the last $500$ and $100$ iterations. Both $\NSGD^{(\lambda)}$ and $\PSGD$ are run for $2000$ iterations with step sizes $\alpha_k = \frac12(k+1)^{-\gamma}$, $\gamma \in \{\frac23,\frac34,1\}$ and $\lambda \in \{1,2\}$. \revise{$\RDA$ is run for $2000$ iterations with $\beta_k = \gamma\sqrt{k}$ and $\gamma \in \{5,25\}$.} The averaged rank \revise{$|\{i:\sigma_i > 0, X^k = U \mathrm{diag}({\sigma})V^\top\}|$} and averaged sparsity level \revise{$100\%\cdot|\{i,j:Y_{ij}^k=0\}|/(mn)$} are reported. Standard deviations are shown in parentheses. (Averaged over $5$ runs).
}
\label{table:one}
\vspace{-2ex}
\end{table}

\section{Conclusion}
In this paper, we develop a novel normal map-based proximal stochastic gradient method, $\NSGD$, for nonconvex composite-type problems. Leveraging the unbiasedness of the stochastic normal map steps, we show that the iterates generated by $\NSGD$ satisfy an approximate descent condition which allows establishing asymptotic global convergence results and non-asymptotic complexity bounds under standard assumptions.
Furthermore, we prove a$.$s$.$ strong limit convergence guarantees if the objective function is definable. Based on these results, $\mathsf{Norm}\textsf{-}$ $\mathsf{SGD}$ is shown to a$.$s$.$ have a manifold identification property. This is in stark contrast to the classical proximal stochastic gradient method, which is known to not possess such a property. To our knowledge, $\NSGD$ appears to be one of the first ``variance reduction-free'', basic proximal stochastic algorithms for nonconvex composite optimization with this broad range of convergence properties. We believe that the normal map-based perspective and the KL-based techniques adopted in this work have general potential and can be likely applied to other families of stochastic approaches. 

\section*{Acknowledgments}
We are grateful for the valuable discussions with Lexiao Lai, Kun Huang, Ziqi Qin, and \revise{Ching-pei Lee}. \revise{Moreover, we would like to thank the associate editor and the anonymous referees for their careful reading and important comments, which greatly helped to improve the quality and presentation of this work.}

\appendix

\section{\texorpdfstring{Comparing $\NSGD$ and $\RDA$}{Comparing Norm-SGD and RDA}} \label{app:comparison}
\revise{In this section, we conduct a formal comparison between $\NSGD$ and $\RDA$, \cite{xiao10a,lee2012manifold}. $\RDA$ is a generalization of Nesterov’s dual averaging method \cite{nesterov2009primal}, which was developed for convex optimization problems. Using our notations, $\RDA$---applied to the composite problem \cref{eq:SO}---takes the following form:
\[ \sx^{k+1} = \argmin_{x \in \Rn}~\frac1k\Big\langle {\sum}_{i=1}^k\, \sg^i, x \Big\rangle + \varphi(x) + \frac{\beta_k}{k} \cdot h(x). \]
Here, the function $h : \Rn \to \Rex$ is supposed to be strongly convex on $\dom{\vp}$ with $\argmin_x h(x) \subseteq \argmin_x \vp(x)$, $\sg^i$ is a stochastic approximation of $\nabla f(\sx^i)$, $i \in \N$, and $\{\beta_k\}_{k \in \N}$ are given parameters. 
Mimicking \cref{eq:fixed-point}, the update of $\NSGD$ can be written as $\sz^{k+1} = (1 -{\alpha_k}{\lambda^{-1}}) \sz^k + {\alpha_k}{\lambda^{-1}} \sy^k$, where $\sy^k := \sx^k - \lambda \sg^k $ and $\sx^{k+1} = \mathrm{prox}_{\lambda \varphi}(\sz^{k+1})$. This implies
\begin{equation} \label{eq:rda-zk} \sz^{k+1} = \eta_{k,0} \cdot \sz^1 + {\sum}_{i=1}^k \eta_{k,i} \cdot  \sy^i, \end{equation}
where $\eta_{k,0} := \prod_{i=1}^k (1 - \frac{\alpha_i}{\lambda})$, $\eta_{k,i} := \frac{\alpha_i}{\lambda} \prod_{j=i+1}^k (1 - \frac{\alpha_j}{\lambda})$ for $1\leq i \leq k-1$, and $\eta_{k,k} = \frac{\alpha_k}{\lambda}$. Furthermore, if we choose $\alpha_1 = \lambda$ and $\alpha_i < \lambda$ for all $i \geq 2$, then it holds that $\eta_{k,0} = 0$ and $\sum_{i=1}^k \eta_{k,i} = 1$ with $\eta_{k,i}>0$ for all $1\leq i \leq k$. In this setting, equation \cref{eq:rda-zk} reduces to 
\[ \sz^{k+1} = {\sum}_{i=1}^k \eta_{k,i} \sx^i - \lambda {\sum}_{i=1}^k \eta_{k,i} \sg^i. \]
%
As a result, the update rule of $\NSGD$ can be expressed as
\begin{align*}
\sx^{k+1} & = \argmin_{x \in\Rn}~\varphi(x) + \frac{1}{2\lambda} \|{x - \sz^{k+1}}\|^2 \\
& = \argmin_{x \in \Rn} \Big\langle {\sum}_{i=1}^k \eta_{k,i} \sg^i, x \Big\rangle + \varphi(x) + \frac{1}{2\lambda} \Big\| x - {\sum}_{i=1}^k \eta_{k,i} \sx^i \Big\|^2.
\end{align*}
Hence, the updates of $\NSGD$ and $\RDA$ are both defined through solving similar subproblems,
\[ \min_{x \in \Rn}~\iprod{\bar{\sg}}{x} + \vp(x) + h_k(x), 
\]
where $\bar{\sg}$ is an appropriately averaged gradient and $h_k$ is an additional regularization term. However, both methods are based on different choices of $\bar{\sg}$ and $h_k$. The gradient averaging schemes for $\NSGD$ and $\RDA$ are given by: 
\[
\text{$\NSGD$:} \;\;\; \bar{\sg}_{\mathsf{n}} = {\sum}_{i=1}^k \eta_{k,i} \sg^i \quad \text{and} \quad \text{$\RDA$:} \;\;\; \bar{\sg}_{\mathsf{r}} = \frac{1}{k}{\sum}_{i=1}^k \sg^i.
\]
Clearly, $\bar{\sg}_{\mathsf{n}}$ can recover $\bar{\sg}_{\mathsf{r}}$ via setting $\alpha_i = \frac{\lambda}{i}$. A more fundamental difference lies in the regularization term $h_k$. In $\RDA$, the regularizer uses a fixed center---typically a minimizer of $\varphi$---while updating only the scaling factor $\frac{\beta_k}{k}$ at each iteration. By contrast, $\NSGD$ maintains a fixed scaling factor but adaptively shifts the center using a moving average of the past iterates.  Whether $\NSGD$ can be extended to accommodate broader classes of regularizers while preserving convergence guarantees in both convex and nonconvex settings remains an open question for future research.}


\section{Stochastic and deterministic tools} \label{sec:tools}
We introduce several stochastic concepts and techniques that are the backbone of our analysis. In order to control the stochastic behavior of the error terms $\se^k$ and to establish a$.$s$.$ convergence of the sequence of iterates $\revise{\{\sx^k\}_{k\in\N}}$ generated by $\NSGD$, we use several results from martingale theory. We first state a standard convergence theorem for vector-valued martingales, see, e.g., \cite[Theorem 5.2.22]{Str11} and \cite[Theorem 5.14]{Bre92}.

\begin{lemma}[Martingale convergence theorem]
    \label{Thm:mart_conv}
    Let $\revise{\{\sw^k\}_{k\in\N}}$ be a given vector-valued martingale with respect to some filtration \revise{$\{\mathcal{U}_k\}_{k\in\N}$}. If $\sup_k \Exp[{\|\sw^k\|}] < \infty$, then \revise{$\{\sw^{k}\}_{k\in\N}$} converges a$.$s$.$ to an integrable random vector $\sw$.
\end{lemma}

We also require the Burkholder-Davis-Gundy inequality \cite{BurDavGun72,Str11} in our analysis. 

\begin{lemma}[Burkholder-Davis-Gundy inequality]
    \label{Thm:BDG}
    Let \revise{$\{\sw^{k}\}_{k\in\N}$} be a given vector-valued martingale with an associated filtration $\revise{\{\mathcal U_k\}_{k\in\N}}$ and $\sw^{0}=0$. Then, for all $p \in (1,\infty)$, there exists $C_p > 0$ such that \vspace{-1ex}
    
    \[ \Exp\Big[{{\sup}_{k \geq 0} \|\sw^k\|^p}\Big] \leq C_p \cdot \Exp\Big[{\Big({{\sum}_{k=1}^\infty \|\sw^{k}-\sw^{k-1}\|^2 }\Big)^{{p}/{2}} }\Big]. \]
\end{lemma}

The constant $C_p$ in \cref{Thm:BDG} is universal and it holds that $C_2 \leq 4$, \cite{BurDavGun72,Str11}.

Next, we introduce Gronwall's inequality \cite{gronwall1919note,borkar2009stochastic,tadic2015convergence}. We state a discrete variant of Gronwall's inequality, which can be found in, e.g., \cite[Appendix B, Lemma 8]{borkar2009stochastic}.

\begin{lemma}[Gronwall's inequality]
    \label{Thm:G}
    Let $\revise{\{a_k\}_{k\in\N}} \subseteq \R_{++}$ and $\revise{\{y_k\}_{k\in\N}} \subseteq \R_+$ be given sequences. Suppose that we have $y_0\leq b_1$, and $y_{t+1} \leq b_1 + b_2 \sum_{k=0}^t a_k y_k$ for all $t$ and some $b_1,b_2\geq0$. 
    Then, it holds that $y_{t+1} \leq b_1 \cdot \exp({b_2\sum_{k=0}^t a_k})$ for all $t\geq 0$.
\end{lemma}
 
\section{Derivation of preparatory lemmas} \label{app:pf-prep-lemma}

\subsection{\texorpdfstring{Proof of \cref{lem:est-err}}{Proof of Lemma 2.1}}\label{proof:lem:est-err}

\begin{proof} Let $\omega \in \Omega$ be arbitrary and let us set $\tau = \tau_{m,n}$, $p = p^{m,n} = \scp^{m,n}(\omega)$, $s = s_{m,n} = \scs_{m,n}(\omega)$, $x^k = \sx^k(\omega)$, $z^k = \sz^k(\omega)$, etc. It follows from the definition of $p^{m,i}$ \revise{and $s_{m,n}$} and the triangle inequality that
\begin{align}
        \nonumber \|p^{m,i}\| &\leq {\sum}_{j=m}^{i-1}\alpha_j\|{\Fnor(z^j) - \Fnor(z^m)}\| + \Big\|{{\sum}_{j=m}^{i-1} \alpha_j e^j}\Big\| \\ &\leq {\sum}_{j=m}^{i-1}\alpha_j\|{\Fnor(z^j) - \Fnor(z^m)}\| + s,\qquad \text{for all $i\in(m,n] \cap \N$.} \label{eq:lem-4-est-1}
\end{align}
Recalling $\Fnor(z) = \nabla f(x) + \lambda^{-1}(z-x)$ and using \ref{A1} and the nonexpansiveness of the proximity operator (i.e., $\|x^i-x^m\| \leq \|z^i-z^m\|$), we have
\begin{equation}
\label{eq:nor is Lip}
     \|\Fnor(z^j) - \Fnor(z^m)\| \leq (\sL+\lambda^{-1})\|x^j-x^m\| + \lambda^{-1}\|z^j-z^m\| \leq  \Llambda \|z^j-z^m\|,
\end{equation} 
where $\Llambda =\sL+2\lambda^{-1}$. Thus, combining the relation $z^i = z^m - \tau_{m,i}\Fnor(z^m) + p^{m,i}$ (cf$.$ \cref{eq:pretty-important})\revise{, $\tau_{m,i} \leq \tau$,} and the inequality \cref{eq:lem-4-est-1}, we can infer
\begin{align}
    \nonumber \|z^i-z^m\| &\leq \tau_{m,i} \|\Fnor(z^m)\| + \|p^{m,i}\| \\
    \nonumber &\leq {\sum}_{j=m}^{i-1} \alpha_j \|\Fnor(z^j)-\Fnor(z^m)\| + \tau \|\Fnor(z^m)\| + s\\ 
    \label{eq:lem-4-est-3} &\leq \Llambda {\sum}_{j=m}^{i-1} \alpha_j  \|z^j-z^m\| + \tau \|\Fnor(z^m)\| + s, \qquad \text{for all $i\in(m,n] \cap \N$}.
\end{align}
We now apply Gronwall's inequality (\cref{Thm:G}) upon setting $b_1 :=\tau \|\Fnor(z^m)\| + s$, $b_2 := \Llambda$, $a_k:=\alpha_{m+k}$, $y_k:=\|z^{m+k} - z^m\|$, and $t:=i-m-1$. For all $i\in(m,n] \cap \N$, this establishes the following upper bound  
\begin{equation}
            \|z^i - z^m\|  \leq \exp(\tau \Llambda) \cdot [\tau\|\Fnor(z^m)\|+s]\leq 2\tau\|\Fnor(z^m)\| + 2s,
            \label{eq:lem-4-est-5}
\end{equation}
where the last inequality is due to $\tau\Llambda \leq \frac12$ and $\exp(\frac12)\leq 2$. Moreover, it follows from $\|x^i-x^m\| \leq \|z^i-z^m\|$ that ${\max}_{m<i\leq n}\|x^i-x^m\| \leq {\max}_{m<i\leq n}\|z^i-z^m\| \leq 2\tau\|\Fnor(z^m)\| + 2s$. This verifies the first bound in \cref{lem:est-err}. 
     
Next, combining \cref{eq:nor is Lip} and \cref{eq:lem-4-est-5}, we obtain
\[
{\max}_{m<i\leq n} \|\Fnor(z^i) - \Fnor(z^m)\|\leq \Llambda{\max}_{m<i\leq n}\|z^i-z^m\| \leq 2\Llambda[\tau\|\Fnor(z^m)\| + s].
\]
Using this bound in \cref{eq:lem-4-est-1} with $i=n$ and recalling $\tau={\sum}_{i=m}^{n-1} \alpha_i$, this yields $\|p\|\leq s+ 2\Llambda (\tau\|\Fnor(z^m)\| + s) {\sum}_{j=m}^{n-1} \alpha_j \leq 2\Llambda\tau^2\|\Fnor(z^m)\| + 2s$.
\end{proof}

\subsection{\texorpdfstring{Proof of \cref{lem:err_estimate}}{Proof of Lemma 2.2}} \label{app:proof-martingale} 

\begin{proof}  Let us define the filtration $\mathcal{U}_\ell:=\mathcal{F}_{t_k+\ell}$ and introduce the sequence \revise{$\{\sy^\ell\}_{\ell\in\N}$} as follows $\sy^0:=0$, $\sy^\ell := {\sum}_{i=t_k}^{\min\{t_k+\ell,t_{k+1}\}-1}\alpha_i \se^i$ for all $\ell \geq 1$.
Then, $\sy_\ell$ is $\mathcal U_\ell$-measurable and for all $\ell \in \{1,\dots,t_{k+1}-t_k-1\}$, we have 
\[ {\Exp[{\sy^{\ell+1} \mid \mathcal U_\ell}] = {\sum}_{i=t_k}^{t_k+\ell} \alpha_i\Exp[{\se^i \mid \mathcal U_\ell}] = \sy^\ell + \alpha_{t_k+\ell} \Exp[{\se^{t_k+\ell} \mid \mathcal{F}_{t_k+\ell}}] = \sy^\ell.} \]
(Similarly for $\ell \geq t_{k+1}-t_k$). Hence, \revise{$\{\sy^\ell\}_{\ell\in\N}$} defines a \revise{$\{\mathcal{U}_{\ell}\}_{\ell\in\N}$}-martingale. Let us further set $\sr_k = \scs_{t_k,t_{k+1}}$. By the Burkholder-Davis-Gundy inequality (\cref{Thm:BDG} with $p=2$ and $C_p\leq 4$) and assumption \ref{B1}, we have $\Exp[{\sr_k^2}] =\Exp[{{\sup}_{\ell>0} \|\sy^\ell\|^2}]$ and
\begin{align*} 
\nonumber\Exp[{\sr_k^2}] \leq 4 \Exp\Big[{ {\sum}_{\ell=1}^{\infty}\|\sy^\ell-\sy^{\ell-1}\|^2}\Big] & = 4 {\sum}_{i=t_k}^{t_{k+1}-1}\alpha_i^2 \Exp[{ \|\se^i\|^2}]  \leq 4 {\sum}_{i=t_k}^{t_{k+1}-1} \alpha_i^2\sigma_i^2. 
\end{align*}
This proves the first claim in \Cref{lem:err_estimate}. Since $\revise{\{\beta_k\}_{k\in\N}}$ is non-decreasing, we obtain%
\begin{equation}\label{eq:app-bdg}
	\beta_{t_k}^2 \Exp[\sr_k^2] \leq 4 \beta_{t_k}^2 {\sum}_{i=t_k}^{t_{k+1}-1} \alpha_i^2\sigma_i^2 \leq 4  {\sum}_{i=t_k}^{t_{k+1}-1} \alpha_i^2\beta_{i}^2\sigma_i^2
\end{equation}
for all $k \in \N$. The monotone convergence theorem and \cref{eq:app-bdg} further imply
\begin{align*}
 \Exp\Big[{ {\sum}_{k=0}^\infty \beta_{t_k}^2{\sr_k^2}}\Big] = {\sum}_{k=0}^{\infty}\beta_{t_k}^2 \Exp[\sr_k^2] 
 \leq 
    {4} {\sum}_{i=t_{0}}^\infty \alpha_i^2\beta_i^2\sigma_i^2 < \infty,
\end{align*}
and consequently, it holds that $\sum_{k=0}^\infty \beta_{t_k}^2{\sr_k^2} < \infty$ a$.$s$.$.
\end{proof}

\subsection{\texorpdfstring{Proof of \cref{lem:basic bounds}}{Proof of Lemma 2.3}}\label{app:lem:basic bounds}
   
\begin{proof}[Proof of (a)]
Let $\omega\in\Omega$ be arbitrary and denote $x^k = \sx^k(\omega)$, $z^k = \sz^k(\omega)$, $p = \scp^{m,n}(\omega)$, $u = \su^{m,n}(\omega)$, $\tau = \tau_{m,n}$, etc$.$. Recalling $\Fnor(z) = \nabla f(x) + \lambda^{-1}(z-x)$ and by the definition \cref{eq:u-def1} of $\su^{m,n}$, we can express $\Fnor(z^n)$ as
\begin{align*} 
    \Fnor(z^n) 
    & = \Fnor(z^m) + [\nabla f(x^n)-\nabla f(x^m)] + \lambda^{-1}\prt{z^n-z^m+x^m-x^n}\\
    & = u + [\nabla f(x^n)-\nabla f(x^m)] - \lambda^{-1}(x^n-x^m).
\end{align*}	
Taking squares and using the Lipschitz continuity of $\nabla f$, this yields
\begin{align}
            \nonumber \norm{\Fnor(z^n)}^2 & = \norm{u}^2 + 2\iprod{u}{(\nabla f(x^n)-\nabla f(x^m)) - \lambda^{-1}(x^n-x^m)} \\
            \nonumber&\hspace{4ex}  +  \norm{(\nabla f(x^n)-\nabla f(x^m)) - \lambda^{-1}(x^n-x^m)}^2\\
            \label{C1} & \hspace{-6ex} \leq \norm{u}^2 + 2\iprod{u}{(\nabla f(x^n)-\nabla f(x^m))-\lambda^{-1}(x^n-x^m)} + (\sL+{\lambda}^{-1})^2 \norm{x^n-x^m}^2.
\end{align}
Using $u=(1-\lambda^{-1}\tau)\Fnor(z^m) + \lambda^{-1}p$ (cf$.$ \cref{eq:u-def3}) and $\norm{a+b}^2 \leq (1+\varepsilon) \norm{a}^2 + (1+\frac{1}{\varepsilon})\norm{b}^2$ with $a= (1-\lambda^{-1}\tau) \Fnor(z^m)$, $b = \lambda^{-1}p$, and $\varepsilon = \frac{\lambda^{-1}\tau}{1-\lambda^{-1}\tau}$, we have
\begin{equation}\label{C2}
        \norm{u}^2 \leq (1-\lambda^{-1}{\tau}) \|{\Fnor(z^m)}\|^2 + (\lambda\tau)^{-1} \|{p}\|^2.
\end{equation}
Utilizing Young's inequality for $2\iprod{u}{\nabla f(x^n)-\nabla f(x^m)}$ and \cref{C2}, we obtain 
\begin{align}
            \nonumber 2\iprod{u}{\nabla f(x^n)-\nabla f(x^m)} & \leq \tfrac{\tau}{2\lambda} \norm{u}^2 + \tfrac{2\lambda}{\tau} \norm{\nabla f(x^n)-\nabla f(x^m)}^2\\
            \label{C3} & \leq  \tfrac{\tau}{2\lambda}\norm{\Fnor(z^m)}^2 + \tfrac{1}{2\lambda^2} \norm{p}^2 + \tfrac{2\lambda \sL^2}{\tau} \norm{x^n-x^m}^2,
\end{align}
Summing up the estimates \cref{C1}--\cref{C3}, and noticing $\tau \leq \lambda$, it follows
\begin{align*}
        \norm{\Fnor(z^n)}^2  & \leq  ({1-\tfrac{\tau}{2\lambda}})\|{\Fnor(z^m)}\|^2 - \tfrac{2}{\lambda}\iprod{u}{x^n-x^m} \\
        & \hspace{4ex} + [(\sL+\lambda^{-1})^2 + {2\sL^2\lambda}{\tau^{-1}}] \norm{x^n-x^m}^2 + \tfrac{3}{2\lambda\tau} \norm{p}^2.
\end{align*}

\emph{Proof of (b).} By \ref{A1}--\ref{A2}, $\nabla f$ is $\sL$-continuous and $\vp$ is convex. Thus, we have
\begin{equation*}
        \psi(x^n)-\psi(x^m) \leq \iprod{\nabla f(x^m)}{x^n-x^m}+\tfrac{\sL}{2}\|x^n-x^m\|^2 + \tfrac{1}{\lambda}\iprod{z^n-x^n}{x^n-x^m},
\end{equation*}
where we used $\lambda^{-1}(z^n-x^n) \in\partial \vp(x^n)$ in the last inequality. Moreover, based on the definition of $\Fnor$, we can reformulate 
the above inequality as
\begin{align*}
\psi(x^n)-\psi(x^m) & \leq \iprod{\nabla f(x^m)+\lambda^{-1}(z^m-x^m)}{x^n-x^m}+\tfrac{\sL}{2}\|{x^n-x^m}\|^2\\
& \hspace{4ex} + \lambda^{-1}\iprod{z^n-x^n}{x^n-x^m} - \lambda^{-1}\iprod{z^m-x^m}{x^n-x^m}\\
& \hspace{-4ex} = [{\tfrac{\sL}{2}-\tfrac{1}{\lambda}}]\|{x^n-x^m}\|^2 + 	\iprod{\Fnor(z^m)}{x^n-x^m} + \tfrac{1}{\lambda}\iprod{z^n-z^m}{x^n-x^m}.
\end{align*}
This concludes the proof by substituting \cref{eq:u-def1} into the last expression.
\end{proof}

\subsection{\texorpdfstring{Proof of \cref{lem:index-time bound}}{Proof of Lemma 3.5}} \label{proof:lem:index-time bound}

\begin{proof}
As shown in \cref{eq:nor is Lip}, $\Fnor$ is Lipschitz with modulus $\Llambda$ and we have 
\begin{align*}
        \|\Fnor(\sz^i)\| &\leq \|\Fnor(\sz^i) - \Fnor(\sz^{t_{\fk(i)}})\| + \| \Fnor(\sz^{t_{\fk(i)}})\|\\ &\leq  \Llambda \|\sz^i-\sz^{t_{\fk(i)}}\|+ \| \Fnor(\sz^{t_{\fk(i)}})\| \leq (2\Llambda\sT+1)\| \Fnor(\sz^{t_{\fk(i)}})\| + 2\Llambda \scs_{t_{\fk(i)},i}\\
        &\leq (2\Llambda\sT+1)\| \Fnor(\sz^{t_{\fk(i)}})\| + 2\Llambda \sr_{\fk(i)},
\end{align*}
where the \revise{second to last} inequality follows from \cref{lem:est-err} by setting $m=t_{\fk(i)}$, $n=t_{\fk(i)+1}$ (the condition $\tau_{t_{\fk(i)},t_{\fk(i)+1}}\leq \sT \leq 1/(2\Llambda)$ is due to \cref{eq:def delta and T}). Squaring both sides of the above inequality and using $(a+b)^2\leq 2a^2+2b^2$, we obtain \cref{eq:tech1}.

Based on \cref{lem:basic bounds} (b) and the definition of $\su^{m,n}$ (see \cref{eq:u-def1}), it holds that
\[
    \psi(\sx^n) \leq \psi(\sx^m) + (\sL/2-\lambda^{-1})\norm{\sx^n-\sx^m}^2 + \iprod{\lambda^{-1}(\sz^n-\sz^m) + \Fnor(\sz^m)}{\sx^n-\sx^m}.
\]
Applying the nonexpansiveness of the proximity operator, \cref{prox-prop}, and $\iprod{\Fnor(\sz^m)}{\sx^n-\sx^m} \leq \frac{\lambda}{4} \norm{\Fnor(\sz^m)}^2 + \lambda^{-1} \norm{\sx^n-\sx^m}^2$, we have
\[
    \psi(\sx^n)\leq  \psi(\sx^m) + \tfrac{\Llambda}{2} \norm{\sz^n-\sz^m}^2 + \tfrac{\lambda}{4} \norm{\Fnor(\sz^m)}^2.
\]
Invoking \Cref{lem:est-err} (for $\norm{\sz^n-\sz^m}$) and $(a+b)^2\leq 2a^2+2b^2$,  we obtain
\[
    \psi(\sx^n)\leq  \psi(\sx^m) + (4\Llambda\tau_{m,n}^2+\tfrac{\lambda}{4}) \norm{\Fnor(\sz^m)}^2 + 4\Llambda \scs_{m,n}^2.
\]
Then, \cref{eq:tech2} readily follows by setting $n = i$, $m = t_{\fk(i)}$ and using $\tau_{t_{\fk(i)},i}\leq \sT$. Similarly, \cref{eq:tech3} can be obtained by setting $n = t_{\fk(i)+1}$, $m = i$ and using  $\tau_{i,t_{\fk(i)+1}}\leq \sT$ and $\scs_{i,t_{\fk(i)+1}} = \max_{i<j\leq t_{\fk(i)+1}}\|{\sum}_{\ell=i}^{\revise{j}-1} \alpha_\ell \se^\ell\| \leq \scs_{t_{\fk(i)},i}+  \scs_{t_{\fk(i)},t_{\fk(i)+1}}  \leq 2 \sr_{\fk(i)}$.
\end{proof}

\subsection{\texorpdfstring{Proof of \cref{lem:subsequence converge}}{Proof of Lemma 4.2}}\label{proof:lem:subsequence converge}  
\begin{proof}
Let $\omega\in\mathcal C$ be arbitrary and let $x^k = \sx^k(\omega)$,  $z^k = \sz^k(\omega)$, etc$.$ be the corresponding realizations. By the definition of the event $\mathcal C$, there exists $\revise{\{p_k\}_{k\in\N}}\subseteq\N$ such that \revise{$\{x^{p_k}\}_{k\in\N}$} is bounded and we have $\Fnor(z^k) \to 0$, $\psi(x^k) \to \psi^*$. Using the definition of the normal map, we can infer $\|z^{p_k}\| \leq \lambda(\|\Fnor(z^{p_k})\|+\|\nabla f(x^{p_k})\|)+\|x^{p_k}\|$. Hence, \revise{$\{z^{p_k}\}_{k\in\N}$} is also bounded and there is a convergent subsequence \revise{$\{z^{q_k}\}_{k\in\N}$} of \revise{$\{z^{p_k}\}_{k\in\N}$} such that $\lim_{k\to\infty} z^{q_k} = z^*$ for some $z^*$ satisfying $\|\Fnor(z^{*})\|=0$ and $x^*:=\proxl(z^*) \in \crit(\psi)$. Due to $\|x^k-x^*\| \leq \|z^k - z^*\|$, we then have $\lim_{k\to\infty} x^{q_k}=x^*$. Using the continuity of the Moreau envelope, this yields
\begin{align*} \psi^* & = {\lim}_{k \to \infty}\,\psi(x^{q_k})  = {\lim}_{k \to \infty}\,\vp(\proxi{\lambda}(z^{q_k})) + f(x^{q_k}) \\ & = {\lim}_{k\to\infty}\,\envi{\lambda}(z^{q_k}) - \tfrac{1}{2\lambda} \|z^{q_k}- x^{q_k}\|^2 + f(x^*) = \psi(x^*). \end{align*}
This readily implies $H(z^k) \to \psi^* = H(z^*)$ and finishes the proof.
\end{proof}

\bibliographystyle{siamplain}
\bibliography{references}

\end{document}

%% file: ex_shared.tex
\usepackage{lipsum}
\usepackage{amsfonts,amsmath,amsopn}
\usepackage{amssymb}
\usepackage{enumitem}
\usepackage{dsfont}
\usepackage{graphicx}
\usepackage{mathtools}
\usepackage{epstopdf}
\usepackage{algorithm}
\usepackage{algorithmic}
\usepackage{tikz}
\usepackage{pgfplots}
\pgfplotsset{compat=1.18}
\usepackage{comment}
\usepackage{bbding}
\usepackage{pifont}
\usepackage[flushleft]{threeparttable}
\usetikzlibrary{decorations.pathreplacing, positioning, arrows.meta,calc,math,shapes.geometric}
\usepackage{nicematrix}
\usepackage{booktabs}

\definecolor{lavender}{rgb}{0.9, 0.9, 0.98}
\usepackage{mdframed}

\ifpdf
  \DeclareGraphicsExtensions{.eps,.pdf,.png,.jpg}
\else
  \DeclareGraphicsExtensions{.eps}
\fi

\newcommand{\mer}{H}

\newcommand{\R}{\mathbb{R}}
\newcommand{\N}{\mathbb{N}}
\newcommand{\Rn}{\mathbb{R}^d}
\newcommand{\Rmn}{\mathbb{R}^{m\times n}}

\newcommand{\Rex}{(-\infty,\infty]}

\newcommand{\Prob}{\mathbb{P}}
\newcommand{\Exp}{\mathbb{E}}
\newcommand{\vp}{\varphi}

\newcommand{\dist}{\mathrm{dist}}
\newcommand{\crit}{\mathrm{crit}}

\newcommand{\Llambda}{{\sf L}_\lambda}

\newcommand{\sL}{{\sf L}}

\newcommand{\sT}{{\sf T}}

\newcommand{\cB}{\mathcal B}

\newcommand{\cF}{\mathcal F}

\newcommand{\cX}{\mathcal{X}}

\newcommand{\cM}{\mathcal{M}}

\newcommand{\sct}[1]{\boldsymbol{#1}}
\newcommand{\se}{\sct{e}}
\newcommand{\sd}{\sct{d}}
\newcommand{\sr}{\sct{r}}

\newcommand{\sk}{\sct{k}}

\newcommand{\sg}{\sct{g}}

\newcommand{\scp}{\sct{p}}

\newcommand{\scs}{\sct{s}}
\newcommand{\su}{\sct{u}}
\newcommand{\sw}{\sct{w}}
\newcommand{\sx}{\sct{x}}
\newcommand{\sy}{\sct{y}}
\newcommand{\sz}{\sct{z}}

\newcommand{\spsi}{\sct{\psi}}

\newcommand{\ssigma}{\sct{\varsigma}}

\newcommand{\revise}[1]{#1}

\newcommand{\fk}{\kappa}

\newcommand{\iprod}[2]{\langle #1, #2 \rangle}

\newcommand{\sxi}{\sct{\xi}}

\newcommand{\SGD}{\mathsf{SGD}}
\newcommand{\RR}{\mathsf{RR}}
\newcommand{\SKM}{\mathsf{SKM}}
\newcommand{\RDA}{\mathsf{RDA}}
\newcommand{\PSGD}{\mathsf{Prox}\text{-}\mathsf{SGD}}
\newcommand{\PSVRG}{\mathsf{Prox}\text{-}\mathsf{SVRG}}

\newcommand{\PSAGA}{(\mathsf{Prox}\text{-})\mathsf{SAGA}}
\newcommand{\STORM}{\mathsf{STORM}}

\newcommand{\PGD}{\mathsf{Prox}\text{-}\mathsf{GD}}
\newcommand{\NSGD}{\mathsf{Norm}\text{-}\mathsf{SGD}}

\newcommand{\prox}[1]{\mathrm{prox}_{#1}}

\newcommand\proxl{\mathrm{prox}_{\lambda\vp}}

\newcommand\proxi[1]{\mathrm{prox}_{#1\vp}}
\newcommand\envi[1]{\mathrm{env}_{#1\vp}}

\newcommand\Fnor{F^{\lambda}_{\mathrm{nor}}}

\newcommand\Fnat{F^{\lambda}_{\mathrm{nat}}}

\newcommand{\rmd}{\mathrm{d}}

\newcommand{\dom}[1]{\mathrm{dom\\}(#1)}

\newcommand\env[1]{\mathrm{env}_{#1}}

\newcommand{\norm}[1]{\|{#1}\|}
\newcommand{\prt}[1]{\left({#1}\right)}
\newcommand{\srt}[1]{\left[{#1}\right]}

\DeclareMathOperator*{\argmin}{argmin}

\definecolor{purp}{RGB}{152,24,147}
\definecolor{bluep}{RGB}{0,128,255}
\definecolor{redp}{RGB}{255,0,0}
\definecolor{orangep}{RGB}{255,128,0}

\definecolor{OxfordBlue}{rgb}{0,0.106,0.329}   
\definecolor{UMRed}{rgb}{0.73,0.09,0.19}   
\definecolor{CUBrown}{RGB}{152,95,42}   
\definecolor{LightBrown}{RGB}{219,199,181}   
\definecolor{MyBlue}{RGB}{123,144,210} 
\definecolor{LightBlue}{RGB}{229,232,247}   
\definecolor{MyGreen}{RGB}{165,222,228}   
\definecolor{LightGreen}{RGB}{228,245,247} 
\definecolor{turql}{RGB}{53,130,134}

\definecolor{MyGrayTwo}{RGB}{180,180,180} 
\definecolor{McBlack}{RGB}{39,37,31} 
\definecolor{McRed}{RGB}{218,41,28} 
\definecolor{McYellow}{RGB}{255,199,44} 
\definecolor{DeepBlue}{RGB}{2,103,178}
\definecolor{DeepRed}{RGB}{192,0,0}
\definecolor{MyOrange}{RGB}{237,156,89}
\definecolor{MyPurple}{RGB}{74,34,93}
\definecolor{MyTurq}{RGB}{89,157,160}

\definecolor{MyGray}{RGB}{180,180,180}   
\definecolor{orangep}{RGB}{255,128,0}

\definecolor{MyRed}{RGB}{192,0,0}

\newcommand{\be}{\begin{equation}}
	\newcommand{\ee}{\end{equation}}

\newsiamremark{remark}{Remark}
\newsiamremark{hypothesis}{Hypothesis}
\crefname{hypothesis}{Hypothesis}{Hypotheses}
\crefname{assumption}{Assumption}{Assumptions}

\newsiamthm{assumption}{Assumption}
\newsiamthm{claim}{Claim}

\headers{A Normal Map-Based Proximal SGD Method}{J. Qiu, L. Jiang, and A. Milzarek}

\title{A Normal Map-Based Proximal Stochastic Gradient Method: Convergence and Identification Properties\thanks{{\textbf{Funding:} Andre Milzarek was partly supported by the National Natural Science Foundation of China (Foreign Young Scholar Research Fund Project) under Grant No$.$ 12150410304, by the Shenzhen Science and Technology Program (No$.$ RCYX20221008093033010), and by the Guangdong Provincial Key Laboratory of Mathematical Foundations for Artificial Intelligence (2023B1212010001).}}}
\author{
Junwen Qiu\thanks{Industrial Systems Engineering and Management, National University of Singapore, 119077, Singapore (\email{jwqiu@nus.edu.sg}).}
\and Li Jiang\thanks{School of Data Science, The Chinese University of Hong Kong, Shenzhen (CUHK-Shenzhen), Guangdong, 518172, P.R. China (\email{lijiang@link.cuhk.edu.cn}, \email{andremilzarek@cuhk.edu.cn}).}
\and Andre Milzarek\footnotemark[3]
}
